\documentclass{article}
\usepackage[top=2cm, bottom=2cm, left=2cm, right=2cm]{geometry}
\usepackage{amssymb}
\usepackage{amsmath}
\usepackage{latexsym}
\usepackage{mathrsfs}%花体
\usepackage{graphicx}%插图
\usepackage[all]{xy}%交换图
\usepackage[CJKbookmarks=true]{hyperref}
\newtheorem{Th}{Theorem}
\newtheorem{Prop}[Th]{Proposition}

\title{Computations of the Comodule Structures \\of the Chow rings of Flag Varieties}
\author{XIONG Rui}

\begin{document}

\def\CH{\operatorname{CH}}
\def\longdash{\mathord-\!\!\!\mathord-\!\!\!\mathord-\!\!\!\mathord-}
\def\longDash{\mathord=\!\!\!\mathord\Rightarrow\!\!\!\mathord=\!\!\!\mathord=}
\def\longDdash{\mathord\equiv\!\!\!\mathord\Rrightarrow\!\!\!\mathord\equiv\!\!\!\mathord\equiv}

\def\ru{\hline\rule[-.5pc]{0pc}{1.5pc}}
\def\Ad{\operatorname{Ad}}

\maketitle

\begin{abstract}
  Let $G$ be a connected reductive group, and $G/B$ be its flag variety.
  Let $\pi:G\to G/B$ be the natural projection.
  In this paper, we developed an algorithm to describe the map
  $\pi^* :\CH^*(G/B;\mathbb{F}_p)\longrightarrow \CH^*(G;\mathbb{F}_p)$ in terms of Schubert cells.
  Taking advantage of the Pieri rule,
  we give an explicit formula for $A$-type, $C$-type, $G_2$, $F_4$ of the cohomology map
  $\pi^* :\CH^*(G/B;\mathbb{F}_p)\longrightarrow \CH^*(G;\mathbb{F}_p)$, and some partial result of $\pi^*$ is given for $E_6$ and $E_7$.
  Denote the group action map $\mu:G\times G/B\to G/B$, we also give an explicit formula for $A$-type, $C$-type, $G_2$, $F_4$ of the cohomology map
  $\mu^*: \CH^*(G/B;\mathbb{F}_p)\longrightarrow \CH^*(G\times G/B;\mathbb{F}_p)$.
\end{abstract}

I would politely express my gratitude to Victor Petrov who pointed out a significant optimization of the algotirhm.
Thanks also Luyu Chen, Shoumin Liu and Arsenty Kambalin for the discussion.

%\tableofcontents

\section{Introduction}
%
%Let $p$ be a prime or $0$.
%Denote $\mathbb{F}_p$ the finite fields of $p$ elements when $p$ is finite.
%When $p=0$, the $\mathbb{F}_p$ is understood as $\mathbb{Q}$.
%
%Let $G$ be a connected reductive group over $\mathbb{C}$, and $B$ be its Borel subgroup.
%The homogenous variety $G/B$ is called the \emph{flag manifold} of $G$. Let $W$ be its Weyl group.
%We will denote the natural projection by $\pi:G\to G/B$.
%
%In this paper, we will give an explicit formula of the map
%$$\pi^* H^*(G/B;\mathbb{F}_p)\longrightarrow H^*(G;\mathbb{F}_p)$$
%in terms of Schubert cells $\{[\Sigma_w]:w\in W\}$.
%
%Let $\mu: G\times G/B\to G/B$ be the map of left multiplication.
%In author's previous preprint \cite{xiong2020comodule}, the comodule structure is computed
%$$[\Sigma_w]\longmapsto \sum_{{w=u\cdot v}\atop{\ell(w)=\ell(u)+\ell(v)}} \pi^*[\Sigma_u]\otimes [\Sigma_v].$$
%So the description of $\pi^*$ will give an explicit formula for this map.

\paragraph{Lie Groups.}
Let $G$ be a connected reductive group over $\mathbb{C}$, and $B$ be its Borel subgroup.
The homogenous variety $G/B$ is called the \emph{flag manifold} of $G$.
We will denote the natural projection by $\pi:G\to G/B$.

Let $K$ be its compact form of $G$, and $T$ be its maximal torus.
The natural inclusion $K\subseteq G$ is a homotopy equivalence,
  and induces a homeomorphism $K/T\cong G/B$ which gives the complex structure over $K/T$.

Let $Z(G)$ be the center of $G$, and the quotient group $G/Z(G)$ is semisimple and has the same flag manifold as $G$'s.
So $\pi$ factors through $G/Z(G)$.

%\paragraph{Weyl groups.}
Let $W=N_K(T)/T$ be the Weyl group of $G$.
Let $\ell:W\to \mathbb{Z}$ be the length function.
We denote the longest word by $w_0$ in $W$.
We write $w=u\odot v$ for three elements $w,u,v\in W$, if $w=uv$ and $\ell(w)=\ell(u)+\ell(v)$.
Similar notation for $w=u\odot v\odot z$, etc.

For any element $w\in W$, we fix a lift in $G$, denoted also by $w$.
We have the \emph{Bruhat decomposition} $G=\bigcup_{w\in W} BwB$.
Let $BwB/B$ be the image of $BwB$ in $G/B$.
It is known that $BwB/B$ is isomorphic to affine space $\mathbb{C}^{\ell(w)}$.
So there is a cellular structure decomposition $G/B$ by $\{BwB/B:w\in W\}$.

Denote $B^-=w_0Bw_0$ the opposite Borel group. We denote $\Sigma_w$ be the Zariski closure of $B^-wB/B$ in $G/B$.
Then $\Sigma_w$ is of complex codimension $\ell(w)$.

\paragraph{Chow ring and cohomology ring.}
Let $p$ be a prime or $0$.
Denote $\mathbb{F}_p$ the finite fields of $p$ elements when $p$ is finite.
When $p=0$, the $\mathbb{F}_p$ is understood as $\mathbb{Q}$.

Let $\CH^*(X;\mathbb{F}_p)$ be the cohomology group with coefficients in $\mathbb{F}_p$.
To work better with the cohomology ring, we \emph{double its degree}, so $\CH^*$ has no nonzero odd degree element.
When $X$ is smooth, there is a commutative ring structure.
When $X$ is a variety, and $Y$ be a closed subvariety,
we will denote $[Y]$ the fundamental class in $\CH^*(X;\mathbb{F}_p)$,
whose degree is $2$ times its codimension.

%We usually omit the symbol of multiplication.
Let $H^*(X;\mathbb{F}_p)$ be the cohomology group with coefficients in $\mathbb{F}_p$.
We will denote the cup product by $\smile$.
When $X$ is a smooth variety, and $Y$ be a closed subvariety,
we will denote $[Y]$ the fundamental class in $H^*(X;\mathbb{F}_p)$, whose degree is $2$ times its codimension.
In this case, there is a natural map $\CH^*(X;\mathbb{F}_p)\to H^{*}(X;\mathbb{F}_p)$.

%\paragraph{Cohomology of flag manifolds.}
Consider the maps induced by $\pi$ above
$$\xymatrix{
\CH^*(G/B;\mathbb{F}_p)\ar[r]^{\pi_*}\ar[d]_{\alpha}&\CH^*(G;\mathbb{F}_p)\ar[d]^{\beta}\\
H^*(G/B;\mathbb{F}_p)\ar[r]_{\pi_*} & H^*(G;\mathbb{F}_p).} $$

Due to the cellular structure over $G/B$, the map $\alpha$ is an isomorphism,
and $H^*(G/B;\mathbb{F}_p)$ is freely generated by $\{[\Sigma_w]: w\in W\}$ as an abelian group,
where $[\Sigma_w]$ is of degree $2\ell(w)$.

It follows from Grothendieck \cite{grothendieck1958torsion} that the above $\pi_*$ is a surjection, and $\beta$ is an injection.
On the other hand, the Chow ring, and the cohomology ring of $G$ is computed in the paper \cite{kac1985torsion} (see the table below).
Let $K/T\to BT$ be the classifying map of $T$-principle bundle $K\to K/T$.
It turns out the kernel of $\pi^*$ is the ideal generated by the image of the induced map $H^2(BT;\mathbb{F}_p)\to H^2(G/B;\mathbb{F}_p)$
(see \cite{kac1985torsion}).
In particular, when $p=0$, $\CH^{*}(G;\mathbb{Q})=\CH^{0}(G;\mathbb{Q})=\mathbb{Q}$.

\paragraph{Comodule structure.}
Let $\mu: G\times G/B\to G/B$ be the map of the left action.
Since $\mu$ is algebraic,
$\mu^*:H^*(G/B;\mathbb{F}_p)\otimes H^*(G/B;\mathbb{F}_p)\otimes H^*(G;\mathbb{F}_p)$ factors through
$\CH^*(G/B;\mathbb{F}_p)\otimes \CH^*(G;\mathbb{F}_p)$.
So it suffices to compute $\CH(G/B;\mathbb{F}_p)\to \CH(G;\mathbb{F}_p)\otimes \CH(G/B;\mathbb{F}_p)$.

In the author's previous preprint \cite{xiong2020comodule}, the comodule structure is computed
$$[\Sigma_w]\longmapsto \sum_{w=u\odot v} \pi^*[\Sigma_u]\otimes [\Sigma_v].$$
So the description of $\pi^*$ will give an explicit formula for this map.
%and as a corollary, the coalgebra structure $\CH^*(G;\mathbb{F}_p)\to \CH^*(\mathbb{G};\mathbb{F}_p)\otimes \CH^*(\mathbb{G};\mathbb{F}_p)$

\paragraph{The layout of this paper.}
The main purpose of this paper is to compute the projection
$$\pi^*: \CH^*(G/B;\mathbb{F}_p)\longrightarrow \CH^*(G;\mathbb{F}_p)$$
and the comodule structure
$$\mu^*:\CH(G/B;\mathbb{F}_p)\to \CH(G;\mathbb{F}_p)\otimes \CH(G/B;\mathbb{F}_p).$$
Since $\pi$ factors through $G/Z(G)$, so we can assume $G$ to be simple and of adjoint type.

%When $p=0$, the above map is trivial, this will be discussed in section \ref{Trivialmod0}.
We will develop an algorithm to compute $\pi^*$ in section \ref{AlgoPieri}.
%We will prove a couple of general facts to be used in section \ref{GeneralFacts}.
For $G$ is of $A$-type or $C$-type, there is a conceptual description, and we will discuss it in section \ref{DescForAC}.
For the exceptional group $G_2,F_4,E_6,E_7$, they will be discussed separately in the sections ahead.
%Some of them are not completely described.

The result of this paper is summarized in the following table.
%\begin{table}[c]
$$\begin{array}{c|c|c|c} \ru
G & p& \CH(G;\mathbb{F}_p) \text{ from table 2 of }\cite{kac1985torsion} & \text{result} \\\ru
\operatorname{PGL}_n & p|n; p^k\|n & \mathbb{F}_p[x_{2}]\big/\big<x_2^{p^k}\big> & \text{Theorem }\ref{Amodp} \\\ru
%\operatorname{Spin}_{2n+1} & p=2  &\text{see \cite{kac1985torsion}}& \text{not yet}\\\ru
%\Ad \operatorname{SO}_{2n+1}
\operatorname{Spin}_{n} & p=2 &\text{see }\cite{kac1985torsion}& \text{not yet}\\\ru
\Ad\operatorname{SO}_{n} & p=2 &\text{see }\cite{kac1985torsion}& \text{not yet}\\\ru
\Ad \operatorname{Sp}_n & p=2; 2^k\|n & \mathbb{F}_2[x_{2}]\big/\big<x_2^{2^k}\big> & \text{Theorem }\ref{Cmodp}\\\ru
%\operatorname{Spin}_{2n} & p=2 &\text{see \cite{kac1985torsion}}& \text{not yet}\\\ru
%\Ad \operatorname{SO}_{2n} & p=2 &\text{see \cite{kac1985torsion}}& \text{not yet}\\\ru
G_2 & p=2 & \mathbb{F}_2[x_6]\big/\big<x_6^2\big> &\text{Theorem }\ref{G2mod2} \\\ru
F_4 & p=2 & \mathbb{F}_2[x_6]\big/\big<x_6^2\big> &\text{Theorem }\ref{F4mod2}\\\ru
    & p=3 & \mathbb{F}_3[x_8]\big/\big<x_8^3\big> &\text{Theorem }\ref{F4mod3}\\\ru
E_6 & p=2 & \mathbb{F}_2[x_6]\big/\big<x_6^2\big> & \text{Theorem }\ref{E6mod2}\\\ru
    & p=3 & \mathbb{F}_3[x_8]\big/\big<x_8^3\big> & \text{Theorem }\ref{scE6mod3}\\\ru
\Ad E_6 & p=2 & \mathbb{F}_2[x_6]\big/\big<x_6^2\big>& %\text{the same as $E_6$}
                                                        \text{Theorem }\ref{E6mod2}\\\ru
    & p=3 & \mathbb{F}_3[x_2,x_8]\big/\big<x_2^9,x_8^3\big>& \text{not yet}\\\ru %\ref{E6mod3}
E_7 & p=2 & \mathbb{F}_2[x_6,x_{10},x_{18}]\big/\big<x_6^2,x_{10}^2,x_{18}^2\big> & \text{not yet}\\\ru
    & p=3 & \mathbb{F}_3[x_8]\big/\big<x_8^3\big> & \text{Theorem }\ref{E7mod3}\\\ru
\Ad E_7 & p=2 & \mathbb{F}_2[x_2,x_6,x_{10},x_{18}]\big/\big<x_2^2,x_6^2,x_{10}^2,x_{18}^2\big>& \text{not yet}\\\ru %\ref{E7mod2}
    & p=3 & \mathbb{F}_3[x_8]\big/\big<x_8^3\big> & %\text{the same as $E_7$}
                                                    \text{Theorem }\ref{E7mod3}\\\ru
E_8 & p=2 & \mathbb{F}_2[x_6,x_{10},x_{18},x_{30}]\big/\big<x_6^8,x_{10}^4,x_{18}^2,x_{30}^2\big>&\text{no yet} \\\ru
    & p=3 & \mathbb{F}_3[x_8,x_{20}]\big/\big<x_8^3,x_{20}^3\big> & \text{no yet}\\\ru
    & p=5 & \mathbb{F}_3[x_{12}]\big/\big<x_{12}^5\big> &\text{no yet}
\\\hline
\end{array}$$
%\caption{Summary of the results}\label{SummofRes}
%\end{table}
%For the case
%%The missing case is $B$-type, $D$-type and $E_8$; . %takes too much time to compute.
%$$\begin{array}{c|c|c|c} \ru
%G & p& \CH(G;\mathbb{F}_p)&  \\\ru
%\Ad E_6 & p=2 & \mathbb{F}_2[x_6]\big/\big<x_6^2\big>& \text{the same as $E_6$}\\\ru
%    & p=3 & \mathbb{F}_3[x_2,x_8]\big/\big<x_2^9,x_8^3\big>& \text{not yet}\\\ru %\ref{E6mod3}
%\Ad E_7 & p=2 & \mathbb{F}_2[x_2,x_6,x_{10},x_{18}]\big/\big<x_2^2,x_6^2,x_{10}^2,x_{18}^2\big>& \text{not yet}\\\ru %\ref{E7mod2}
%    & p=3 & \mathbb{F}_3[x_8]\big/\big<x_8^3\big> & \text{the same as $E_7$}\\\hline
%\end{array}$$

%\section{General Facts}\label{GeneralFacts}

\section{The Algorithm}\label{AlgoPieri}

\paragraph{Pieri rule.}
%\paragraph{Root system.}
Let $\{\alpha_i:i\in I\}$ be the set of simple roots,
$\{\omega_i:i\in I\}$ the corresponding fundamental weights.
Let $\{s_i: i\in I\}$ be the set of simple reflections.
Assume $K$ is a semisimple compact group.
We can identify $H^2(BT;\mathbb{Z})$ with the character of $T$.

When $K$ is simply-connected, it can be further recognized as the weight lattice,
and $[\Sigma_{s_i}]\in H^*(G/B;\mathbb{Z})$ is presented by the image of $\omega_i\in H^*(BT:\mathbb{Z})$.
In particular, $H^2(BT;\mathbb{Z})\to H^2(G/B;\mathbb{Z})$ is an isomorphism.
For any $p$, the kernel of $\pi^*$ is the ideal generated by $H^2(G/B;\mathbb{F}_p)$.
On the contrary, when $K$ is of adjoint type, it can be further recognized as the root lattice.

%Let $\{\partial_i: i \in I\}$ be the set of Demazure operators over $H^*(G/B;\mathbb{Z})$.
%It is classic that $\partial_i[\Sigma_w]=\begin{cases}
%{}[\Sigma_{ws_i}], & \ell(ws_i)=\ell(w)-1,\\
%0, & \text{otherwise},
%\end{cases}$.
%If $w=s_{i_1}\odot \cdots \odot s_{i_k}$ for $i_\bullet\in I$,
%we can define the \emph{Demazure operator}
%$\partial_w=\partial_{i_1}\cdots \partial_{i_k}$.
%It turns out it is independent of the choice of the reduced word.
%Then for any $\alpha\in H^{2k}(G/B;\mathbb{Z})$,
%$$\alpha=\sum_{\ell(w)=k}\partial_w\alpha \cdot [\Sigma_w], $$
%where $\partial_w \alpha\in H^0(G/B)\cong \mathbb{Z}$.

%The Demazure operators satisfy $\partial_i(\alpha\beta)=\partial_i\alpha\smile \beta + s_i \alpha\smile \partial_i \beta$,
%where $s_i$ stands the map induced by the conjugation action of $W$ on $K/T$.
%As a result, if $w=s_{i_1}\odot \cdots \odot s_{i_k}$ for $i_\bullet\in I$,
%$\partial_w$ acts on product as
%$\big(\partial_{i_1}\otimes 1+s_{i_1}\otimes \partial_{i_1}\big)\cdots \big(\partial_{i_k}\otimes 1+s_{i_k}\otimes \partial_{i_k}\big)$.
For any element $\alpha \in H^2(G/B;\mathbb{Z})$ and $u\in W$, we have the Pieri formula \cite{hiller1982geometry}
$$\alpha\smile [\Sigma_u]=\sum_{} \left<\beta^\vee,\alpha\right>[\Sigma_w]$$
where the sum is taken over all $w$ such that $\ell(w)=\ell(u)+1$ and there is a positive root $\beta$ with $w=us_\beta$,
where $\beta^\vee=\beta/\left<\beta,\beta\right>$ for the standard inner product $\left<\cdot,\cdot\right>$ over the root space.

\paragraph{The algorithm. }The author used \textsf{Sage} to realize the algorithm.
\begin{verbatim}
Ct=CartanType("G2"); DynkinDiagram(Ct)      # Input

Phi=RootSystem(Ct); hvee=Phi.ambient_space();
Phi_plus=hvee.positive_roots(); beta=Phi_plus;
Delta=hvee.simple_roots(); alpha=Delta; omega=hvee.fundamental_weights();
I=Delta.keys();
W=WeylGroup(Ct,prefix="s"); s=W.simple_reflections();
alphavee=hvee.simple_coroots();

def B(x,y): return x.inner_product(y)
def vee(x): return 2*x/B(x,x)

refl=W.reflections(); Refl=[refl[beta[k]] for k in [0..len(beta)-1]];

p=2;                                        # Input
k=6; k=k/2;                                 # Input

Lengthk=list(W.elements_of_length(k));
Lengthkk=list(W.elements_of_length(k-1));
print(len(Lengthkk),len(Lengthk))

V=VectorSpace(GF(p),len(Lengthk))

Lengthkinv={}
for b in [0..len(Lengthk)-1]:
    Lengthkinv[Lengthk[b]]=b

Relation = {}
for a in [0..len(Lengthkk)-1]:
    for i in I:
        Relation[(a,i)]=V(0);

for a in [0..len(Lengthkk)-1]:
    for j in [0..len(beta)-1]:
        z=Lengthkk[a]* Refl[j];
        if z in Lengthk:
            for i in I:
                Relation[(a,i)][Lengthkinv[z]]=B(vee(beta[j]),alpha[i]);    #(*)

Relation=Relation.values()

Q,pi,lift=V.quotient_abstract(V.span(Relation))

dim(Q)

Image = [0 for b in [0..len(Lengthk)-1]]
for b in [0..len(Lengthk)-1]:
    v=V(0); v[b]=1;
    Image[b]=pi(v)

for q in Q:
    if q != 0:
        print(q,"=")
        for b in [0..len(Lengthk)-1]:
            if Image[b]==q:
                print("$",latex(Lengthk[b]),"$,")
\end{verbatim}

%
%Then \verb|Ct| is the type with the index set \verb|I|;
%%$\mathtt{hvee}=\mathfrak{h}^\vee$ is the dual of the Cartan algebra;
%\verb|Phi| is the root system; \verb|Phi_plus|=[\verb|beta[j]|] is the list of positive roots;
%[\verb|omega[i]|: \verb|i| $\in$ \verb|I|] the list of fundamental weights.
%\verb|Delta|=[\verb|alpha[i]|:\verb|i|$\in$ \verb|I|] is the list of simple roots;
%[\verb|alphavee[i]|:\verb|i|$\in$ \verb|I|] the list of simple coroots, with
%\verb|B|(,) the pairing.
%
%$\mathtt{W}=W$ is the Weyl group with $\{\mathtt{s}[i]:i\in \mathtt{I}\}=\{s_i:i\in I\}$ the set of simple reflections.
%
%For any root $\beta$, the coroot $\mathtt{vee}(\beta)=\beta^\vee$.

\noindent\textbf{Remark.}
If one only wants the result for the simply-connected case,
it suffices to exchange the above \verb|#(*)| to
\begin{verbatim}
                Relation[(a,i)][Lengthkinv[z]]=B(vee(beta[j]),omega[i]);    #(*)
\end{verbatim}

\section{Computations of $A$-type and $C$-type}\label{DescForAC}

\paragraph{A general fact.}
Let $X,Y,Z$ be three path connected spaces.
Pick an $x_0\in X$ and a $y_0\in Y$ arbitrarily.
Let $f: X\times Y\to Z$ be a continuous map.
%It induces map $H^*(Z;\mathbb{F}_p)\to H^*(X;\mathbb{F}_p)\otimes H^*(Y;\mathbb{F}_p)$.
If $H^1(Y;\mathbb{F}_p)=0$, then the map
$$H^2(Z)\to H^2(X)\otimes H^0(Y)\oplus H^0(X)\otimes H^2(Y)$$
is given by $x\mapsto f_1^*(x)\otimes 1+1\otimes f_2^*(x)$ where
$f_1:X\to Z$ maps $x$ to $f(x,y_0)$ and
$f_2:Y\to Z$ maps $y$ to $f(x_0,y)$.

\paragraph{The $A$-type case.}
Let $G=\operatorname{PGL}_n$.
In type A, it is well-known that
$$H^*(G/B;\mathbb{Z})=\frac{\mathbb{Z}[x_1,\ldots,x_n]}{\left<f(x_1,\ldots,x_n)\in \mathbb{Z}[x_1,\ldots,x_n]^W_{+}\right>}$$
where $\mathbb{Z}[x]^W_{+}$ stands the set of symmetric polynomials without constant term.
In particular, $H^2(G/B;\mathbb{Z})$ is generated by $H^2$-elements.

Note that the root system is given by $\{x_i-x_j:i\neq j\}$, so
$$\CH^*(G;\mathbb{F}_p)
=\frac{\CH^*(G/B;\mathbb{F}_p)}{\left<x_i-x_j\right>}
\cong\frac{\mathbb{F}_p[t]}{\left<f(t,\ldots,t)\in \mathbb{Z}[x_1,\ldots,x_n]^W_{+}\right>}. $$
Actually, $\mathbb{Z}[x_1,\ldots,x_n]^W$ is generated by the elementary symmetric polynomials.
Therefore, $\CH^*(G;\mathbb{F}_p)=\mathbb{F}_p[t]/\left<\binom{n}{i}t^i: i=1,\ldots,n\right>$.
Assume $p^k|n$ but $p^{k+1}\nmid n$.
Elementary number theory shows that $\binom{n}{i}=0$ for $i<p^k$,
and $\binom{n}{p^k}\neq 0$ in $\mathbb{F}_p$. As a result,
$$\CH(G;\mathbb{F}_p)=\mathbb{F}_p[t]/\big<t^{p^k}\big>\qquad \text{$p^k|n$ but $p^{k+1}\nmid n$}$$
with $\deg t=2$. This is also summarized in \cite{kac1985torsion}.
%
%Denote $\tau: H^*(G/B;\mathbb{F}_p)\to H^*(G/B;\mathbb{F}_p)\otimes \mathbb{F}_p[t]/\big<t^{p^k}\big>$
%given by the shift by $t$, $f(x_1,\ldots,x_n)\mapsto f(x_1+t,\ldots,x_n+t)$.
%This is well-defined, since for $f\in \mathbb{Z}[x_i]^W_+$,
%$f(x_i+t)$ has no constant term by above argument, and other terms are symmetric.

\begin{Th}[$A$-type]\label{Amodp}For the group $G=\operatorname{PGL}_n$, we have
\begin{description}
  \item[(1)] The projection $\CH^*(G/B;\mathbb{F}_p)\to \CH^*(G;\mathbb{F}_p)$ is given by
$$x_i \longmapsto t\qquad \forall i\in \{1,\ldots,n\}.$$
  \item[(2)] The comoudle structure $\CH(G/B;\mathbb{F}_p)\to \CH(G\times G/B;\mathbb{F}_p)$ is given by
$$x_i\longmapsto x_i+t\qquad \forall i \in \{1,\ldots,n\}.$$
\end{description}
\end{Th}

%To prove this, note that $H^*(G/B;\mathbb{F}_p)$ is generated by $H^2(G/B;\mathbb{F}_p)$,
%and $x_i\mapsto x_i+t$ for all $i$.

\paragraph{The $C$-type case.}
Let $G=\Ad\operatorname{Sp}_n$.
In type C, it is also known that
$$H^*(G/B;\mathbb{Z})=\frac{\mathbb{Z}[x_1,\ldots,x_n]}{\left<f(x_1,\ldots,x_n)\in \mathbb{Z}[x_1,\ldots,x_n]^W_{+}\right>}$$
where $\mathbb{Z}[x]^W_{+}$ the set of even symmetric polynomials without constant term and actually generated by
the elementary symmetric polynomials in $x_i^2$.
In particular, $H^2(G/B;\mathbb{Z})$ is generated by $H^2$-degree.

Note that the root system is given by $\{x_i-x_j:i\neq j\}\cup \{\pm 2 x_i:i=1,\ldots,b\}$.
The analysis is the same as the $A$-type case. We can conclude that
$$\CH^*(G;\mathbb{F}_2)=\mathbb{F}_2[t]/\big<t^{2^k}\big>\qquad \text{$2^k|n$ but $2^{k+1}\nmid n$}, $$
with $\deg t=2$.

\begin{Th}[$C$-type]\label{Cmodp}For the group $G=\Ad\operatorname{Sp}_n$, we have
\begin{description}
  \item[(1)] The projection $\CH^*(G/B;\mathbb{F}_2)\to \CH^*(G;\mathbb{F}_2)$ is given by
$$x_i \longmapsto t\qquad \forall i\in \{1,\ldots,n\}.$$
  \item[(2)] The comoudle structure $\CH(G/B;\mathbb{F}_2)\to \CH(G\times G/B;\mathbb{F}_2)$ is given by
$$x_i\longmapsto x_i+t\qquad \forall i \in \{1,\ldots,n\}.$$
\end{description}
\end{Th}

\section{Computations of $G_2$}

Denote $G_2$ the complex Lie group. Denote $G/B$ the flag variety of it, and $W$ the Weyl group of it.
Now $\CH^{*}(G/B;\mathbb{F}_2)$ is of order $12$.
Denote the short (resp. long) simple root by $\alpha$ (resp. $\beta$),
and $s$ (resp. $t$) the corresponding simple reflection.
$$\mathop{\circ}\limits^{\beta}_{t}\longDdash\mathop{\circ}\limits^{\alpha}_{s}$$
The list of elements of Weyl group is given by
\def\bi#1#2{\displaystyle{{#1}\atop{#2}}}
$$\begin{array}{c|c|c|c|c|c|c|c}\hline\rule{0pc}{1.5pc}
\textrm{element} & 1 & \bi{s}{t}& \bi{st}{ts} & \bi{sts}{tst} & \bi{stst}{tsts} & \bi{ststs}{tstst} & \bi{ststst}{=tststs}\\ [2ex]\ru
\textrm{length} & 0 & 1 & 2 & 3 & 4 & 5 & 6 \\\hline
\end{array}$$

\paragraph{The case $p=2$. }
Now $\CH^*(G_2;\mathbb{F}_2)=\mathbb{F}_2[x_6]/\left<x_6^2\right>$.

There are only two elements of Weyl group of length $3$, say $sts$ and $tst$,
and only two for length $2$, say $st$ and $ts$. So there are four relations
$$\begin{cases}
\alpha\smile [\Sigma_{st}]= \left<\alpha^\vee,\alpha\right>[\Sigma_{sts}]
+\left<(3\alpha+2\beta)^\vee,\alpha\right>[\Sigma_{tst}]
&=2[\Sigma_{sts}]+0 [\Sigma_{tst}]\\
\beta\smile [\Sigma_{st}]= \left<\alpha^\vee,\beta\right>[\Sigma_{sts}]
+\left<(3\alpha+2\beta)^\vee,\beta\right>[\Sigma_{tst}]
&=-[\Sigma_{sts}]+[\Sigma_{tst}]\\
\alpha\smile [\Sigma_{ts}]= \left<(2\alpha+\beta)^\vee,\alpha\right>[\Sigma_{sts}]
+\left<\beta^\vee,\alpha\right>[\Sigma_{tst}]
&=[\Sigma_{sts}]-[\Sigma_{tst}]\\
\beta\smile [\Sigma_{ts}]= \left<(2\alpha+\beta)^\vee,\beta\right>[\Sigma_{sts}]
+\left<\beta^\vee,\beta\right>[\Sigma_{tst}]
&=0[\Sigma_{sts}]+2[\Sigma_{tst}]\\
\end{cases}$$
So $\pi^*[\Sigma_{tst}]=\pi^*[\Sigma_{sts}]=x_6$.

\begin{Th}[The case $G_2$]\label{G2mod2}For the Chow ring of the exceptional group $G_2$ over $\mathbb{F}_2$, we have
\begin{description}
\item[(1)] The projection $\CH^*(G/B;\mathbb{F}_2)\to \CH^*(G_2;\mathbb{F}_2)$ is given by
$$[\Sigma_w]\longmapsto \begin{cases}
1, & w=1\\
x_6, & w=sts,tst,\\
0, & \textrm{otherwise}.
\end{cases}$$
\item[(2)] The comodule structure $\CH^*(G/B;\mathbb{F}_2)\to \CH^*(G_2;\mathbb{F}_2)\otimes \CH^*(G/B;\mathbb{F}_2)$ is given by
$$[\Sigma_{w}]\longmapsto 1\otimes [\Sigma_w]+\begin{cases}
0, & w=1,s,t,st,ts,\\
x_6\otimes [\Sigma_{stsw}],& w=sts,stst,ststs,\\
x_6\otimes [\Sigma_{tstw}],& w=tst,tsts,tstst,\\
x_6\otimes ([\Sigma_{tst}]+[\Sigma_{sts}]),& w=ststst.\\
\end{cases}
$$
%(3) The coalgebra structure $\CH^*(G_2;\mathbb{F}_2)\to \CH^*(G_2;\mathbb{F}_2)\otimes \CH^*(G_2;\mathbb{F}_2)$ is given by
%$$x_6\longmapsto x_6\otimes 1+1\otimes x_6. $$
\end{description}
\end{Th}

\noindent\textbf{Remark.} It seems that $G_2$ is the only case which can be computed by hand.

%Note that the last assertion follows easily from the fact that $H^*(G_2;\mathbb{F}_2)$ is monogenic.

\section{Computations of $F_4$}

The Dynkin diagram of $F_4$ is labeled as follows.
$$\underset{1}{\circ}\longdash \underset{2}{\circ}\longDash \underset{3}{\circ}\longdash \underset{4}{\circ}$$
We denote $s_i$ the corresponding simple reflection with $i\in \{1,2,3,4\}$.

Denote $F_4$ the complex Lie group. Denote $G/B$ the flag variety of it, and $W$ the Weyl group of $F_4$.

\paragraph{The case $p=2$.}
Now $\CH^*(F_4; \mathbb{F}_2)=\mathbb{F}_2[x_6]\big/\big<x_6^2\big>$.
Let
$$\Pi(x_6)=\{s_{1}s_{2}s_{3},
s_{2}s_{3}s_{1},
s_{2}s_{3}s_{2},
s_{3}s_{1}s_{2},
s_{3}s_{2}s_{1},
s_{3}s_{2}s_{3}\}\subseteq W. $$

\begin{Th}[The case $F_4$, $p=2$]\label{F4mod2}For the Chow ring of the exceptional group $F_4$ over $\mathbb{F}_2$,
we have
\begin{description}
\item[(1)] The projection $\pi^*:\CH^*(G/B;\mathbb{F}_2)\to \CH^*(F_4;\mathbb{F}_2)$ is given by
$$[\Sigma_w]\longmapsto \begin{cases}
1, & w=1\\
x_6,  & w\in \Pi(x_6),\\
0, & \text{otherwise}.
\end{cases}$$
\item[(2)] The comodule structure $\CH^*(G/B;\mathbb{F}_2)\to \CH^*(F_4;\mathbb{F}_2)\otimes \CH^*(G/B;\mathbb{F}_2)$ is given by
$$[\Sigma_w]\longmapsto 1\otimes [\Sigma_w]+
x_6\otimes \sum_{{w=u\odot v,} \atop {u\in \Pi(x_6)}} [\Sigma_{v}]. $$
%(3) The coalgebra structure $\CH^*(F_4;\mathbb{F}_2)\to \CH^*(F_4;\mathbb{F}_2)\otimes \CH^*(F_4;\mathbb{F}_2)$ is given by
%$$x_6\longmapsto x_6\otimes 1+1\otimes x_6. $$
\end{description}
\end{Th}

\noindent\textbf{Remark.} Not like the case $G_2$, there are $16$ many elements in Weyl group $W$ of length $3$
but only $6$ of them are mapped to $x_6$ by $\pi^*$.

\paragraph{The case $p=3$.}
Now
$\CH^*(F_4; \mathbb{F}_3)=\mathbb{F}_3[x_8]\big/\big<x_8^3\big>$.
Let
$$\Pi(-x_8)=\left\{\begin{array}{c}
s_{1}s_{2}s_{3}s_{2},
s_{1}s_{2}s_{3}s_{4},
s_{2}s_{3}s_{4}s_{2},
s_{3}s_{1}s_{2}s_{3},\\
s_{3}s_{4}s_{2}s_{1},
s_{4}s_{2}s_{3}s_{1},
s_{4}s_{3}s_{1}s_{2},
s_{4}s_{3}s_{2}s_{3}
\end{array}\right\}. $$
Denote also
$$\Pi(x_8)=\left\{\begin{array}{c}
s_{2}s_{3}s_{2}s_{1} ,
s_{2}s_{3}s_{4}s_{1} ,
s_{3}s_{2}s_{3}s_{1} ,
s_{3}s_{2}s_{3}s_{4} ,\\
s_{3}s_{4}s_{1}s_{2} ,
s_{4}s_{1}s_{2}s_{3} ,
s_{4}s_{2}s_{3}s_{2} ,
s_{4}s_{3}s_{2}s_{1}
\end{array}\right\}. $$
Note that $\Pi(-x_8)=\{w^{-1}\in W: w\in \Pi(x_8)\}$.

Let $\Pi(x_8^2)$ be the subset of $W$ consisting of
\begin{quote}
$ s_{1}s_{2}s_{3}s_{4}s_{1}s_{2}s_{3}s_{1} $,
$ s_{1}s_{2}s_{3}s_{4}s_{3}s_{1}s_{2}s_{3} $,
$ s_{1}s_{2}s_{3}s_{4}s_{3}s_{2}s_{3}s_{2} $,
$ s_{2}s_{3}s_{1}s_{2}s_{3}s_{1}s_{2}s_{1} $,
$ s_{2}s_{3}s_{1}s_{2}s_{3}s_{4}s_{1}s_{2} $,
$ s_{2}s_{3}s_{1}s_{2}s_{3}s_{4}s_{3}s_{1} $,
$ s_{2}s_{3}s_{1}s_{2}s_{3}s_{4}s_{3}s_{2} $,
$ s_{2}s_{3}s_{4}s_{1}s_{2}s_{3}s_{1}s_{2} $,
$ s_{2}s_{3}s_{4}s_{3}s_{1}s_{2}s_{3}s_{2} $,
$ s_{2}s_{3}s_{4}s_{3}s_{2}s_{3}s_{1}s_{2} $,
$ s_{3}s_{1}s_{2}s_{3}s_{4}s_{1}s_{2}s_{1} $,
$ s_{3}s_{1}s_{2}s_{3}s_{4}s_{1}s_{2}s_{3} $,
$ s_{3}s_{2}s_{3}s_{1}s_{2}s_{3}s_{2}s_{1} $,
$ s_{3}s_{2}s_{3}s_{1}s_{2}s_{3}s_{4}s_{1} $,
$ s_{3}s_{2}s_{3}s_{1}s_{2}s_{3}s_{4}s_{3} $,
$ s_{3}s_{2}s_{3}s_{4}s_{1}s_{2}s_{3}s_{1} $,
$ s_{3}s_{2}s_{3}s_{4}s_{2}s_{3}s_{1}s_{2} $,
$ s_{3}s_{2}s_{3}s_{4}s_{2}s_{3}s_{2}s_{1} $,
$ s_{3}s_{2}s_{3}s_{4}s_{3}s_{1}s_{2}s_{1} $,
$ s_{3}s_{4}s_{1}s_{2}s_{3}s_{1}s_{2}s_{1} $,
$ s_{3}s_{4}s_{2}s_{3}s_{1}s_{2}s_{3}s_{1} $,
$ s_{3}s_{4}s_{2}s_{3}s_{1}s_{2}s_{3}s_{2} $,
$ s_{3}s_{4}s_{2}s_{3}s_{1}s_{2}s_{3}s_{4} $,
$ s_{3}s_{4}s_{3}s_{1}s_{2}s_{3}s_{1}s_{2} $,
$ s_{3}s_{4}s_{3}s_{2}s_{3}s_{1}s_{2}s_{3} $,
$ s_{4}s_{2}s_{3}s_{1}s_{2}s_{3}s_{1}s_{2} $,
$ s_{4}s_{2}s_{3}s_{1}s_{2}s_{3}s_{2}s_{1} $,
$ s_{4}s_{3}s_{1}s_{2}s_{3}s_{4}s_{2}s_{1} $,
$ s_{4}s_{3}s_{1}s_{2}s_{3}s_{4}s_{2}s_{3} $,
$ s_{4}s_{3}s_{2}s_{3}s_{1}s_{2}s_{3}s_{2} $,
$ s_{4}s_{3}s_{2}s_{3}s_{1}s_{2}s_{3}s_{4} $,
$ s_{4}s_{3}s_{2}s_{3}s_{4}s_{1}s_{2}s_{3} $,
$ s_{4}s_{3}s_{2}s_{3}s_{4}s_{3}s_{2}s_{1} $,
$ s_{4}s_{3}s_{2}s_{3}s_{4}s_{3}s_{2}s_{3} $,
\end{quote}
Let $\Pi(-x_8^2)$ be the subset of $W$ whose elements are
\begin{quote}
  $ s_{1}s_{2}s_{3}s_{4}s_{1}s_{2}s_{3}s_{2} $,
$ s_{1}s_{2}s_{3}s_{4}s_{2}s_{3}s_{1}s_{2} $,
$ s_{1}s_{2}s_{3}s_{4}s_{2}s_{3}s_{2}s_{1} $,
$ s_{1}s_{2}s_{3}s_{4}s_{3}s_{2}s_{3}s_{1} $,
$ s_{2}s_{3}s_{1}s_{2}s_{3}s_{4}s_{2}s_{1} $,
$ s_{2}s_{3}s_{1}s_{2}s_{3}s_{4}s_{2}s_{3} $,
$ s_{2}s_{3}s_{4}s_{1}s_{2}s_{3}s_{2}s_{1} $,
$ s_{2}s_{3}s_{4}s_{2}s_{3}s_{1}s_{2}s_{1} $,
$ s_{3}s_{1}s_{2}s_{3}s_{4}s_{3}s_{2}s_{1} $,
$ s_{3}s_{1}s_{2}s_{3}s_{4}s_{3}s_{2}s_{3} $,
$ s_{3}s_{2}s_{3}s_{4}s_{3}s_{1}s_{2}s_{3} $,
$ s_{3}s_{2}s_{3}s_{4}s_{3}s_{2}s_{3}s_{1} $,
$ s_{3}s_{4}s_{3}s_{2}s_{3}s_{1}s_{2}s_{1} $,
$ s_{4}s_{2}s_{3}s_{1}s_{2}s_{3}s_{4}s_{1} $,
$ s_{4}s_{2}s_{3}s_{1}s_{2}s_{3}s_{4}s_{3} $,
$ s_{4}s_{3}s_{1}s_{2}s_{3}s_{1}s_{2}s_{1} $,
$ s_{4}s_{3}s_{1}s_{2}s_{3}s_{4}s_{1}s_{2} $,
$ s_{4}s_{3}s_{1}s_{2}s_{3}s_{4}s_{3}s_{2} $,
$ s_{4}s_{3}s_{2}s_{3}s_{4}s_{2}s_{3}s_{1} $,
$ s_{4}s_{3}s_{2}s_{3}s_{4}s_{2}s_{3}s_{2} $,
$ s_{4}s_{3}s_{2}s_{3}s_{4}s_{3}s_{1}s_{2} $.
\end{quote}

\begin{Th}[The case $F_4$, $p=3$]\label{F4mod3}For the Chow ring of the exceptional group $F_4$ over $\mathbb{F}_3$,
we have
\begin{description}
\item[(1)] The projection $\pi^*:\CH^*(G/B;\mathbb{F}_3)\to \CH^*(F_4;\mathbb{F}_3)$ is given by
$$[\Sigma_w]\longmapsto
\begin{cases}
1, & w=1\\
x_8, & w\in \Pi(x_8),\\
-x_8,& w\in \Pi(-x_8),\\
x_8^2,& w\in \Pi(x_8^2),\\
-x_8^2,& w\in \Pi(-x_8^2),\\
0, & \text{otherwise}
\end{cases}$$
\item[(2)] The comodule structure $\CH^*(G/B;\mathbb{F}_3)\to \CH^*(F_4;\mathbb{F}_3)\otimes \CH^*(G/B;\mathbb{F}_3)$ is given by
$$[\Sigma_w]\longmapsto 1\otimes [\Sigma_w]+
x_8\otimes \sum_{{w=u\odot v,} \atop {u\in \Pi(\pm x_8)}} \pm [\Sigma_{v}]+
x_8^2\otimes \sum_{{w=u\odot v,} \atop {u\in \Pi(\pm x_8^2)}} \pm [\Sigma_{v}]. $$
%(3) The coalgebra structure $\Delta:\CH^*(F_4;\mathbb{F}_3)\to \CH^*(F_4;\mathbb{F}_3)\otimes \CH^*(F_4;\mathbb{F}_3)$ is given by
%$$x_8\longmapsto x_8\otimes 1+1\otimes x_8. $$
\end{description}
\end{Th}

\noindent\textbf{Remark.}
Here we take the convention from the paper \cite{duan2014schubert}%%%%%%%%%%%%%%%%%%%%%%%%%%%%%%%55% The Chow rings of generalized Grassmannians.
that $x_8$ is the image of $[\Sigma_{s_4s_3s_2s_1}]$.

\smallbreak
\noindent\textbf{Remark.}
By the algorithm, we can only know that $\Pi(x_8^2)$ and $\Pi(-x_8^2)$ are mapped to different nonzero values under $\pi^*$.
Consider the coalgebra structure $\Delta:\CH^*(F_4;\mathbb{F}_3)\to \CH^*(F_4;\mathbb{F}_3)\otimes \CH^*(F_4;\mathbb{F}_3)$.
It is given by
$$x_8\longmapsto x_8\otimes 1+1\otimes x_8. $$
The coefficient of $x_8\otimes x_8$ in $\Delta(x_8^2)$ should be 2.
So we can use this to find which $[\Sigma_w]$ is mapped $x_8^2$ for $w\in \Pi(x_8^2)$ or $w\in \Pi(-x_8^2)$.

\section{Computations of $E_6$}

The Dynkin diagram of $E_6$ is labeled as the follows.
$$\underset{1}{\circ}
\longdash
\underset{3}{\circ}
\longdash
\stackrel{\begin{array}{@{}c@{}}
\stackrel{2}{\circ}\\[-0.3pc]\mid\\[-1.5pc]\\\mid
\end{array}}{\underset{4}{\circ}}
\longdash
\underset{5}{\circ}
\longdash
\underset{6}{\circ}$$
Denote $E_6$ the simply-connected complex group, and $\Ad E_6$ the adjoint type complex Lie group.
Denote $G/B$ the flag variety of it, and $W$ the Weyl group of $E_6$.

\paragraph{The case $p=2$.}
Now
$\CH^*(\Ad E_6;\mathbb{F}_3)=\CH^*(E_6;\mathbb{F}_3)=\mathbb{F}_2[x_6]\big/\big<x_6^2\big>$.
Denote
$$\Pi(x_6)=\left\{\begin{array}{c}
s_{2}s_{4}s_{3},
s_{2}s_{4}s_{5},
s_{3}s_{2}s_{4},
s_{3}s_{4}s_{2},\\
s_{3}s_{4}s_{5},
s_{4}s_{3}s_{2},
s_{4}s_{5}s_{2},
s_{4}s_{5}s_{3},\\
s_{5}s_{2}s_{4},
s_{5}s_{3}s_{4},
s_{5}s_{4}s_{2},
s_{5}s_{4}s_{3}
\end{array}\right\}. $$

\begin{Th}[The case $\Ad E_6$, $p=2$]\label{E6mod2}For the Chow ring of the exceptional group $\Ad E_6$ over $\mathbb{F}_2$,
we have
\begin{description}
\item[(1)] The projection $\pi^*:\CH^*(G/B;\mathbb{F}_2)\to \CH^*(\Ad E_6;\mathbb{F}_2)$ is given by
$$[\Sigma_w]\longmapsto \begin{cases}
1, & w=1\\
x_6,  & w\in \Pi(x_6),\\
0, & \text{otherwise}
\end{cases}$$
\item[(2)] The comodule structure $\CH^*(G/B;\mathbb{F}_2)\to \CH^*(\Ad E_6;\mathbb{F}_2)\otimes \CH^*(G/B;\mathbb{F}_2)$ is given by
$$[\Sigma_w]\longmapsto 1\otimes [\Sigma_w]+
x_6\otimes \sum_{{w=u\odot v,} \atop {u\in \Pi(x_6)}} [\Sigma_{v}]. $$
%
%(3) The coalgebra structure $\CH^*(F_4;\mathbb{F}_2)\to \CH^*(F_4;\mathbb{F}_2)\otimes \CH^*(F_4;\mathbb{F}_2)$ is given by
%$$x_6\longmapsto x_6\otimes 1+1\otimes x_6. $$
\end{description}
\end{Th}

\noindent\textbf{Remark.}
Here we take advantage of the work of Duan and Zhao \cite{duan2014schubert}, for one choice of $w=s_6s_5s_4s_2\in W$, such that $\pi^*[\Sigma_w]=x_8$.
%(https://arxiv.org/pdf/0801.2444v9.pdf)

\paragraph{Simply connected case, $p=3$.}
For the simply-connected case, we can get a fully description for $p=3$.
Now
$\CH^*(E_6;\mathbb{F}_3)=\mathbb{F}_3[x_8]\big/\big<x_8^3\big>$.
Let $\Pi(x_8)$ be the set of
\begin{quote}
  $ s_{1}s_{3}s_{4}s_{2} $,
$ s_{2}s_{4}s_{1}s_{3} $,
$ s_{2}s_{4}s_{5}s_{3} $,
$ s_{2}s_{4}s_{5}s_{4} $,
$ s_{3}s_{2}s_{4}s_{1} $,
$ s_{3}s_{2}s_{4}s_{3} $,
$ s_{3}s_{4}s_{5}s_{1} $,
$ s_{3}s_{4}s_{5}s_{2} $,
$ s_{3}s_{4}s_{5}s_{3} $,
$ s_{3}s_{4}s_{5}s_{6} $,
$ s_{4}s_{3}s_{2}s_{1} $,
$ s_{4}s_{5}s_{1}s_{3} $,
$ s_{4}s_{5}s_{4}s_{3} $,
$ s_{4}s_{5}s_{6}s_{2} $,
$ s_{5}s_{1}s_{3}s_{4} $,
$ s_{5}s_{3}s_{2}s_{4} $,
$ s_{5}s_{4}s_{3}s_{1} $,
$ s_{5}s_{4}s_{3}s_{2} $,
$ s_{5}s_{6}s_{2}s_{4} $,
$ s_{5}s_{6}s_{4}s_{3} $,
$ s_{6}s_{2}s_{4}s_{5} $,
$ s_{6}s_{4}s_{5}s_{3} $,
$ s_{6}s_{5}s_{3}s_{4} $,
$ s_{6}s_{5}s_{4}s_{2} $.
\end{quote}
Let $\Pi(-x_8)$ be the set of
\begin{quote}
$ s_{1}s_{3}s_{2}s_{4} $,
$ s_{1}s_{3}s_{4}s_{5} $,
$ s_{2}s_{4}s_{3}s_{1} $,
$ s_{2}s_{4}s_{5}s_{6} $,
$ s_{3}s_{2}s_{4}s_{5} $,
$ s_{3}s_{4}s_{2}s_{1} $,
$ s_{3}s_{4}s_{3}s_{2} $,
$ s_{3}s_{4}s_{5}s_{4} $,
$ s_{4}s_{1}s_{3}s_{2} $,
$ s_{4}s_{5}s_{3}s_{1} $,
$ s_{4}s_{5}s_{3}s_{2} $,
$ s_{4}s_{5}s_{4}s_{2} $,
$ s_{4}s_{5}s_{6}s_{3} $,
$ s_{5}s_{2}s_{4}s_{3} $,
$ s_{5}s_{3}s_{4}s_{1} $,
$ s_{5}s_{3}s_{4}s_{2} $,
$ s_{5}s_{3}s_{4}s_{3} $,
$ s_{5}s_{4}s_{1}s_{3} $,
$ s_{5}s_{6}s_{3}s_{4} $,
$ s_{5}s_{6}s_{4}s_{2} $,
$ s_{6}s_{3}s_{4}s_{5} $,
$ s_{6}s_{4}s_{5}s_{2} $,
$ s_{6}s_{5}s_{2}s_{4} $,
$ s_{6}s_{5}s_{4}s_{3} $.
\end{quote}
Let $\Pi(x_8^2)$ be the set of
\begin{quote}
$ s_{1}s_{3}s_{2}s_{4}s_{5}s_{1}s_{3}s_{2} $,
$ s_{1}s_{3}s_{2}s_{4}s_{5}s_{1}s_{3}s_{4} $,
$ s_{1}s_{3}s_{2}s_{4}s_{5}s_{3}s_{2}s_{4} $,
$ s_{1}s_{3}s_{2}s_{4}s_{5}s_{6}s_{3}s_{1} $,
$ s_{1}s_{3}s_{2}s_{4}s_{5}s_{6}s_{4}s_{3} $,
$ s_{1}s_{3}s_{4}s_{5}s_{1}s_{3}s_{2}s_{4} $,
$ s_{1}s_{3}s_{4}s_{5}s_{2}s_{4}s_{3}s_{1} $,
$ s_{1}s_{3}s_{4}s_{5}s_{3}s_{2}s_{4}s_{2} $,
$ s_{1}s_{3}s_{4}s_{5}s_{3}s_{4}s_{3}s_{1} $,
$ s_{1}s_{3}s_{4}s_{5}s_{6}s_{2}s_{4}s_{5} $,
$ s_{1}s_{3}s_{4}s_{5}s_{6}s_{3}s_{2}s_{4} $,
$ s_{1}s_{3}s_{4}s_{5}s_{6}s_{3}s_{4}s_{1} $,
$ s_{1}s_{3}s_{4}s_{5}s_{6}s_{4}s_{1}s_{3} $,
$ s_{1}s_{3}s_{4}s_{5}s_{6}s_{4}s_{3}s_{2} $,
$ s_{1}s_{3}s_{4}s_{5}s_{6}s_{4}s_{5}s_{3} $,
$ s_{1}s_{3}s_{4}s_{5}s_{6}s_{5}s_{3}s_{4} $,
$ s_{1}s_{3}s_{4}s_{5}s_{6}s_{5}s_{4}s_{2} $,
$ s_{2}s_{4}s_{5}s_{1}s_{3}s_{2}s_{4}s_{3} $,
$ s_{2}s_{4}s_{5}s_{1}s_{3}s_{4}s_{2}s_{1} $,
$ s_{2}s_{4}s_{5}s_{1}s_{3}s_{4}s_{3}s_{2} $,
$ s_{2}s_{4}s_{5}s_{3}s_{2}s_{4}s_{1}s_{3} $,
$ s_{2}s_{4}s_{5}s_{3}s_{2}s_{4}s_{3}s_{1} $,
$ s_{2}s_{4}s_{5}s_{4}s_{1}s_{3}s_{2}s_{4} $,
$ s_{2}s_{4}s_{5}s_{4}s_{3}s_{2}s_{4}s_{3} $,
$ s_{2}s_{4}s_{5}s_{6}s_{1}s_{3}s_{2}s_{4} $,
$ s_{2}s_{4}s_{5}s_{6}s_{2}s_{4}s_{1}s_{3} $,
$ s_{2}s_{4}s_{5}s_{6}s_{2}s_{4}s_{5}s_{2} $,
$ s_{2}s_{4}s_{5}s_{6}s_{3}s_{2}s_{4}s_{3} $,
$ s_{2}s_{4}s_{5}s_{6}s_{3}s_{2}s_{4}s_{5} $,
$ s_{2}s_{4}s_{5}s_{6}s_{3}s_{4}s_{1}s_{3} $,
$ s_{2}s_{4}s_{5}s_{6}s_{3}s_{4}s_{2}s_{1} $,
$ s_{2}s_{4}s_{5}s_{6}s_{3}s_{4}s_{3}s_{2} $,
$ s_{2}s_{4}s_{5}s_{6}s_{3}s_{4}s_{5}s_{1} $,
$ s_{2}s_{4}s_{5}s_{6}s_{3}s_{4}s_{5}s_{3} $,
$ s_{2}s_{4}s_{5}s_{6}s_{4}s_{3}s_{2}s_{1} $,
$ s_{2}s_{4}s_{5}s_{6}s_{4}s_{3}s_{2}s_{4} $,
$ s_{2}s_{4}s_{5}s_{6}s_{4}s_{5}s_{1}s_{3} $,
$ s_{2}s_{4}s_{5}s_{6}s_{4}s_{5}s_{3}s_{2} $,
$ s_{2}s_{4}s_{5}s_{6}s_{5}s_{1}s_{3}s_{4} $,
$ s_{2}s_{4}s_{5}s_{6}s_{5}s_{2}s_{4}s_{3} $,
$ s_{2}s_{4}s_{5}s_{6}s_{5}s_{3}s_{4}s_{2} $,
$ s_{2}s_{4}s_{5}s_{6}s_{5}s_{4}s_{3}s_{1} $,
$ s_{3}s_{2}s_{4}s_{5}s_{1}s_{3}s_{2}s_{1} $,
$ s_{3}s_{2}s_{4}s_{5}s_{1}s_{3}s_{2}s_{4} $,
$ s_{3}s_{2}s_{4}s_{5}s_{2}s_{4}s_{1}s_{3} $,
$ s_{3}s_{2}s_{4}s_{5}s_{2}s_{4}s_{3}s_{2} $,
$ s_{3}s_{2}s_{4}s_{5}s_{3}s_{2}s_{4}s_{1} $,
$ s_{3}s_{2}s_{4}s_{5}s_{3}s_{4}s_{1}s_{3} $,
$ s_{3}s_{2}s_{4}s_{5}s_{3}s_{4}s_{2}s_{1} $,
$ s_{3}s_{2}s_{4}s_{5}s_{3}s_{4}s_{3}s_{2} $,
$ s_{3}s_{2}s_{4}s_{5}s_{6}s_{1}s_{3}s_{4} $,
$ s_{3}s_{2}s_{4}s_{5}s_{6}s_{2}s_{4}s_{3} $,
$ s_{3}s_{2}s_{4}s_{5}s_{6}s_{3}s_{2}s_{1} $,
$ s_{3}s_{2}s_{4}s_{5}s_{6}s_{4}s_{5}s_{2} $,
$ s_{3}s_{2}s_{4}s_{5}s_{6}s_{5}s_{2}s_{4} $,
$ s_{3}s_{4}s_{1}s_{3}s_{2}s_{4}s_{3}s_{2} $,
$ s_{3}s_{4}s_{1}s_{3}s_{2}s_{4}s_{5}s_{4} $,
$ s_{3}s_{4}s_{1}s_{3}s_{2}s_{4}s_{5}s_{6} $,
$ s_{3}s_{4}s_{5}s_{1}s_{3}s_{2}s_{4}s_{2} $,
$ s_{3}s_{4}s_{5}s_{1}s_{3}s_{2}s_{4}s_{3} $,
$ s_{3}s_{4}s_{5}s_{1}s_{3}s_{4}s_{3}s_{2} $,
$ s_{3}s_{4}s_{5}s_{2}s_{4}s_{1}s_{3}s_{2} $,
$ s_{3}s_{4}s_{5}s_{3}s_{2}s_{4}s_{1}s_{3} $,
$ s_{3}s_{4}s_{5}s_{3}s_{2}s_{4}s_{2}s_{1} $,
$ s_{3}s_{4}s_{5}s_{4}s_{1}s_{3}s_{2}s_{4} $,
$ s_{3}s_{4}s_{5}s_{6}s_{2}s_{4}s_{3}s_{1} $,
$ s_{3}s_{4}s_{5}s_{6}s_{2}s_{4}s_{3}s_{2} $,
$ s_{3}s_{4}s_{5}s_{6}s_{2}s_{4}s_{5}s_{3} $,
$ s_{3}s_{4}s_{5}s_{6}s_{3}s_{2}s_{4}s_{5} $,
$ s_{3}s_{4}s_{5}s_{6}s_{3}s_{4}s_{5}s_{2} $,
$ s_{3}s_{4}s_{5}s_{6}s_{4}s_{1}s_{3}s_{1} $,
$ s_{3}s_{4}s_{5}s_{6}s_{4}s_{1}s_{3}s_{2} $,
$ s_{3}s_{4}s_{5}s_{6}s_{4}s_{3}s_{2}s_{1} $,
$ s_{3}s_{4}s_{5}s_{6}s_{4}s_{5}s_{1}s_{3} $,
$ s_{3}s_{4}s_{5}s_{6}s_{4}s_{5}s_{2}s_{1} $,
$ s_{3}s_{4}s_{5}s_{6}s_{4}s_{5}s_{3}s_{1} $,
$ s_{3}s_{4}s_{5}s_{6}s_{5}s_{1}s_{3}s_{4} $,
$ s_{3}s_{4}s_{5}s_{6}s_{5}s_{2}s_{4}s_{1} $,
$ s_{3}s_{4}s_{5}s_{6}s_{5}s_{2}s_{4}s_{2} $,
$ s_{3}s_{4}s_{5}s_{6}s_{5}s_{3}s_{4}s_{1} $,
$ s_{3}s_{4}s_{5}s_{6}s_{5}s_{3}s_{4}s_{3} $,
$ s_{4}s_{1}s_{3}s_{2}s_{4}s_{5}s_{2}s_{1} $,
$ s_{4}s_{1}s_{3}s_{2}s_{4}s_{5}s_{3}s_{1} $,
$ s_{4}s_{1}s_{3}s_{2}s_{4}s_{5}s_{6}s_{4} $,
$ s_{4}s_{3}s_{2}s_{4}s_{1}s_{3}s_{2}s_{1} $,
$ s_{4}s_{3}s_{2}s_{4}s_{5}s_{2}s_{4}s_{2} $,
$ s_{4}s_{3}s_{2}s_{4}s_{5}s_{3}s_{2}s_{4} $,
$ s_{4}s_{3}s_{2}s_{4}s_{5}s_{3}s_{4}s_{1} $,
$ s_{4}s_{3}s_{2}s_{4}s_{5}s_{4}s_{1}s_{3} $,
$ s_{4}s_{3}s_{2}s_{4}s_{5}s_{4}s_{2}s_{1} $,
$ s_{4}s_{3}s_{2}s_{4}s_{5}s_{6}s_{1}s_{3} $,
$ s_{4}s_{3}s_{2}s_{4}s_{5}s_{6}s_{3}s_{4} $,
$ s_{4}s_{3}s_{2}s_{4}s_{5}s_{6}s_{4}s_{2} $,
$ s_{4}s_{3}s_{2}s_{4}s_{5}s_{6}s_{4}s_{5} $,
$ s_{4}s_{5}s_{1}s_{3}s_{2}s_{4}s_{3}s_{2} $,
$ s_{4}s_{5}s_{3}s_{2}s_{4}s_{1}s_{3}s_{1} $,
$ s_{4}s_{5}s_{3}s_{2}s_{4}s_{1}s_{3}s_{2} $,
$ s_{4}s_{5}s_{3}s_{2}s_{4}s_{3}s_{2}s_{1} $,
$ s_{4}s_{5}s_{3}s_{4}s_{1}s_{3}s_{2}s_{1} $,
$ s_{4}s_{5}s_{4}s_{1}s_{3}s_{2}s_{4}s_{2} $,
$ s_{4}s_{5}s_{4}s_{3}s_{2}s_{4}s_{1}s_{3} $,
$ s_{4}s_{5}s_{4}s_{3}s_{2}s_{4}s_{3}s_{2} $,
$ s_{4}s_{5}s_{6}s_{1}s_{3}s_{2}s_{4}s_{2} $,
$ s_{4}s_{5}s_{6}s_{1}s_{3}s_{2}s_{4}s_{3} $,
$ s_{4}s_{5}s_{6}s_{1}s_{3}s_{4}s_{2}s_{1} $,
$ s_{4}s_{5}s_{6}s_{1}s_{3}s_{4}s_{3}s_{1} $,
$ s_{4}s_{5}s_{6}s_{1}s_{3}s_{4}s_{5}s_{1} $,
$ s_{4}s_{5}s_{6}s_{1}s_{3}s_{4}s_{5}s_{2} $,
$ s_{4}s_{5}s_{6}s_{2}s_{4}s_{1}s_{3}s_{2} $,
$ s_{4}s_{5}s_{6}s_{2}s_{4}s_{5}s_{3}s_{1} $,
$ s_{4}s_{5}s_{6}s_{3}s_{2}s_{4}s_{3}s_{2} $,
$ s_{4}s_{5}s_{6}s_{3}s_{2}s_{4}s_{5}s_{2} $,
$ s_{4}s_{5}s_{6}s_{3}s_{2}s_{4}s_{5}s_{3} $,
$ s_{4}s_{5}s_{6}s_{4}s_{1}s_{3}s_{2}s_{4} $,
$ s_{4}s_{5}s_{6}s_{4}s_{3}s_{2}s_{4}s_{2} $,
$ s_{4}s_{5}s_{6}s_{4}s_{3}s_{2}s_{4}s_{5} $,
$ s_{4}s_{5}s_{6}s_{5}s_{1}s_{3}s_{4}s_{3} $,
$ s_{4}s_{5}s_{6}s_{5}s_{2}s_{4}s_{1}s_{3} $,
$ s_{4}s_{5}s_{6}s_{5}s_{2}s_{4}s_{3}s_{2} $,
$ s_{4}s_{5}s_{6}s_{5}s_{3}s_{4}s_{1}s_{3} $,
$ s_{4}s_{5}s_{6}s_{5}s_{3}s_{4}s_{2}s_{1} $,
$ s_{4}s_{5}s_{6}s_{5}s_{3}s_{4}s_{3}s_{1} $,
$ s_{5}s_{1}s_{3}s_{2}s_{4}s_{1}s_{3}s_{2} $,
$ s_{5}s_{3}s_{2}s_{4}s_{1}s_{3}s_{2}s_{1} $,
$ s_{5}s_{3}s_{4}s_{1}s_{3}s_{2}s_{4}s_{5} $,
$ s_{5}s_{4}s_{1}s_{3}s_{2}s_{4}s_{2}s_{1} $,
$ s_{5}s_{4}s_{1}s_{3}s_{2}s_{4}s_{3}s_{1} $,
$ s_{5}s_{4}s_{1}s_{3}s_{2}s_{4}s_{5}s_{2} $,
$ s_{5}s_{4}s_{1}s_{3}s_{2}s_{4}s_{5}s_{3} $,
$ s_{5}s_{4}s_{3}s_{2}s_{4}s_{5}s_{1}s_{3} $,
$ s_{5}s_{4}s_{3}s_{2}s_{4}s_{5}s_{2}s_{1} $,
$ s_{5}s_{4}s_{3}s_{2}s_{4}s_{5}s_{3}s_{1} $,
$ s_{5}s_{4}s_{3}s_{2}s_{4}s_{5}s_{3}s_{2} $,
$ s_{5}s_{4}s_{3}s_{2}s_{4}s_{5}s_{6}s_{1} $,
$ s_{5}s_{6}s_{1}s_{3}s_{2}s_{4}s_{5}s_{3} $,
$ s_{5}s_{6}s_{1}s_{3}s_{4}s_{1}s_{3}s_{2} $,
$ s_{5}s_{6}s_{1}s_{3}s_{4}s_{5}s_{1}s_{3} $,
$ s_{5}s_{6}s_{1}s_{3}s_{4}s_{5}s_{2}s_{4} $,
$ s_{5}s_{6}s_{1}s_{3}s_{4}s_{5}s_{4}s_{3} $,
$ s_{5}s_{6}s_{2}s_{4}s_{5}s_{1}s_{3}s_{2} $,
$ s_{5}s_{6}s_{2}s_{4}s_{5}s_{2}s_{4}s_{2} $,
$ s_{5}s_{6}s_{2}s_{4}s_{5}s_{3}s_{2}s_{4} $,
$ s_{5}s_{6}s_{2}s_{4}s_{5}s_{3}s_{4}s_{1} $,
$ s_{5}s_{6}s_{2}s_{4}s_{5}s_{3}s_{4}s_{3} $,
$ s_{5}s_{6}s_{2}s_{4}s_{5}s_{4}s_{1}s_{3} $,
$ s_{5}s_{6}s_{2}s_{4}s_{5}s_{4}s_{3}s_{2} $,
$ s_{5}s_{6}s_{3}s_{2}s_{4}s_{1}s_{3}s_{2} $,
$ s_{5}s_{6}s_{3}s_{2}s_{4}s_{5}s_{2}s_{1} $,
$ s_{5}s_{6}s_{3}s_{2}s_{4}s_{5}s_{4}s_{2} $,
$ s_{5}s_{6}s_{3}s_{4}s_{1}s_{3}s_{2}s_{1} $,
$ s_{5}s_{6}s_{3}s_{4}s_{5}s_{1}s_{3}s_{1} $,
$ s_{5}s_{6}s_{3}s_{4}s_{5}s_{2}s_{4}s_{3} $,
$ s_{5}s_{6}s_{3}s_{4}s_{5}s_{3}s_{2}s_{1} $,
$ s_{5}s_{6}s_{3}s_{4}s_{5}s_{3}s_{2}s_{4} $,
$ s_{5}s_{6}s_{3}s_{4}s_{5}s_{3}s_{4}s_{2} $,
$ s_{5}s_{6}s_{3}s_{4}s_{5}s_{4}s_{1}s_{3} $,
$ s_{5}s_{6}s_{3}s_{4}s_{5}s_{4}s_{2}s_{1} $,
$ s_{5}s_{6}s_{3}s_{4}s_{5}s_{4}s_{3}s_{1} $,
$ s_{5}s_{6}s_{4}s_{1}s_{3}s_{2}s_{4}s_{1} $,
$ s_{5}s_{6}s_{4}s_{1}s_{3}s_{2}s_{4}s_{2} $,
$ s_{5}s_{6}s_{4}s_{1}s_{3}s_{2}s_{4}s_{3} $,
$ s_{5}s_{6}s_{4}s_{3}s_{2}s_{4}s_{5}s_{2} $,
$ s_{5}s_{6}s_{4}s_{3}s_{2}s_{4}s_{5}s_{3} $,
$ s_{5}s_{6}s_{4}s_{3}s_{2}s_{4}s_{5}s_{4} $,
$ s_{5}s_{6}s_{4}s_{5}s_{1}s_{3}s_{4}s_{1} $,
$ s_{5}s_{6}s_{4}s_{5}s_{1}s_{3}s_{4}s_{2} $,
$ s_{5}s_{6}s_{4}s_{5}s_{2}s_{4}s_{3}s_{1} $,
$ s_{5}s_{6}s_{4}s_{5}s_{3}s_{2}s_{4}s_{2} $,
$ s_{5}s_{6}s_{4}s_{5}s_{3}s_{2}s_{4}s_{3} $,
$ s_{5}s_{6}s_{4}s_{5}s_{4}s_{3}s_{2}s_{4} $,
$ s_{6}s_{1}s_{3}s_{2}s_{4}s_{5}s_{1}s_{3} $,
$ s_{6}s_{1}s_{3}s_{2}s_{4}s_{5}s_{4}s_{2} $,
$ s_{6}s_{1}s_{3}s_{2}s_{4}s_{5}s_{4}s_{3} $,
$ s_{6}s_{1}s_{3}s_{4}s_{5}s_{2}s_{4}s_{3} $,
$ s_{6}s_{1}s_{3}s_{4}s_{5}s_{3}s_{4}s_{1} $,
$ s_{6}s_{1}s_{3}s_{4}s_{5}s_{4}s_{1}s_{3} $,
$ s_{6}s_{2}s_{4}s_{5}s_{1}s_{3}s_{4}s_{1} $,
$ s_{6}s_{2}s_{4}s_{5}s_{1}s_{3}s_{4}s_{2} $,
$ s_{6}s_{2}s_{4}s_{5}s_{3}s_{2}s_{4}s_{1} $,
$ s_{6}s_{2}s_{4}s_{5}s_{4}s_{3}s_{2}s_{4} $,
$ s_{6}s_{3}s_{2}s_{4}s_{5}s_{1}s_{3}s_{1} $,
$ s_{6}s_{3}s_{2}s_{4}s_{5}s_{1}s_{3}s_{2} $,
$ s_{6}s_{3}s_{2}s_{4}s_{5}s_{2}s_{4}s_{1} $,
$ s_{6}s_{3}s_{2}s_{4}s_{5}s_{2}s_{4}s_{3} $,
$ s_{6}s_{3}s_{2}s_{4}s_{5}s_{3}s_{2}s_{4} $,
$ s_{6}s_{3}s_{2}s_{4}s_{5}s_{3}s_{4}s_{2} $,
$ s_{6}s_{3}s_{2}s_{4}s_{5}s_{4}s_{1}s_{3} $,
$ s_{6}s_{3}s_{2}s_{4}s_{5}s_{4}s_{3}s_{1} $,
$ s_{6}s_{3}s_{4}s_{5}s_{1}s_{3}s_{4}s_{2} $,
$ s_{6}s_{3}s_{4}s_{5}s_{2}s_{4}s_{2}s_{1} $,
$ s_{6}s_{3}s_{4}s_{5}s_{2}s_{4}s_{3}s_{2} $,
$ s_{6}s_{3}s_{4}s_{5}s_{3}s_{2}s_{4}s_{1} $,
$ s_{6}s_{3}s_{4}s_{5}s_{3}s_{2}s_{4}s_{2} $,
$ s_{6}s_{3}s_{4}s_{5}s_{4}s_{1}s_{3}s_{1} $,
$ s_{6}s_{3}s_{4}s_{5}s_{4}s_{3}s_{2}s_{1} $,
$ s_{6}s_{4}s_{1}s_{3}s_{2}s_{4}s_{5}s_{1} $,
$ s_{6}s_{4}s_{3}s_{2}s_{4}s_{5}s_{2}s_{1} $,
$ s_{6}s_{4}s_{3}s_{2}s_{4}s_{5}s_{3}s_{1} $,
$ s_{6}s_{4}s_{3}s_{2}s_{4}s_{5}s_{3}s_{4} $,
$ s_{6}s_{4}s_{3}s_{2}s_{4}s_{5}s_{4}s_{1} $,
$ s_{6}s_{4}s_{3}s_{2}s_{4}s_{5}s_{4}s_{2} $,
$ s_{6}s_{4}s_{5}s_{1}s_{3}s_{4}s_{3}s_{1} $,
$ s_{6}s_{4}s_{5}s_{3}s_{2}s_{4}s_{2}s_{1} $,
$ s_{6}s_{4}s_{5}s_{3}s_{2}s_{4}s_{3}s_{1} $,
$ s_{6}s_{4}s_{5}s_{3}s_{4}s_{1}s_{3}s_{2} $,
$ s_{6}s_{4}s_{5}s_{4}s_{3}s_{2}s_{4}s_{1} $,
$ s_{6}s_{4}s_{5}s_{4}s_{3}s_{2}s_{4}s_{2} $,
$ s_{6}s_{5}s_{1}s_{3}s_{2}s_{4}s_{3}s_{2} $,
$ s_{6}s_{5}s_{1}s_{3}s_{4}s_{3}s_{2}s_{1} $,
$ s_{6}s_{5}s_{3}s_{4}s_{1}s_{3}s_{2}s_{4} $,
$ s_{6}s_{5}s_{4}s_{1}s_{3}s_{2}s_{4}s_{5} $,
$ s_{6}s_{5}s_{4}s_{3}s_{2}s_{4}s_{1}s_{3} $,
$ s_{6}s_{5}s_{4}s_{3}s_{2}s_{4}s_{5}s_{6} $.
\end{quote}
and $\Pi(-x_8^2)$ be the set of
\begin{quote}
$ s_{1}s_{3}s_{2}s_{4}s_{5}s_{2}s_{4}s_{3} $,
$ s_{1}s_{3}s_{2}s_{4}s_{5}s_{3}s_{2}s_{1} $,
$ s_{1}s_{3}s_{2}s_{4}s_{5}s_{3}s_{4}s_{1} $,
$ s_{1}s_{3}s_{2}s_{4}s_{5}s_{3}s_{4}s_{3} $,
$ s_{1}s_{3}s_{2}s_{4}s_{5}s_{6}s_{1}s_{3} $,
$ s_{1}s_{3}s_{2}s_{4}s_{5}s_{6}s_{3}s_{4} $,
$ s_{1}s_{3}s_{4}s_{5}s_{1}s_{3}s_{4}s_{3} $,
$ s_{1}s_{3}s_{4}s_{5}s_{2}s_{4}s_{1}s_{3} $,
$ s_{1}s_{3}s_{4}s_{5}s_{3}s_{2}s_{4}s_{1} $,
$ s_{1}s_{3}s_{4}s_{5}s_{3}s_{2}s_{4}s_{3} $,
$ s_{1}s_{3}s_{4}s_{5}s_{4}s_{3}s_{2}s_{4} $,
$ s_{1}s_{3}s_{4}s_{5}s_{6}s_{1}s_{3}s_{4} $,
$ s_{1}s_{3}s_{4}s_{5}s_{6}s_{2}s_{4}s_{3} $,
$ s_{1}s_{3}s_{4}s_{5}s_{6}s_{3}s_{4}s_{2} $,
$ s_{1}s_{3}s_{4}s_{5}s_{6}s_{3}s_{4}s_{3} $,
$ s_{1}s_{3}s_{4}s_{5}s_{6}s_{3}s_{4}s_{5} $,
$ s_{1}s_{3}s_{4}s_{5}s_{6}s_{4}s_{3}s_{1} $,
$ s_{1}s_{3}s_{4}s_{5}s_{6}s_{4}s_{5}s_{2} $,
$ s_{1}s_{3}s_{4}s_{5}s_{6}s_{5}s_{2}s_{4} $,
$ s_{1}s_{3}s_{4}s_{5}s_{6}s_{5}s_{4}s_{3} $,
$ s_{2}s_{4}s_{5}s_{1}s_{3}s_{2}s_{4}s_{1} $,
$ s_{2}s_{4}s_{5}s_{3}s_{4}s_{1}s_{3}s_{2} $,
$ s_{2}s_{4}s_{5}s_{4}s_{3}s_{2}s_{4}s_{1} $,
$ s_{2}s_{4}s_{5}s_{6}s_{1}s_{3}s_{4}s_{1} $,
$ s_{2}s_{4}s_{5}s_{6}s_{1}s_{3}s_{4}s_{2} $,
$ s_{2}s_{4}s_{5}s_{6}s_{1}s_{3}s_{4}s_{3} $,
$ s_{2}s_{4}s_{5}s_{6}s_{1}s_{3}s_{4}s_{5} $,
$ s_{2}s_{4}s_{5}s_{6}s_{2}s_{4}s_{3}s_{1} $,
$ s_{2}s_{4}s_{5}s_{6}s_{2}s_{4}s_{3}s_{2} $,
$ s_{2}s_{4}s_{5}s_{6}s_{2}s_{4}s_{5}s_{3} $,
$ s_{2}s_{4}s_{5}s_{6}s_{3}s_{2}s_{4}s_{1} $,
$ s_{2}s_{4}s_{5}s_{6}s_{3}s_{4}s_{5}s_{2} $,
$ s_{2}s_{4}s_{5}s_{6}s_{4}s_{1}s_{3}s_{2} $,
$ s_{2}s_{4}s_{5}s_{6}s_{4}s_{5}s_{3}s_{1} $,
$ s_{2}s_{4}s_{5}s_{6}s_{5}s_{2}s_{4}s_{2} $,
$ s_{2}s_{4}s_{5}s_{6}s_{5}s_{3}s_{2}s_{4} $,
$ s_{2}s_{4}s_{5}s_{6}s_{5}s_{3}s_{4}s_{1} $,
$ s_{2}s_{4}s_{5}s_{6}s_{5}s_{3}s_{4}s_{3} $,
$ s_{2}s_{4}s_{5}s_{6}s_{5}s_{4}s_{1}s_{3} $,
$ s_{2}s_{4}s_{5}s_{6}s_{5}s_{4}s_{3}s_{2} $,
$ s_{3}s_{2}s_{4}s_{5}s_{1}s_{3}s_{4}s_{2} $,
$ s_{3}s_{2}s_{4}s_{5}s_{3}s_{2}s_{4}s_{3} $,
$ s_{3}s_{2}s_{4}s_{5}s_{4}s_{3}s_{2}s_{4} $,
$ s_{3}s_{2}s_{4}s_{5}s_{6}s_{1}s_{3}s_{1} $,
$ s_{3}s_{2}s_{4}s_{5}s_{6}s_{1}s_{3}s_{2} $,
$ s_{3}s_{2}s_{4}s_{5}s_{6}s_{2}s_{4}s_{5} $,
$ s_{3}s_{2}s_{4}s_{5}s_{6}s_{3}s_{2}s_{4} $,
$ s_{3}s_{2}s_{4}s_{5}s_{6}s_{4}s_{1}s_{3} $,
$ s_{3}s_{2}s_{4}s_{5}s_{6}s_{5}s_{4}s_{2} $,
$ s_{3}s_{4}s_{1}s_{3}s_{2}s_{4}s_{5}s_{2} $,
$ s_{3}s_{4}s_{1}s_{3}s_{2}s_{4}s_{5}s_{3} $,
$ s_{3}s_{4}s_{5}s_{2}s_{4}s_{1}s_{3}s_{1} $,
$ s_{3}s_{4}s_{5}s_{2}s_{4}s_{3}s_{2}s_{1} $,
$ s_{3}s_{4}s_{5}s_{3}s_{2}s_{4}s_{3}s_{1} $,
$ s_{3}s_{4}s_{5}s_{3}s_{2}s_{4}s_{3}s_{2} $,
$ s_{3}s_{4}s_{5}s_{3}s_{4}s_{1}s_{3}s_{1} $,
$ s_{3}s_{4}s_{5}s_{3}s_{4}s_{3}s_{2}s_{1} $,
$ s_{3}s_{4}s_{5}s_{4}s_{3}s_{2}s_{4}s_{2} $,
$ s_{3}s_{4}s_{5}s_{4}s_{3}s_{2}s_{4}s_{3} $,
$ s_{3}s_{4}s_{5}s_{6}s_{1}s_{3}s_{4}s_{2} $,
$ s_{3}s_{4}s_{5}s_{6}s_{1}s_{3}s_{4}s_{3} $,
$ s_{3}s_{4}s_{5}s_{6}s_{1}s_{3}s_{4}s_{5} $,
$ s_{3}s_{4}s_{5}s_{6}s_{2}s_{4}s_{5}s_{1} $,
$ s_{3}s_{4}s_{5}s_{6}s_{2}s_{4}s_{5}s_{2} $,
$ s_{3}s_{4}s_{5}s_{6}s_{3}s_{2}s_{4}s_{1} $,
$ s_{3}s_{4}s_{5}s_{6}s_{3}s_{2}s_{4}s_{2} $,
$ s_{3}s_{4}s_{5}s_{6}s_{3}s_{4}s_{1}s_{3} $,
$ s_{3}s_{4}s_{5}s_{6}s_{3}s_{4}s_{2}s_{1} $,
$ s_{3}s_{4}s_{5}s_{6}s_{3}s_{4}s_{3}s_{2} $,
$ s_{3}s_{4}s_{5}s_{6}s_{3}s_{4}s_{5}s_{1} $,
$ s_{3}s_{4}s_{5}s_{6}s_{3}s_{4}s_{5}s_{3} $,
$ s_{3}s_{4}s_{5}s_{6}s_{5}s_{2}s_{4}s_{3} $,
$ s_{3}s_{4}s_{5}s_{6}s_{5}s_{3}s_{2}s_{4} $,
$ s_{3}s_{4}s_{5}s_{6}s_{5}s_{3}s_{4}s_{2} $,
$ s_{3}s_{4}s_{5}s_{6}s_{5}s_{4}s_{1}s_{3} $,
$ s_{3}s_{4}s_{5}s_{6}s_{5}s_{4}s_{2}s_{1} $,
$ s_{3}s_{4}s_{5}s_{6}s_{5}s_{4}s_{3}s_{1} $,
$ s_{4}s_{1}s_{3}s_{2}s_{4}s_{3}s_{2}s_{1} $,
$ s_{4}s_{1}s_{3}s_{2}s_{4}s_{5}s_{3}s_{4} $,
$ s_{4}s_{1}s_{3}s_{2}s_{4}s_{5}s_{4}s_{1} $,
$ s_{4}s_{1}s_{3}s_{2}s_{4}s_{5}s_{4}s_{2} $,
$ s_{4}s_{1}s_{3}s_{2}s_{4}s_{5}s_{6}s_{1} $,
$ s_{4}s_{3}s_{2}s_{4}s_{5}s_{1}s_{3}s_{1} $,
$ s_{4}s_{3}s_{2}s_{4}s_{5}s_{1}s_{3}s_{2} $,
$ s_{4}s_{3}s_{2}s_{4}s_{5}s_{2}s_{4}s_{1} $,
$ s_{4}s_{3}s_{2}s_{4}s_{5}s_{2}s_{4}s_{3} $,
$ s_{4}s_{3}s_{2}s_{4}s_{5}s_{3}s_{2}s_{1} $,
$ s_{4}s_{3}s_{2}s_{4}s_{5}s_{3}s_{4}s_{3} $,
$ s_{4}s_{3}s_{2}s_{4}s_{5}s_{4}s_{3}s_{1} $,
$ s_{4}s_{3}s_{2}s_{4}s_{5}s_{4}s_{3}s_{2} $,
$ s_{4}s_{3}s_{2}s_{4}s_{5}s_{6}s_{2}s_{1} $,
$ s_{4}s_{3}s_{2}s_{4}s_{5}s_{6}s_{3}s_{1} $,
$ s_{4}s_{3}s_{2}s_{4}s_{5}s_{6}s_{4}s_{1} $,
$ s_{4}s_{3}s_{2}s_{4}s_{5}s_{6}s_{4}s_{3} $,
$ s_{4}s_{3}s_{2}s_{4}s_{5}s_{6}s_{5}s_{4} $,
$ s_{4}s_{5}s_{1}s_{3}s_{2}s_{4}s_{2}s_{1} $,
$ s_{4}s_{5}s_{1}s_{3}s_{2}s_{4}s_{3}s_{1} $,
$ s_{4}s_{5}s_{1}s_{3}s_{4}s_{3}s_{2}s_{1} $,
$ s_{4}s_{5}s_{4}s_{1}s_{3}s_{2}s_{4}s_{1} $,
$ s_{4}s_{5}s_{4}s_{1}s_{3}s_{2}s_{4}s_{3} $,
$ s_{4}s_{5}s_{4}s_{3}s_{2}s_{4}s_{2}s_{1} $,
$ s_{4}s_{5}s_{6}s_{1}s_{3}s_{4}s_{3}s_{2} $,
$ s_{4}s_{5}s_{6}s_{1}s_{3}s_{4}s_{5}s_{3} $,
$ s_{4}s_{5}s_{6}s_{2}s_{4}s_{3}s_{2}s_{1} $,
$ s_{4}s_{5}s_{6}s_{2}s_{4}s_{5}s_{1}s_{3} $,
$ s_{4}s_{5}s_{6}s_{2}s_{4}s_{5}s_{3}s_{2} $,
$ s_{4}s_{5}s_{6}s_{3}s_{2}s_{4}s_{2}s_{1} $,
$ s_{4}s_{5}s_{6}s_{3}s_{2}s_{4}s_{3}s_{1} $,
$ s_{4}s_{5}s_{6}s_{3}s_{4}s_{1}s_{3}s_{1} $,
$ s_{4}s_{5}s_{6}s_{3}s_{4}s_{1}s_{3}s_{2} $,
$ s_{4}s_{5}s_{6}s_{3}s_{4}s_{3}s_{2}s_{1} $,
$ s_{4}s_{5}s_{6}s_{3}s_{4}s_{5}s_{1}s_{3} $,
$ s_{4}s_{5}s_{6}s_{3}s_{4}s_{5}s_{2}s_{1} $,
$ s_{4}s_{5}s_{6}s_{3}s_{4}s_{5}s_{3}s_{1} $,
$ s_{4}s_{5}s_{6}s_{4}s_{3}s_{2}s_{4}s_{1} $,
$ s_{4}s_{5}s_{6}s_{4}s_{3}s_{2}s_{4}s_{3} $,
$ s_{4}s_{5}s_{6}s_{5}s_{1}s_{3}s_{4}s_{1} $,
$ s_{4}s_{5}s_{6}s_{5}s_{1}s_{3}s_{4}s_{2} $,
$ s_{4}s_{5}s_{6}s_{5}s_{2}s_{4}s_{3}s_{1} $,
$ s_{4}s_{5}s_{6}s_{5}s_{3}s_{2}s_{4}s_{2} $,
$ s_{4}s_{5}s_{6}s_{5}s_{3}s_{2}s_{4}s_{3} $,
$ s_{4}s_{5}s_{6}s_{5}s_{4}s_{3}s_{2}s_{4} $,
$ s_{5}s_{1}s_{3}s_{2}s_{4}s_{3}s_{2}s_{1} $,
$ s_{5}s_{3}s_{4}s_{1}s_{3}s_{2}s_{4}s_{2} $,
$ s_{5}s_{3}s_{4}s_{1}s_{3}s_{2}s_{4}s_{3} $,
$ s_{5}s_{4}s_{1}s_{3}s_{2}s_{4}s_{3}s_{2} $,
$ s_{5}s_{4}s_{1}s_{3}s_{2}s_{4}s_{5}s_{1} $,
$ s_{5}s_{4}s_{1}s_{3}s_{2}s_{4}s_{5}s_{6} $,
$ s_{5}s_{4}s_{3}s_{2}s_{4}s_{1}s_{3}s_{1} $,
$ s_{5}s_{4}s_{3}s_{2}s_{4}s_{1}s_{3}s_{2} $,
$ s_{5}s_{4}s_{3}s_{2}s_{4}s_{5}s_{2}s_{4} $,
$ s_{5}s_{4}s_{3}s_{2}s_{4}s_{5}s_{3}s_{4} $,
$ s_{5}s_{4}s_{3}s_{2}s_{4}s_{5}s_{4}s_{2} $,
$ s_{5}s_{4}s_{3}s_{2}s_{4}s_{5}s_{4}s_{3} $,
$ s_{5}s_{4}s_{3}s_{2}s_{4}s_{5}s_{6}s_{2} $,
$ s_{5}s_{4}s_{3}s_{2}s_{4}s_{5}s_{6}s_{4} $,
$ s_{5}s_{4}s_{3}s_{2}s_{4}s_{5}s_{6}s_{5} $,
$ s_{5}s_{6}s_{1}s_{3}s_{2}s_{4}s_{3}s_{2} $,
$ s_{5}s_{6}s_{1}s_{3}s_{2}s_{4}s_{5}s_{2} $,
$ s_{5}s_{6}s_{1}s_{3}s_{4}s_{3}s_{2}s_{1} $,
$ s_{5}s_{6}s_{1}s_{3}s_{4}s_{5}s_{3}s_{1} $,
$ s_{5}s_{6}s_{1}s_{3}s_{4}s_{5}s_{3}s_{2} $,
$ s_{5}s_{6}s_{1}s_{3}s_{4}s_{5}s_{3}s_{4} $,
$ s_{5}s_{6}s_{1}s_{3}s_{4}s_{5}s_{4}s_{2} $,
$ s_{5}s_{6}s_{2}s_{4}s_{5}s_{1}s_{3}s_{4} $,
$ s_{5}s_{6}s_{2}s_{4}s_{5}s_{2}s_{4}s_{3} $,
$ s_{5}s_{6}s_{2}s_{4}s_{5}s_{3}s_{2}s_{1} $,
$ s_{5}s_{6}s_{2}s_{4}s_{5}s_{3}s_{4}s_{2} $,
$ s_{5}s_{6}s_{2}s_{4}s_{5}s_{4}s_{3}s_{1} $,
$ s_{5}s_{6}s_{3}s_{2}s_{4}s_{5}s_{2}s_{4} $,
$ s_{5}s_{6}s_{3}s_{2}s_{4}s_{5}s_{3}s_{1} $,
$ s_{5}s_{6}s_{3}s_{4}s_{1}s_{3}s_{2}s_{4} $,
$ s_{5}s_{6}s_{3}s_{4}s_{5}s_{1}s_{3}s_{4} $,
$ s_{5}s_{6}s_{3}s_{4}s_{5}s_{2}s_{4}s_{1} $,
$ s_{5}s_{6}s_{3}s_{4}s_{5}s_{2}s_{4}s_{2} $,
$ s_{5}s_{6}s_{3}s_{4}s_{5}s_{3}s_{4}s_{1} $,
$ s_{5}s_{6}s_{3}s_{4}s_{5}s_{3}s_{4}s_{3} $,
$ s_{5}s_{6}s_{4}s_{3}s_{2}s_{4}s_{1}s_{3} $,
$ s_{5}s_{6}s_{4}s_{5}s_{1}s_{3}s_{4}s_{3} $,
$ s_{5}s_{6}s_{4}s_{5}s_{2}s_{4}s_{1}s_{3} $,
$ s_{5}s_{6}s_{4}s_{5}s_{2}s_{4}s_{3}s_{2} $,
$ s_{5}s_{6}s_{4}s_{5}s_{3}s_{4}s_{1}s_{3} $,
$ s_{5}s_{6}s_{4}s_{5}s_{3}s_{4}s_{2}s_{1} $,
$ s_{5}s_{6}s_{4}s_{5}s_{3}s_{4}s_{3}s_{1} $,
$ s_{5}s_{6}s_{5}s_{4}s_{3}s_{2}s_{4}s_{5} $,
$ s_{6}s_{1}s_{3}s_{2}s_{4}s_{5}s_{2}s_{4} $,
$ s_{6}s_{1}s_{3}s_{2}s_{4}s_{5}s_{3}s_{1} $,
$ s_{6}s_{1}s_{3}s_{2}s_{4}s_{5}s_{3}s_{4} $,
$ s_{6}s_{1}s_{3}s_{4}s_{5}s_{1}s_{3}s_{4} $,
$ s_{6}s_{1}s_{3}s_{4}s_{5}s_{2}s_{4}s_{2} $,
$ s_{6}s_{1}s_{3}s_{4}s_{5}s_{3}s_{2}s_{4} $,
$ s_{6}s_{1}s_{3}s_{4}s_{5}s_{4}s_{3}s_{1} $,
$ s_{6}s_{2}s_{4}s_{5}s_{1}s_{3}s_{2}s_{4} $,
$ s_{6}s_{2}s_{4}s_{5}s_{1}s_{3}s_{4}s_{3} $,
$ s_{6}s_{2}s_{4}s_{5}s_{3}s_{2}s_{4}s_{3} $,
$ s_{6}s_{2}s_{4}s_{5}s_{3}s_{4}s_{1}s_{3} $,
$ s_{6}s_{2}s_{4}s_{5}s_{3}s_{4}s_{2}s_{1} $,
$ s_{6}s_{2}s_{4}s_{5}s_{3}s_{4}s_{3}s_{1} $,
$ s_{6}s_{3}s_{2}s_{4}s_{5}s_{1}s_{3}s_{4} $,
$ s_{6}s_{3}s_{2}s_{4}s_{5}s_{2}s_{4}s_{2} $,
$ s_{6}s_{3}s_{2}s_{4}s_{5}s_{3}s_{2}s_{1} $,
$ s_{6}s_{3}s_{2}s_{4}s_{5}s_{3}s_{4}s_{1} $,
$ s_{6}s_{3}s_{2}s_{4}s_{5}s_{3}s_{4}s_{3} $,
$ s_{6}s_{3}s_{2}s_{4}s_{5}s_{4}s_{2}s_{1} $,
$ s_{6}s_{3}s_{4}s_{1}s_{3}s_{2}s_{4}s_{5} $,
$ s_{6}s_{3}s_{4}s_{5}s_{1}s_{3}s_{4}s_{3} $,
$ s_{6}s_{3}s_{4}s_{5}s_{2}s_{4}s_{3}s_{1} $,
$ s_{6}s_{3}s_{4}s_{5}s_{3}s_{4}s_{1}s_{3} $,
$ s_{6}s_{3}s_{4}s_{5}s_{3}s_{4}s_{2}s_{1} $,
$ s_{6}s_{3}s_{4}s_{5}s_{3}s_{4}s_{3}s_{1} $,
$ s_{6}s_{3}s_{4}s_{5}s_{4}s_{1}s_{3}s_{2} $,
$ s_{6}s_{3}s_{4}s_{5}s_{4}s_{3}s_{2}s_{4} $,
$ s_{6}s_{4}s_{1}s_{3}s_{2}s_{4}s_{5}s_{4} $,
$ s_{6}s_{4}s_{3}s_{2}s_{4}s_{5}s_{1}s_{3} $,
$ s_{6}s_{4}s_{3}s_{2}s_{4}s_{5}s_{2}s_{4} $,
$ s_{6}s_{4}s_{5}s_{1}s_{3}s_{2}s_{4}s_{2} $,
$ s_{6}s_{4}s_{5}s_{1}s_{3}s_{2}s_{4}s_{3} $,
$ s_{6}s_{4}s_{5}s_{1}s_{3}s_{4}s_{2}s_{1} $,
$ s_{6}s_{4}s_{5}s_{3}s_{2}s_{4}s_{3}s_{2} $,
$ s_{6}s_{4}s_{5}s_{3}s_{4}s_{1}s_{3}s_{1} $,
$ s_{6}s_{4}s_{5}s_{3}s_{4}s_{3}s_{2}s_{1} $,
$ s_{6}s_{4}s_{5}s_{4}s_{1}s_{3}s_{2}s_{4} $,
$ s_{6}s_{5}s_{1}s_{3}s_{4}s_{1}s_{3}s_{2} $,
$ s_{6}s_{5}s_{3}s_{2}s_{4}s_{1}s_{3}s_{2} $,
$ s_{6}s_{5}s_{3}s_{4}s_{1}s_{3}s_{2}s_{1} $,
$ s_{6}s_{5}s_{4}s_{1}s_{3}s_{2}s_{4}s_{1} $,
$ s_{6}s_{5}s_{4}s_{1}s_{3}s_{2}s_{4}s_{2} $,
$ s_{6}s_{5}s_{4}s_{1}s_{3}s_{2}s_{4}s_{3} $,
$ s_{6}s_{5}s_{4}s_{3}s_{2}s_{4}s_{5}s_{1} $,
$ s_{6}s_{5}s_{4}s_{3}s_{2}s_{4}s_{5}s_{3} $,
$ s_{6}s_{5}s_{4}s_{3}s_{2}s_{4}s_{5}s_{4} $.
\end{quote}

\begin{Th}[The case $E_6$, $p=3$]\label{scE6mod3}For the Chow ring of the exceptional group $E_6$ over $\mathbb{F}_3$,
we have
\begin{description}
\item[(1)] The projection $\pi^*:\CH^*(G/B;\mathbb{F}_3)\to \CH^*(E_6;\mathbb{F}_3)$ is given by
$$[\Sigma_w]\longmapsto
\begin{cases}
1, & w=1\\
x_8, & w\in \Pi(x_8),\\
-x_8,& w\in \Pi(-x_8),\\
x_8^2,& w\in \Pi(x_8^2),\\
-x_8^2,& w\in \Pi(-x_8^2),\\
0, & \text{otherwise}
\end{cases}$$
\item[(2)] The comodule structure $\CH^*(G/B;\mathbb{F}_3)\to \CH^*(F_4;\mathbb{F}_3)\otimes \CH^*(G/B;\mathbb{F}_3)$ is given by
$$[\Sigma_w]\longmapsto 1\otimes [\Sigma_w]+
x_8\otimes \sum_{{w=u\odot v,} \atop {u\in \Pi(\pm x_8)}} \pm [\Sigma_{v}]+
x_8^2\otimes \sum_{{w=u\odot v,} \atop {u\in \Pi(\pm x_8^2)}} \pm [\Sigma_{v}]. $$
%(3) The coalgebra structure $\Delta:\CH^*(F_4;\mathbb{F}_3)\to \CH^*(F_4;\mathbb{F}_3)\otimes \CH^*(F_4;\mathbb{F}_3)$ is given by
%$$x_8\longmapsto x_8\otimes 1+1\otimes x_8. $$
\end{description}
\end{Th}

%p=W(s[1]*s[3]*s[2]*s[4]*s[5]*s[1]*s[3]*s[2]) #E6p3
%k=0
%for b1 in [0..len(Lengthk)-1]:
%    for b2 in [0..len(Lengthk)-1]:
%        if p==Lengthk[b1] * Lengthk[b2]:
%            k=k+Image[b1]*Image[b2]
%k

\section{Computations of $E_7$}

The Dynkin diagram of $E_7$ is labeled as follows.
$$\underset{1}{\circ}
\longdash
\underset{3}{\circ}
\longdash
\stackrel{\begin{array}{@{}c@{}}
\stackrel{2}{\circ}\\[-0.3pc]\mid\\[-1.5pc]\\\mid
\end{array}}{\underset{4}{\circ}}
\longdash
\underset{5}{\circ}
\longdash
\underset{6}{\circ}
\longdash
\underset{7}{\circ}$$
We denote $s_i$ the corresponding simple reflection with $i\in \{1,\ldots,7\}$.

Denote $\Ad E_7$ the adjoint type complex Lie group.
Denote $G/B$ the flag variety of it, and $W$ the Weyl group of $E_7$.

\paragraph{The case $p=3$. }
Now
$\CH^*(\Ad E_7;\mathbb{F}_3)=\CH^*(E_7;\mathbb{F}_3)=\mathbb{F}_3[x_8]\big/\big<x_8^3\big> $.
Denote $\Pi(x_8)$ the set of
\begin{quote}
$ s_{1}s_{3}s_{4}s_{2} $,
$ s_{2}s_{4}s_{1}s_{3} $,
$ s_{2}s_{4}s_{5}s_{3} $,
$ s_{2}s_{4}s_{5}s_{4} $,
$ s_{3}s_{2}s_{4}s_{1} $,
$ s_{3}s_{2}s_{4}s_{3} $,
$ s_{3}s_{4}s_{5}s_{1} $,
$ s_{3}s_{4}s_{5}s_{2} $,
$ s_{3}s_{4}s_{5}s_{3} $,
$ s_{3}s_{4}s_{5}s_{6} $,
$ s_{4}s_{3}s_{2}s_{1} $,
$ s_{4}s_{5}s_{1}s_{3} $,
$ s_{4}s_{5}s_{4}s_{3} $,
$ s_{4}s_{5}s_{6}s_{2} $,
$ s_{5}s_{1}s_{3}s_{4} $,
$ s_{5}s_{3}s_{2}s_{4} $,
$ s_{5}s_{4}s_{3}s_{1} $,
$ s_{5}s_{4}s_{3}s_{2} $,
$ s_{5}s_{6}s_{2}s_{4} $,
$ s_{5}s_{6}s_{4}s_{3} $,
$ s_{6}s_{2}s_{4}s_{5} $,
$ s_{6}s_{4}s_{5}s_{3} $,
$ s_{6}s_{5}s_{3}s_{4} $,
$ s_{6}s_{5}s_{4}s_{2} $;
\end{quote}
denote $\Pi(-x_8)$ the set of
\begin{quote}
$ s_{1}s_{3}s_{2}s_{4} $,
$ s_{1}s_{3}s_{4}s_{5} $,
$ s_{2}s_{4}s_{3}s_{1} $,
$ s_{2}s_{4}s_{5}s_{6} $,
$ s_{3}s_{2}s_{4}s_{5} $,
$ s_{3}s_{4}s_{2}s_{1} $,
$ s_{3}s_{4}s_{3}s_{2} $,
$ s_{3}s_{4}s_{5}s_{4} $,
$ s_{4}s_{1}s_{3}s_{2} $,
$ s_{4}s_{5}s_{3}s_{1} $,
$ s_{4}s_{5}s_{3}s_{2} $,
$ s_{4}s_{5}s_{4}s_{2} $,
$ s_{4}s_{5}s_{6}s_{3} $,
$ s_{5}s_{2}s_{4}s_{3} $,
$ s_{5}s_{3}s_{4}s_{1} $,
$ s_{5}s_{3}s_{4}s_{2} $,
$ s_{5}s_{3}s_{4}s_{3} $,
$ s_{5}s_{4}s_{1}s_{3} $,
$ s_{5}s_{6}s_{3}s_{4} $,
$ s_{5}s_{6}s_{4}s_{2} $,
$ s_{6}s_{3}s_{4}s_{5} $,
$ s_{6}s_{4}s_{5}s_{2} $,
$ s_{6}s_{5}s_{2}s_{4} $,
$ s_{6}s_{5}s_{4}s_{3} $.
\end{quote}
We also denote $\Pi(x_8^2)$ the set of
\begin{quote}
$ s_{1}s_{3}s_{2}s_{4}s_{5}s_{1}s_{3}s_{2} $,
$ s_{1}s_{3}s_{2}s_{4}s_{5}s_{1}s_{3}s_{4} $,
$ s_{1}s_{3}s_{2}s_{4}s_{5}s_{3}s_{2}s_{4} $,
$ s_{1}s_{3}s_{2}s_{4}s_{5}s_{6}s_{3}s_{1} $,
$ s_{1}s_{3}s_{2}s_{4}s_{5}s_{6}s_{4}s_{3} $,
$ s_{1}s_{3}s_{4}s_{5}s_{1}s_{3}s_{2}s_{4} $,
$ s_{1}s_{3}s_{4}s_{5}s_{2}s_{4}s_{3}s_{1} $,
$ s_{1}s_{3}s_{4}s_{5}s_{3}s_{2}s_{4}s_{2} $,
$ s_{1}s_{3}s_{4}s_{5}s_{3}s_{4}s_{3}s_{1} $,
$ s_{1}s_{3}s_{4}s_{5}s_{6}s_{2}s_{4}s_{5} $,
$ s_{1}s_{3}s_{4}s_{5}s_{6}s_{3}s_{2}s_{4} $,
$ s_{1}s_{3}s_{4}s_{5}s_{6}s_{3}s_{4}s_{1} $,
$ s_{1}s_{3}s_{4}s_{5}s_{6}s_{4}s_{1}s_{3} $,
$ s_{1}s_{3}s_{4}s_{5}s_{6}s_{4}s_{3}s_{2} $,
$ s_{1}s_{3}s_{4}s_{5}s_{6}s_{4}s_{5}s_{3} $,
$ s_{1}s_{3}s_{4}s_{5}s_{6}s_{5}s_{3}s_{4} $,
$ s_{1}s_{3}s_{4}s_{5}s_{6}s_{5}s_{4}s_{2} $,
$ s_{2}s_{4}s_{5}s_{1}s_{3}s_{2}s_{4}s_{3} $,
$ s_{2}s_{4}s_{5}s_{1}s_{3}s_{4}s_{2}s_{1} $,
$ s_{2}s_{4}s_{5}s_{1}s_{3}s_{4}s_{3}s_{2} $,
$ s_{2}s_{4}s_{5}s_{3}s_{2}s_{4}s_{1}s_{3} $,
$ s_{2}s_{4}s_{5}s_{3}s_{2}s_{4}s_{3}s_{1} $,
$ s_{2}s_{4}s_{5}s_{4}s_{1}s_{3}s_{2}s_{4} $,
$ s_{2}s_{4}s_{5}s_{4}s_{3}s_{2}s_{4}s_{3} $,
$ s_{2}s_{4}s_{5}s_{6}s_{1}s_{3}s_{2}s_{4} $,
$ s_{2}s_{4}s_{5}s_{6}s_{2}s_{4}s_{1}s_{3} $,
$ s_{2}s_{4}s_{5}s_{6}s_{2}s_{4}s_{5}s_{2} $,
$ s_{2}s_{4}s_{5}s_{6}s_{3}s_{2}s_{4}s_{3} $,
$ s_{2}s_{4}s_{5}s_{6}s_{3}s_{2}s_{4}s_{5} $,
$ s_{2}s_{4}s_{5}s_{6}s_{3}s_{4}s_{1}s_{3} $,
$ s_{2}s_{4}s_{5}s_{6}s_{3}s_{4}s_{2}s_{1} $,
$ s_{2}s_{4}s_{5}s_{6}s_{3}s_{4}s_{3}s_{2} $,
$ s_{2}s_{4}s_{5}s_{6}s_{3}s_{4}s_{5}s_{1} $,
$ s_{2}s_{4}s_{5}s_{6}s_{3}s_{4}s_{5}s_{3} $,
$ s_{2}s_{4}s_{5}s_{6}s_{4}s_{3}s_{2}s_{1} $,
$ s_{2}s_{4}s_{5}s_{6}s_{4}s_{3}s_{2}s_{4} $,
$ s_{2}s_{4}s_{5}s_{6}s_{4}s_{5}s_{1}s_{3} $,
$ s_{2}s_{4}s_{5}s_{6}s_{4}s_{5}s_{3}s_{2} $,
$ s_{2}s_{4}s_{5}s_{6}s_{5}s_{1}s_{3}s_{4} $,
$ s_{2}s_{4}s_{5}s_{6}s_{5}s_{2}s_{4}s_{3} $,
$ s_{2}s_{4}s_{5}s_{6}s_{5}s_{3}s_{4}s_{2} $,
$ s_{2}s_{4}s_{5}s_{6}s_{5}s_{4}s_{3}s_{1} $,
$ s_{3}s_{2}s_{4}s_{5}s_{1}s_{3}s_{2}s_{1} $,
$ s_{3}s_{2}s_{4}s_{5}s_{1}s_{3}s_{2}s_{4} $,
$ s_{3}s_{2}s_{4}s_{5}s_{2}s_{4}s_{1}s_{3} $,
$ s_{3}s_{2}s_{4}s_{5}s_{2}s_{4}s_{3}s_{2} $,
$ s_{3}s_{2}s_{4}s_{5}s_{3}s_{2}s_{4}s_{1} $,
$ s_{3}s_{2}s_{4}s_{5}s_{3}s_{4}s_{1}s_{3} $,
$ s_{3}s_{2}s_{4}s_{5}s_{3}s_{4}s_{2}s_{1} $,
$ s_{3}s_{2}s_{4}s_{5}s_{3}s_{4}s_{3}s_{2} $,
$ s_{3}s_{2}s_{4}s_{5}s_{6}s_{1}s_{3}s_{4} $,
$ s_{3}s_{2}s_{4}s_{5}s_{6}s_{2}s_{4}s_{3} $,
$ s_{3}s_{2}s_{4}s_{5}s_{6}s_{3}s_{2}s_{1} $,
$ s_{3}s_{2}s_{4}s_{5}s_{6}s_{4}s_{5}s_{2} $,
$ s_{3}s_{2}s_{4}s_{5}s_{6}s_{5}s_{2}s_{4} $,
$ s_{3}s_{4}s_{1}s_{3}s_{2}s_{4}s_{3}s_{2} $,
$ s_{3}s_{4}s_{1}s_{3}s_{2}s_{4}s_{5}s_{4} $,
$ s_{3}s_{4}s_{1}s_{3}s_{2}s_{4}s_{5}s_{6} $,
$ s_{3}s_{4}s_{5}s_{1}s_{3}s_{2}s_{4}s_{2} $,
$ s_{3}s_{4}s_{5}s_{1}s_{3}s_{2}s_{4}s_{3} $,
$ s_{3}s_{4}s_{5}s_{1}s_{3}s_{4}s_{3}s_{2} $,
$ s_{3}s_{4}s_{5}s_{2}s_{4}s_{1}s_{3}s_{2} $,
$ s_{3}s_{4}s_{5}s_{3}s_{2}s_{4}s_{1}s_{3} $,
$ s_{3}s_{4}s_{5}s_{3}s_{2}s_{4}s_{2}s_{1} $,
$ s_{3}s_{4}s_{5}s_{4}s_{1}s_{3}s_{2}s_{4} $,
$ s_{3}s_{4}s_{5}s_{6}s_{2}s_{4}s_{3}s_{1} $,
$ s_{3}s_{4}s_{5}s_{6}s_{2}s_{4}s_{3}s_{2} $,
$ s_{3}s_{4}s_{5}s_{6}s_{2}s_{4}s_{5}s_{3} $,
$ s_{3}s_{4}s_{5}s_{6}s_{3}s_{2}s_{4}s_{5} $,
$ s_{3}s_{4}s_{5}s_{6}s_{3}s_{4}s_{5}s_{2} $,
$ s_{3}s_{4}s_{5}s_{6}s_{4}s_{1}s_{3}s_{1} $,
$ s_{3}s_{4}s_{5}s_{6}s_{4}s_{1}s_{3}s_{2} $,
$ s_{3}s_{4}s_{5}s_{6}s_{4}s_{3}s_{2}s_{1} $,
$ s_{3}s_{4}s_{5}s_{6}s_{4}s_{5}s_{1}s_{3} $,
$ s_{3}s_{4}s_{5}s_{6}s_{4}s_{5}s_{2}s_{1} $,
$ s_{3}s_{4}s_{5}s_{6}s_{4}s_{5}s_{3}s_{1} $,
$ s_{3}s_{4}s_{5}s_{6}s_{5}s_{1}s_{3}s_{4} $,
$ s_{3}s_{4}s_{5}s_{6}s_{5}s_{2}s_{4}s_{1} $,
$ s_{3}s_{4}s_{5}s_{6}s_{5}s_{2}s_{4}s_{2} $,
$ s_{3}s_{4}s_{5}s_{6}s_{5}s_{3}s_{4}s_{1} $,
$ s_{3}s_{4}s_{5}s_{6}s_{5}s_{3}s_{4}s_{3} $,
$ s_{4}s_{1}s_{3}s_{2}s_{4}s_{5}s_{2}s_{1} $,
$ s_{4}s_{1}s_{3}s_{2}s_{4}s_{5}s_{3}s_{1} $,
$ s_{4}s_{1}s_{3}s_{2}s_{4}s_{5}s_{6}s_{4} $,
$ s_{4}s_{3}s_{2}s_{4}s_{1}s_{3}s_{2}s_{1} $,
$ s_{4}s_{3}s_{2}s_{4}s_{5}s_{2}s_{4}s_{2} $,
$ s_{4}s_{3}s_{2}s_{4}s_{5}s_{3}s_{2}s_{4} $,
$ s_{4}s_{3}s_{2}s_{4}s_{5}s_{3}s_{4}s_{1} $,
$ s_{4}s_{3}s_{2}s_{4}s_{5}s_{4}s_{1}s_{3} $,
$ s_{4}s_{3}s_{2}s_{4}s_{5}s_{4}s_{2}s_{1} $,
$ s_{4}s_{3}s_{2}s_{4}s_{5}s_{6}s_{1}s_{3} $,
$ s_{4}s_{3}s_{2}s_{4}s_{5}s_{6}s_{3}s_{4} $,
$ s_{4}s_{3}s_{2}s_{4}s_{5}s_{6}s_{4}s_{2} $,
$ s_{4}s_{3}s_{2}s_{4}s_{5}s_{6}s_{4}s_{5} $,
$ s_{4}s_{5}s_{1}s_{3}s_{2}s_{4}s_{3}s_{2} $,
$ s_{4}s_{5}s_{3}s_{2}s_{4}s_{1}s_{3}s_{1} $,
$ s_{4}s_{5}s_{3}s_{2}s_{4}s_{1}s_{3}s_{2} $,
$ s_{4}s_{5}s_{3}s_{2}s_{4}s_{3}s_{2}s_{1} $,
$ s_{4}s_{5}s_{3}s_{4}s_{1}s_{3}s_{2}s_{1} $,
$ s_{4}s_{5}s_{4}s_{1}s_{3}s_{2}s_{4}s_{2} $,
$ s_{4}s_{5}s_{4}s_{3}s_{2}s_{4}s_{1}s_{3} $,
$ s_{4}s_{5}s_{4}s_{3}s_{2}s_{4}s_{3}s_{2} $,
$ s_{4}s_{5}s_{6}s_{1}s_{3}s_{2}s_{4}s_{2} $,
$ s_{4}s_{5}s_{6}s_{1}s_{3}s_{2}s_{4}s_{3} $,
$ s_{4}s_{5}s_{6}s_{1}s_{3}s_{4}s_{2}s_{1} $,
$ s_{4}s_{5}s_{6}s_{1}s_{3}s_{4}s_{3}s_{1} $,
$ s_{4}s_{5}s_{6}s_{1}s_{3}s_{4}s_{5}s_{1} $,
$ s_{4}s_{5}s_{6}s_{1}s_{3}s_{4}s_{5}s_{2} $,
$ s_{4}s_{5}s_{6}s_{2}s_{4}s_{1}s_{3}s_{2} $,
$ s_{4}s_{5}s_{6}s_{2}s_{4}s_{5}s_{3}s_{1} $,
$ s_{4}s_{5}s_{6}s_{3}s_{2}s_{4}s_{3}s_{2} $,
$ s_{4}s_{5}s_{6}s_{3}s_{2}s_{4}s_{5}s_{2} $,
$ s_{4}s_{5}s_{6}s_{3}s_{2}s_{4}s_{5}s_{3} $,
$ s_{4}s_{5}s_{6}s_{4}s_{1}s_{3}s_{2}s_{4} $,
$ s_{4}s_{5}s_{6}s_{4}s_{3}s_{2}s_{4}s_{2} $,
$ s_{4}s_{5}s_{6}s_{4}s_{3}s_{2}s_{4}s_{5} $,
$ s_{4}s_{5}s_{6}s_{5}s_{1}s_{3}s_{4}s_{3} $,
$ s_{4}s_{5}s_{6}s_{5}s_{2}s_{4}s_{1}s_{3} $,
$ s_{4}s_{5}s_{6}s_{5}s_{2}s_{4}s_{3}s_{2} $,
$ s_{4}s_{5}s_{6}s_{5}s_{3}s_{4}s_{1}s_{3} $,
$ s_{4}s_{5}s_{6}s_{5}s_{3}s_{4}s_{2}s_{1} $,
$ s_{4}s_{5}s_{6}s_{5}s_{3}s_{4}s_{3}s_{1} $,
$ s_{5}s_{1}s_{3}s_{2}s_{4}s_{1}s_{3}s_{2} $,
$ s_{5}s_{3}s_{2}s_{4}s_{1}s_{3}s_{2}s_{1} $,
$ s_{5}s_{3}s_{4}s_{1}s_{3}s_{2}s_{4}s_{5} $,
$ s_{5}s_{4}s_{1}s_{3}s_{2}s_{4}s_{2}s_{1} $,
$ s_{5}s_{4}s_{1}s_{3}s_{2}s_{4}s_{3}s_{1} $,
$ s_{5}s_{4}s_{1}s_{3}s_{2}s_{4}s_{5}s_{2} $,
$ s_{5}s_{4}s_{1}s_{3}s_{2}s_{4}s_{5}s_{3} $,
$ s_{5}s_{4}s_{3}s_{2}s_{4}s_{5}s_{1}s_{3} $,
$ s_{5}s_{4}s_{3}s_{2}s_{4}s_{5}s_{2}s_{1} $,
$ s_{5}s_{4}s_{3}s_{2}s_{4}s_{5}s_{3}s_{1} $,
$ s_{5}s_{4}s_{3}s_{2}s_{4}s_{5}s_{3}s_{2} $,
$ s_{5}s_{4}s_{3}s_{2}s_{4}s_{5}s_{6}s_{1} $,
$ s_{5}s_{6}s_{1}s_{3}s_{2}s_{4}s_{5}s_{3} $,
$ s_{5}s_{6}s_{1}s_{3}s_{4}s_{1}s_{3}s_{2} $,
$ s_{5}s_{6}s_{1}s_{3}s_{4}s_{5}s_{1}s_{3} $,
$ s_{5}s_{6}s_{1}s_{3}s_{4}s_{5}s_{2}s_{4} $,
$ s_{5}s_{6}s_{1}s_{3}s_{4}s_{5}s_{4}s_{3} $,
$ s_{5}s_{6}s_{2}s_{4}s_{5}s_{1}s_{3}s_{2} $,
$ s_{5}s_{6}s_{2}s_{4}s_{5}s_{2}s_{4}s_{2} $,
$ s_{5}s_{6}s_{2}s_{4}s_{5}s_{3}s_{2}s_{4} $,
$ s_{5}s_{6}s_{2}s_{4}s_{5}s_{3}s_{4}s_{1} $,
$ s_{5}s_{6}s_{2}s_{4}s_{5}s_{3}s_{4}s_{3} $,
$ s_{5}s_{6}s_{2}s_{4}s_{5}s_{4}s_{1}s_{3} $,
$ s_{5}s_{6}s_{2}s_{4}s_{5}s_{4}s_{3}s_{2} $,
$ s_{5}s_{6}s_{3}s_{2}s_{4}s_{1}s_{3}s_{2} $,
$ s_{5}s_{6}s_{3}s_{2}s_{4}s_{5}s_{2}s_{1} $,
$ s_{5}s_{6}s_{3}s_{2}s_{4}s_{5}s_{4}s_{2} $,
$ s_{5}s_{6}s_{3}s_{4}s_{1}s_{3}s_{2}s_{1} $,
$ s_{5}s_{6}s_{3}s_{4}s_{5}s_{1}s_{3}s_{1} $,
$ s_{5}s_{6}s_{3}s_{4}s_{5}s_{2}s_{4}s_{3} $,
$ s_{5}s_{6}s_{3}s_{4}s_{5}s_{3}s_{2}s_{1} $,
$ s_{5}s_{6}s_{3}s_{4}s_{5}s_{3}s_{2}s_{4} $,
$ s_{5}s_{6}s_{3}s_{4}s_{5}s_{3}s_{4}s_{2} $,
$ s_{5}s_{6}s_{3}s_{4}s_{5}s_{4}s_{1}s_{3} $,
$ s_{5}s_{6}s_{3}s_{4}s_{5}s_{4}s_{2}s_{1} $,
$ s_{5}s_{6}s_{3}s_{4}s_{5}s_{4}s_{3}s_{1} $,
$ s_{5}s_{6}s_{4}s_{1}s_{3}s_{2}s_{4}s_{1} $,
$ s_{5}s_{6}s_{4}s_{1}s_{3}s_{2}s_{4}s_{2} $,
$ s_{5}s_{6}s_{4}s_{1}s_{3}s_{2}s_{4}s_{3} $,
$ s_{5}s_{6}s_{4}s_{3}s_{2}s_{4}s_{5}s_{2} $,
$ s_{5}s_{6}s_{4}s_{3}s_{2}s_{4}s_{5}s_{3} $,
$ s_{5}s_{6}s_{4}s_{3}s_{2}s_{4}s_{5}s_{4} $,
$ s_{5}s_{6}s_{4}s_{5}s_{1}s_{3}s_{4}s_{1} $,
$ s_{5}s_{6}s_{4}s_{5}s_{1}s_{3}s_{4}s_{2} $,
$ s_{5}s_{6}s_{4}s_{5}s_{2}s_{4}s_{3}s_{1} $,
$ s_{5}s_{6}s_{4}s_{5}s_{3}s_{2}s_{4}s_{2} $,
$ s_{5}s_{6}s_{4}s_{5}s_{3}s_{2}s_{4}s_{3} $,
$ s_{5}s_{6}s_{4}s_{5}s_{4}s_{3}s_{2}s_{4} $,
$ s_{6}s_{1}s_{3}s_{2}s_{4}s_{5}s_{1}s_{3} $,
$ s_{6}s_{1}s_{3}s_{2}s_{4}s_{5}s_{4}s_{2} $,
$ s_{6}s_{1}s_{3}s_{2}s_{4}s_{5}s_{4}s_{3} $,
$ s_{6}s_{1}s_{3}s_{4}s_{5}s_{2}s_{4}s_{3} $,
$ s_{6}s_{1}s_{3}s_{4}s_{5}s_{3}s_{4}s_{1} $,
$ s_{6}s_{1}s_{3}s_{4}s_{5}s_{4}s_{1}s_{3} $,
$ s_{6}s_{2}s_{4}s_{5}s_{1}s_{3}s_{4}s_{1} $,
$ s_{6}s_{2}s_{4}s_{5}s_{1}s_{3}s_{4}s_{2} $,
$ s_{6}s_{2}s_{4}s_{5}s_{3}s_{2}s_{4}s_{1} $,
$ s_{6}s_{2}s_{4}s_{5}s_{4}s_{3}s_{2}s_{4} $,
$ s_{6}s_{3}s_{2}s_{4}s_{5}s_{1}s_{3}s_{1} $,
$ s_{6}s_{3}s_{2}s_{4}s_{5}s_{1}s_{3}s_{2} $,
$ s_{6}s_{3}s_{2}s_{4}s_{5}s_{2}s_{4}s_{1} $,
$ s_{6}s_{3}s_{2}s_{4}s_{5}s_{2}s_{4}s_{3} $,
$ s_{6}s_{3}s_{2}s_{4}s_{5}s_{3}s_{2}s_{4} $,
$ s_{6}s_{3}s_{2}s_{4}s_{5}s_{3}s_{4}s_{2} $,
$ s_{6}s_{3}s_{2}s_{4}s_{5}s_{4}s_{1}s_{3} $,
$ s_{6}s_{3}s_{2}s_{4}s_{5}s_{4}s_{3}s_{1} $,
$ s_{6}s_{3}s_{4}s_{5}s_{1}s_{3}s_{4}s_{2} $,
$ s_{6}s_{3}s_{4}s_{5}s_{2}s_{4}s_{2}s_{1} $,
$ s_{6}s_{3}s_{4}s_{5}s_{2}s_{4}s_{3}s_{2} $,
$ s_{6}s_{3}s_{4}s_{5}s_{3}s_{2}s_{4}s_{1} $,
$ s_{6}s_{3}s_{4}s_{5}s_{3}s_{2}s_{4}s_{2} $,
$ s_{6}s_{3}s_{4}s_{5}s_{4}s_{1}s_{3}s_{1} $,
$ s_{6}s_{3}s_{4}s_{5}s_{4}s_{3}s_{2}s_{1} $,
$ s_{6}s_{4}s_{1}s_{3}s_{2}s_{4}s_{5}s_{1} $,
$ s_{6}s_{4}s_{3}s_{2}s_{4}s_{5}s_{2}s_{1} $,
$ s_{6}s_{4}s_{3}s_{2}s_{4}s_{5}s_{3}s_{1} $,
$ s_{6}s_{4}s_{3}s_{2}s_{4}s_{5}s_{3}s_{4} $,
$ s_{6}s_{4}s_{3}s_{2}s_{4}s_{5}s_{4}s_{1} $,
$ s_{6}s_{4}s_{3}s_{2}s_{4}s_{5}s_{4}s_{2} $,
$ s_{6}s_{4}s_{5}s_{1}s_{3}s_{4}s_{3}s_{1} $,
$ s_{6}s_{4}s_{5}s_{3}s_{2}s_{4}s_{2}s_{1} $,
$ s_{6}s_{4}s_{5}s_{3}s_{2}s_{4}s_{3}s_{1} $,
$ s_{6}s_{4}s_{5}s_{3}s_{4}s_{1}s_{3}s_{2} $,
$ s_{6}s_{4}s_{5}s_{4}s_{3}s_{2}s_{4}s_{1} $,
$ s_{6}s_{4}s_{5}s_{4}s_{3}s_{2}s_{4}s_{2} $,
$ s_{6}s_{5}s_{1}s_{3}s_{2}s_{4}s_{3}s_{2} $,
$ s_{6}s_{5}s_{1}s_{3}s_{4}s_{3}s_{2}s_{1} $,
$ s_{6}s_{5}s_{3}s_{4}s_{1}s_{3}s_{2}s_{4} $,
$ s_{6}s_{5}s_{4}s_{1}s_{3}s_{2}s_{4}s_{5} $,
$ s_{6}s_{5}s_{4}s_{3}s_{2}s_{4}s_{1}s_{3} $,
$ s_{6}s_{5}s_{4}s_{3}s_{2}s_{4}s_{5}s_{6} $;
\end{quote}
and denote $\Pi(-x_8^2)$ the set of
\begin{quote}
$ s_{1}s_{3}s_{2}s_{4}s_{5}s_{2}s_{4}s_{3} $,
$ s_{1}s_{3}s_{2}s_{4}s_{5}s_{3}s_{2}s_{1} $,
$ s_{1}s_{3}s_{2}s_{4}s_{5}s_{3}s_{4}s_{1} $,
$ s_{1}s_{3}s_{2}s_{4}s_{5}s_{3}s_{4}s_{3} $,
$ s_{1}s_{3}s_{2}s_{4}s_{5}s_{6}s_{1}s_{3} $,
$ s_{1}s_{3}s_{2}s_{4}s_{5}s_{6}s_{3}s_{4} $,
$ s_{1}s_{3}s_{4}s_{5}s_{1}s_{3}s_{4}s_{3} $,
$ s_{1}s_{3}s_{4}s_{5}s_{2}s_{4}s_{1}s_{3} $,
$ s_{1}s_{3}s_{4}s_{5}s_{3}s_{2}s_{4}s_{1} $,
$ s_{1}s_{3}s_{4}s_{5}s_{3}s_{2}s_{4}s_{3} $,
$ s_{1}s_{3}s_{4}s_{5}s_{4}s_{3}s_{2}s_{4} $,
$ s_{1}s_{3}s_{4}s_{5}s_{6}s_{1}s_{3}s_{4} $,
$ s_{1}s_{3}s_{4}s_{5}s_{6}s_{2}s_{4}s_{3} $,
$ s_{1}s_{3}s_{4}s_{5}s_{6}s_{3}s_{4}s_{2} $,
$ s_{1}s_{3}s_{4}s_{5}s_{6}s_{3}s_{4}s_{3} $,
$ s_{1}s_{3}s_{4}s_{5}s_{6}s_{3}s_{4}s_{5} $,
$ s_{1}s_{3}s_{4}s_{5}s_{6}s_{4}s_{3}s_{1} $,
$ s_{1}s_{3}s_{4}s_{5}s_{6}s_{4}s_{5}s_{2} $,
$ s_{1}s_{3}s_{4}s_{5}s_{6}s_{5}s_{2}s_{4} $,
$ s_{1}s_{3}s_{4}s_{5}s_{6}s_{5}s_{4}s_{3} $,
$ s_{2}s_{4}s_{5}s_{1}s_{3}s_{2}s_{4}s_{1} $,
$ s_{2}s_{4}s_{5}s_{3}s_{4}s_{1}s_{3}s_{2} $,
$ s_{2}s_{4}s_{5}s_{4}s_{3}s_{2}s_{4}s_{1} $,
$ s_{2}s_{4}s_{5}s_{6}s_{1}s_{3}s_{4}s_{1} $,
$ s_{2}s_{4}s_{5}s_{6}s_{1}s_{3}s_{4}s_{2} $,
$ s_{2}s_{4}s_{5}s_{6}s_{1}s_{3}s_{4}s_{3} $,
$ s_{2}s_{4}s_{5}s_{6}s_{1}s_{3}s_{4}s_{5} $,
$ s_{2}s_{4}s_{5}s_{6}s_{2}s_{4}s_{3}s_{1} $,
$ s_{2}s_{4}s_{5}s_{6}s_{2}s_{4}s_{3}s_{2} $,
$ s_{2}s_{4}s_{5}s_{6}s_{2}s_{4}s_{5}s_{3} $,
$ s_{2}s_{4}s_{5}s_{6}s_{3}s_{2}s_{4}s_{1} $,
$ s_{2}s_{4}s_{5}s_{6}s_{3}s_{4}s_{5}s_{2} $,
$ s_{2}s_{4}s_{5}s_{6}s_{4}s_{1}s_{3}s_{2} $,
$ s_{2}s_{4}s_{5}s_{6}s_{4}s_{5}s_{3}s_{1} $,
$ s_{2}s_{4}s_{5}s_{6}s_{5}s_{2}s_{4}s_{2} $,
$ s_{2}s_{4}s_{5}s_{6}s_{5}s_{3}s_{2}s_{4} $,
$ s_{2}s_{4}s_{5}s_{6}s_{5}s_{3}s_{4}s_{1} $,
$ s_{2}s_{4}s_{5}s_{6}s_{5}s_{3}s_{4}s_{3} $,
$ s_{2}s_{4}s_{5}s_{6}s_{5}s_{4}s_{1}s_{3} $,
$ s_{2}s_{4}s_{5}s_{6}s_{5}s_{4}s_{3}s_{2} $,
$ s_{3}s_{2}s_{4}s_{5}s_{1}s_{3}s_{4}s_{2} $,
$ s_{3}s_{2}s_{4}s_{5}s_{3}s_{2}s_{4}s_{3} $,
$ s_{3}s_{2}s_{4}s_{5}s_{4}s_{3}s_{2}s_{4} $,
$ s_{3}s_{2}s_{4}s_{5}s_{6}s_{1}s_{3}s_{1} $,
$ s_{3}s_{2}s_{4}s_{5}s_{6}s_{1}s_{3}s_{2} $,
$ s_{3}s_{2}s_{4}s_{5}s_{6}s_{2}s_{4}s_{5} $,
$ s_{3}s_{2}s_{4}s_{5}s_{6}s_{3}s_{2}s_{4} $,
$ s_{3}s_{2}s_{4}s_{5}s_{6}s_{4}s_{1}s_{3} $,
$ s_{3}s_{2}s_{4}s_{5}s_{6}s_{5}s_{4}s_{2} $,
$ s_{3}s_{4}s_{1}s_{3}s_{2}s_{4}s_{5}s_{2} $,
$ s_{3}s_{4}s_{1}s_{3}s_{2}s_{4}s_{5}s_{3} $,
$ s_{3}s_{4}s_{5}s_{2}s_{4}s_{1}s_{3}s_{1} $,
$ s_{3}s_{4}s_{5}s_{2}s_{4}s_{3}s_{2}s_{1} $,
$ s_{3}s_{4}s_{5}s_{3}s_{2}s_{4}s_{3}s_{1} $,
$ s_{3}s_{4}s_{5}s_{3}s_{2}s_{4}s_{3}s_{2} $,
$ s_{3}s_{4}s_{5}s_{3}s_{4}s_{1}s_{3}s_{1} $,
$ s_{3}s_{4}s_{5}s_{3}s_{4}s_{3}s_{2}s_{1} $,
$ s_{3}s_{4}s_{5}s_{4}s_{3}s_{2}s_{4}s_{2} $,
$ s_{3}s_{4}s_{5}s_{4}s_{3}s_{2}s_{4}s_{3} $,
$ s_{3}s_{4}s_{5}s_{6}s_{1}s_{3}s_{4}s_{2} $,
$ s_{3}s_{4}s_{5}s_{6}s_{1}s_{3}s_{4}s_{3} $,
$ s_{3}s_{4}s_{5}s_{6}s_{1}s_{3}s_{4}s_{5} $,
$ s_{3}s_{4}s_{5}s_{6}s_{2}s_{4}s_{5}s_{1} $,
$ s_{3}s_{4}s_{5}s_{6}s_{2}s_{4}s_{5}s_{2} $,
$ s_{3}s_{4}s_{5}s_{6}s_{3}s_{2}s_{4}s_{1} $,
$ s_{3}s_{4}s_{5}s_{6}s_{3}s_{2}s_{4}s_{2} $,
$ s_{3}s_{4}s_{5}s_{6}s_{3}s_{4}s_{1}s_{3} $,
$ s_{3}s_{4}s_{5}s_{6}s_{3}s_{4}s_{2}s_{1} $,
$ s_{3}s_{4}s_{5}s_{6}s_{3}s_{4}s_{3}s_{2} $,
$ s_{3}s_{4}s_{5}s_{6}s_{3}s_{4}s_{5}s_{1} $,
$ s_{3}s_{4}s_{5}s_{6}s_{3}s_{4}s_{5}s_{3} $,
$ s_{3}s_{4}s_{5}s_{6}s_{5}s_{2}s_{4}s_{3} $,
$ s_{3}s_{4}s_{5}s_{6}s_{5}s_{3}s_{2}s_{4} $,
$ s_{3}s_{4}s_{5}s_{6}s_{5}s_{3}s_{4}s_{2} $,
$ s_{3}s_{4}s_{5}s_{6}s_{5}s_{4}s_{1}s_{3} $,
$ s_{3}s_{4}s_{5}s_{6}s_{5}s_{4}s_{2}s_{1} $,
$ s_{3}s_{4}s_{5}s_{6}s_{5}s_{4}s_{3}s_{1} $,
$ s_{4}s_{1}s_{3}s_{2}s_{4}s_{3}s_{2}s_{1} $,
$ s_{4}s_{1}s_{3}s_{2}s_{4}s_{5}s_{3}s_{4} $,
$ s_{4}s_{1}s_{3}s_{2}s_{4}s_{5}s_{4}s_{1} $,
$ s_{4}s_{1}s_{3}s_{2}s_{4}s_{5}s_{4}s_{2} $,
$ s_{4}s_{1}s_{3}s_{2}s_{4}s_{5}s_{6}s_{1} $,
$ s_{4}s_{3}s_{2}s_{4}s_{5}s_{1}s_{3}s_{1} $,
$ s_{4}s_{3}s_{2}s_{4}s_{5}s_{1}s_{3}s_{2} $,
$ s_{4}s_{3}s_{2}s_{4}s_{5}s_{2}s_{4}s_{1} $,
$ s_{4}s_{3}s_{2}s_{4}s_{5}s_{2}s_{4}s_{3} $,
$ s_{4}s_{3}s_{2}s_{4}s_{5}s_{3}s_{2}s_{1} $,
$ s_{4}s_{3}s_{2}s_{4}s_{5}s_{3}s_{4}s_{3} $,
$ s_{4}s_{3}s_{2}s_{4}s_{5}s_{4}s_{3}s_{1} $,
$ s_{4}s_{3}s_{2}s_{4}s_{5}s_{4}s_{3}s_{2} $,
$ s_{4}s_{3}s_{2}s_{4}s_{5}s_{6}s_{2}s_{1} $,
$ s_{4}s_{3}s_{2}s_{4}s_{5}s_{6}s_{3}s_{1} $,
$ s_{4}s_{3}s_{2}s_{4}s_{5}s_{6}s_{4}s_{1} $,
$ s_{4}s_{3}s_{2}s_{4}s_{5}s_{6}s_{4}s_{3} $,
$ s_{4}s_{3}s_{2}s_{4}s_{5}s_{6}s_{5}s_{4} $,
$ s_{4}s_{5}s_{1}s_{3}s_{2}s_{4}s_{2}s_{1} $,
$ s_{4}s_{5}s_{1}s_{3}s_{2}s_{4}s_{3}s_{1} $,
$ s_{4}s_{5}s_{1}s_{3}s_{4}s_{3}s_{2}s_{1} $,
$ s_{4}s_{5}s_{4}s_{1}s_{3}s_{2}s_{4}s_{1} $,
$ s_{4}s_{5}s_{4}s_{1}s_{3}s_{2}s_{4}s_{3} $,
$ s_{4}s_{5}s_{4}s_{3}s_{2}s_{4}s_{2}s_{1} $,
$ s_{4}s_{5}s_{6}s_{1}s_{3}s_{4}s_{3}s_{2} $,
$ s_{4}s_{5}s_{6}s_{1}s_{3}s_{4}s_{5}s_{3} $,
$ s_{4}s_{5}s_{6}s_{2}s_{4}s_{3}s_{2}s_{1} $,
$ s_{4}s_{5}s_{6}s_{2}s_{4}s_{5}s_{1}s_{3} $,
$ s_{4}s_{5}s_{6}s_{2}s_{4}s_{5}s_{3}s_{2} $,
$ s_{4}s_{5}s_{6}s_{3}s_{2}s_{4}s_{2}s_{1} $,
$ s_{4}s_{5}s_{6}s_{3}s_{2}s_{4}s_{3}s_{1} $,
$ s_{4}s_{5}s_{6}s_{3}s_{4}s_{1}s_{3}s_{1} $,
$ s_{4}s_{5}s_{6}s_{3}s_{4}s_{1}s_{3}s_{2} $,
$ s_{4}s_{5}s_{6}s_{3}s_{4}s_{3}s_{2}s_{1} $,
$ s_{4}s_{5}s_{6}s_{3}s_{4}s_{5}s_{1}s_{3} $,
$ s_{4}s_{5}s_{6}s_{3}s_{4}s_{5}s_{2}s_{1} $,
$ s_{4}s_{5}s_{6}s_{3}s_{4}s_{5}s_{3}s_{1} $,
$ s_{4}s_{5}s_{6}s_{4}s_{3}s_{2}s_{4}s_{1} $,
$ s_{4}s_{5}s_{6}s_{4}s_{3}s_{2}s_{4}s_{3} $,
$ s_{4}s_{5}s_{6}s_{5}s_{1}s_{3}s_{4}s_{1} $,
$ s_{4}s_{5}s_{6}s_{5}s_{1}s_{3}s_{4}s_{2} $,
$ s_{4}s_{5}s_{6}s_{5}s_{2}s_{4}s_{3}s_{1} $,
$ s_{4}s_{5}s_{6}s_{5}s_{3}s_{2}s_{4}s_{2} $,
$ s_{4}s_{5}s_{6}s_{5}s_{3}s_{2}s_{4}s_{3} $,
$ s_{4}s_{5}s_{6}s_{5}s_{4}s_{3}s_{2}s_{4} $,
$ s_{5}s_{1}s_{3}s_{2}s_{4}s_{3}s_{2}s_{1} $,
$ s_{5}s_{3}s_{4}s_{1}s_{3}s_{2}s_{4}s_{2} $,
$ s_{5}s_{3}s_{4}s_{1}s_{3}s_{2}s_{4}s_{3} $,
$ s_{5}s_{4}s_{1}s_{3}s_{2}s_{4}s_{3}s_{2} $,
$ s_{5}s_{4}s_{1}s_{3}s_{2}s_{4}s_{5}s_{1} $,
$ s_{5}s_{4}s_{1}s_{3}s_{2}s_{4}s_{5}s_{6} $,
$ s_{5}s_{4}s_{3}s_{2}s_{4}s_{1}s_{3}s_{1} $,
$ s_{5}s_{4}s_{3}s_{2}s_{4}s_{1}s_{3}s_{2} $,
$ s_{5}s_{4}s_{3}s_{2}s_{4}s_{5}s_{2}s_{4} $,
$ s_{5}s_{4}s_{3}s_{2}s_{4}s_{5}s_{3}s_{4} $,
$ s_{5}s_{4}s_{3}s_{2}s_{4}s_{5}s_{4}s_{2} $,
$ s_{5}s_{4}s_{3}s_{2}s_{4}s_{5}s_{4}s_{3} $,
$ s_{5}s_{4}s_{3}s_{2}s_{4}s_{5}s_{6}s_{2} $,
$ s_{5}s_{4}s_{3}s_{2}s_{4}s_{5}s_{6}s_{4} $,
$ s_{5}s_{4}s_{3}s_{2}s_{4}s_{5}s_{6}s_{5} $,
$ s_{5}s_{6}s_{1}s_{3}s_{2}s_{4}s_{3}s_{2} $,
$ s_{5}s_{6}s_{1}s_{3}s_{2}s_{4}s_{5}s_{2} $,
$ s_{5}s_{6}s_{1}s_{3}s_{4}s_{3}s_{2}s_{1} $,
$ s_{5}s_{6}s_{1}s_{3}s_{4}s_{5}s_{3}s_{1} $,
$ s_{5}s_{6}s_{1}s_{3}s_{4}s_{5}s_{3}s_{2} $,
$ s_{5}s_{6}s_{1}s_{3}s_{4}s_{5}s_{3}s_{4} $,
$ s_{5}s_{6}s_{1}s_{3}s_{4}s_{5}s_{4}s_{2} $,
$ s_{5}s_{6}s_{2}s_{4}s_{5}s_{1}s_{3}s_{4} $,
$ s_{5}s_{6}s_{2}s_{4}s_{5}s_{2}s_{4}s_{3} $,
$ s_{5}s_{6}s_{2}s_{4}s_{5}s_{3}s_{2}s_{1} $,
$ s_{5}s_{6}s_{2}s_{4}s_{5}s_{3}s_{4}s_{2} $,
$ s_{5}s_{6}s_{2}s_{4}s_{5}s_{4}s_{3}s_{1} $,
$ s_{5}s_{6}s_{3}s_{2}s_{4}s_{5}s_{2}s_{4} $,
$ s_{5}s_{6}s_{3}s_{2}s_{4}s_{5}s_{3}s_{1} $,
$ s_{5}s_{6}s_{3}s_{4}s_{1}s_{3}s_{2}s_{4} $,
$ s_{5}s_{6}s_{3}s_{4}s_{5}s_{1}s_{3}s_{4} $,
$ s_{5}s_{6}s_{3}s_{4}s_{5}s_{2}s_{4}s_{1} $,
$ s_{5}s_{6}s_{3}s_{4}s_{5}s_{2}s_{4}s_{2} $,
$ s_{5}s_{6}s_{3}s_{4}s_{5}s_{3}s_{4}s_{1} $,
$ s_{5}s_{6}s_{3}s_{4}s_{5}s_{3}s_{4}s_{3} $,
$ s_{5}s_{6}s_{4}s_{3}s_{2}s_{4}s_{1}s_{3} $,
$ s_{5}s_{6}s_{4}s_{5}s_{1}s_{3}s_{4}s_{3} $,
$ s_{5}s_{6}s_{4}s_{5}s_{2}s_{4}s_{1}s_{3} $,
$ s_{5}s_{6}s_{4}s_{5}s_{2}s_{4}s_{3}s_{2} $,
$ s_{5}s_{6}s_{4}s_{5}s_{3}s_{4}s_{1}s_{3} $,
$ s_{5}s_{6}s_{4}s_{5}s_{3}s_{4}s_{2}s_{1} $,
$ s_{5}s_{6}s_{4}s_{5}s_{3}s_{4}s_{3}s_{1} $,
$ s_{5}s_{6}s_{5}s_{4}s_{3}s_{2}s_{4}s_{5} $,
$ s_{6}s_{1}s_{3}s_{2}s_{4}s_{5}s_{2}s_{4} $,
$ s_{6}s_{1}s_{3}s_{2}s_{4}s_{5}s_{3}s_{1} $,
$ s_{6}s_{1}s_{3}s_{2}s_{4}s_{5}s_{3}s_{4} $,
$ s_{6}s_{1}s_{3}s_{4}s_{5}s_{1}s_{3}s_{4} $,
$ s_{6}s_{1}s_{3}s_{4}s_{5}s_{2}s_{4}s_{2} $,
$ s_{6}s_{1}s_{3}s_{4}s_{5}s_{3}s_{2}s_{4} $,
$ s_{6}s_{1}s_{3}s_{4}s_{5}s_{4}s_{3}s_{1} $,
$ s_{6}s_{2}s_{4}s_{5}s_{1}s_{3}s_{2}s_{4} $,
$ s_{6}s_{2}s_{4}s_{5}s_{1}s_{3}s_{4}s_{3} $,
$ s_{6}s_{2}s_{4}s_{5}s_{3}s_{2}s_{4}s_{3} $,
$ s_{6}s_{2}s_{4}s_{5}s_{3}s_{4}s_{1}s_{3} $,
$ s_{6}s_{2}s_{4}s_{5}s_{3}s_{4}s_{2}s_{1} $,
$ s_{6}s_{2}s_{4}s_{5}s_{3}s_{4}s_{3}s_{1} $,
$ s_{6}s_{3}s_{2}s_{4}s_{5}s_{1}s_{3}s_{4} $,
$ s_{6}s_{3}s_{2}s_{4}s_{5}s_{2}s_{4}s_{2} $,
$ s_{6}s_{3}s_{2}s_{4}s_{5}s_{3}s_{2}s_{1} $,
$ s_{6}s_{3}s_{2}s_{4}s_{5}s_{3}s_{4}s_{1} $,
$ s_{6}s_{3}s_{2}s_{4}s_{5}s_{3}s_{4}s_{3} $,
$ s_{6}s_{3}s_{2}s_{4}s_{5}s_{4}s_{2}s_{1} $,
$ s_{6}s_{3}s_{4}s_{1}s_{3}s_{2}s_{4}s_{5} $,
$ s_{6}s_{3}s_{4}s_{5}s_{1}s_{3}s_{4}s_{3} $,
$ s_{6}s_{3}s_{4}s_{5}s_{2}s_{4}s_{3}s_{1} $,
$ s_{6}s_{3}s_{4}s_{5}s_{3}s_{4}s_{1}s_{3} $,
$ s_{6}s_{3}s_{4}s_{5}s_{3}s_{4}s_{2}s_{1} $,
$ s_{6}s_{3}s_{4}s_{5}s_{3}s_{4}s_{3}s_{1} $,
$ s_{6}s_{3}s_{4}s_{5}s_{4}s_{1}s_{3}s_{2} $,
$ s_{6}s_{3}s_{4}s_{5}s_{4}s_{3}s_{2}s_{4} $,
$ s_{6}s_{4}s_{1}s_{3}s_{2}s_{4}s_{5}s_{4} $,
$ s_{6}s_{4}s_{3}s_{2}s_{4}s_{5}s_{1}s_{3} $,
$ s_{6}s_{4}s_{3}s_{2}s_{4}s_{5}s_{2}s_{4} $,
$ s_{6}s_{4}s_{5}s_{1}s_{3}s_{2}s_{4}s_{2} $,
$ s_{6}s_{4}s_{5}s_{1}s_{3}s_{2}s_{4}s_{3} $,
$ s_{6}s_{4}s_{5}s_{1}s_{3}s_{4}s_{2}s_{1} $,
$ s_{6}s_{4}s_{5}s_{3}s_{2}s_{4}s_{3}s_{2} $,
$ s_{6}s_{4}s_{5}s_{3}s_{4}s_{1}s_{3}s_{1} $,
$ s_{6}s_{4}s_{5}s_{3}s_{4}s_{3}s_{2}s_{1} $,
$ s_{6}s_{4}s_{5}s_{4}s_{1}s_{3}s_{2}s_{4} $,
$ s_{6}s_{5}s_{1}s_{3}s_{4}s_{1}s_{3}s_{2} $,
$ s_{6}s_{5}s_{3}s_{2}s_{4}s_{1}s_{3}s_{2} $,
$ s_{6}s_{5}s_{3}s_{4}s_{1}s_{3}s_{2}s_{1} $,
$ s_{6}s_{5}s_{4}s_{1}s_{3}s_{2}s_{4}s_{1} $,
$ s_{6}s_{5}s_{4}s_{1}s_{3}s_{2}s_{4}s_{2} $,
$ s_{6}s_{5}s_{4}s_{1}s_{3}s_{2}s_{4}s_{3} $,
$ s_{6}s_{5}s_{4}s_{3}s_{2}s_{4}s_{5}s_{1} $,
$ s_{6}s_{5}s_{4}s_{3}s_{2}s_{4}s_{5}s_{3} $,
$ s_{6}s_{5}s_{4}s_{3}s_{2}s_{4}s_{5}s_{4} $.
\end{quote}

\begin{Th}[The case $E_7$, $p=3$]\label{E7mod3}For the Chow ring of the exceptional group $\Ad E_7$ over $\mathbb{F}_3$,
we have
\begin{description}
\item[(1)] The projection $\pi^*:\CH^*(G/B;\mathbb{F}_3)\to \CH^*(\Ad E_7;\mathbb{F}_3)$ is given by
$$[\Sigma_w]\longmapsto \begin{cases}
1, & w=1\\
x_8,  & w\in \Pi(x_8),\\
-x_8,  & w\in \Pi(-x_8),\\
x_8^2,  & w\in \Pi(x_8^2),\\
-x_8^2,  & w\in \Pi(-x_8^2),\\
0, & \text{otherwise}
\end{cases}$$
\item[(2)] The comodule structure $\CH^*(G/B;\mathbb{F}_3)\to \CH^*(\Ad E_7;\mathbb{F}_3)\otimes \CH^*(G/B;\mathbb{F}_3)$ is given by
$$[\Sigma_w]\longmapsto 1\otimes [\Sigma_w]+
x_8\otimes \sum_{{w=u\odot v,} \atop {u\in \Pi(\pm x_8)}} \pm [\Sigma_{v}]
+x_8^2\otimes \sum_{{w=u\odot v,} \atop {u\in \Pi(\pm x_8^2)}} \pm [\Sigma_{v}]. $$
%
%(3) The coalgebra structure $\CH^*(F_4;\mathbb{F}_2)\to \CH^*(F_4;\mathbb{F}_2)\otimes \CH^*(F_4;\mathbb{F}_2)$ is given by
%$$x_6\longmapsto x_6\otimes 1+1\otimes x_6. $$
\end{description}
\end{Th}

\noindent\textbf{Remark.} It took roughly 40 minutes to get the result on author's personal computer.

\smallbreak

\noindent\textbf{Remark.} For the case $p=2$,
$\CH^*(\Ad E_7;\mathbb{F}_2)=\mathbb{F}_2[x_2,x_6,x_{10},x_{18}]\big/\big<x_2^2,x_6^2,x_{10}^2,x_{18}^2\big>$.
The author also computed the elements $w\in W$ such that $\pi^*[\Sigma_w]$ equals to $x_2$, $x_6$, $x_{10}$ or $x_{18}$
(but not yet for $x_2x_6,\ldots,x_2x_6x_{10}x_{18}$).
But it took too much pages to list them.

\long\def\crazy#1{%#1

\bibliographystyle{plain}
\bibliography{bibfile}

\begin{thebibliography}{1}

\bibitem{duan2014schubert}
Haibao Duan and Xuezhi Zhao.
\newblock Schubert calculus and the hopf algebra structures of exceptional lie
  groups.
\newblock In {\em Forum Mathematicum}, volume~26, pages 113--139. De Gruyter,
  2014.

\bibitem{grothendieck1958torsion}
Alexandre Grothendieck.
\newblock Torsion homologique et sections rationnelles.
\newblock {\em S{\'e}minaire Claude Chevalley}, 3:1--29, 1958.

\bibitem{hiller1982geometry}
Howard Hiller.
\newblock {\em Geometry of Coxeter groups}, volume~54.
\newblock Pitman Publishing, 1982.

\bibitem{kac1985torsion}
VG~Kac.
\newblock Torsion in cohomology of compact lie groups and chow rings of
  reductive algebraic groups.
\newblock {\em Inventiones mathematicae}, 80(1):69--79, 1985.

\bibitem{xiong2020comodule}
Rui Xiong.
\newblock Comodule structures, equivariant hopf structures, and generalized
  schubert polynomials, 2020.

\end{thebibliography}

\noindent

\hfill\rule{0.6\linewidth}{0.4pt}\hfill\mbox{}
\begin{center}
XIONG Rui, master stundent\smallbreak
\textit{\href{http://english.spbu.ru/}{Saint Petersburg State University}}
\textit{\href{https://math-cs.spbu.ru/en/}{Department of Mathematics and Computer Science}}
Saint Petersburg, 199178, Russia, Line 14th (Vasilyevsky Island),
\par
e-mail: \url{XiongRui_Math@126.com},\par
homepage: \url{www.cnblogs.com/XiongRuiMath}.
\end{center}

}

\crazy{

\section{Partial results}

\paragraph{The case $\Ad E_6$, $p=3$.}
Now
$\CH^*(\Ad E_6;\mathbb{F}_3)=\mathbb{F}_3[x_2,x_8]\big/\big<x_2^9,x_8^3\big>$.
The nonzero degree of $\CH^*(\Ad E_6;\mathbb{F}_3)$ is
$$\begin{array}{c|c|c|c|c|c|c|c|c}\hline\rule{0pc}{1.5pc}
\textrm{degree} & 2 & 4 & 6 & 8 & 10 & 12 & 14 & 16 \\ [2ex]\ru
\textrm{dimension} & 1 & 1 & 1 & 2 & 1 & 1 & 1 & 3  \\ \hline
\end{array}$$
Due to the restriction of the author's computer, in this case, the computation is not fully conducted.
But we can say something on the generators.

\begin{Prop}[The case $\Ad E_6$, $p=3$]\label{E6mod3}
The projection $\pi^*:\CH^*(G/B;\mathbb{F}_3)\to \CH^*(\Ad E_6;\mathbb{F}_3)$ satisfies
\begin{description}
\item[(1)] the set of $w\in W$ such that $\pi^*[\Sigma_w]=x_2$ is $\{s_1,s_5\}$;
\item[(2)] the set of $w\in W$ such that $\pi^*[\Sigma_w]=-x_2$ is $\{s_3,s_6\}$;
\item[(3)] the set of $w\in W$ such that $\pi^*[\Sigma_w]=x_8$ is
$$\left\{\begin{array}{c}
s_{1}s_{3}s_{4}s_{2},
s_{2}s_{4}s_{1}s_{3},
s_{3}s_{4}s_{5}s_{6},
s_{5}s_{1}s_{3}s_{4},
s_{5}s_{3}s_{2}s_{4},\\
s_{5}s_{4}s_{3}s_{1},
s_{6}s_{2}s_{4}s_{5},
s_{6}s_{5}s_{3}s_{4},
s_{6}s_{5}s_{4}s_{2}
\end{array}\right\};$$
\item[(4)] the set of $w\in W$ such that $\pi^*[\Sigma_w]=-x_8$ is
$$\left\{\begin{array}{c}
s_{1}s_{3}s_{4}s_{5},
s_{2}s_{4}s_{3}s_{1},
s_{2}s_{4}s_{5}s_{6},
s_{3}s_{4}s_{2}s_{1},
s_{4}s_{5}s_{3}s_{1},\\
s_{4}s_{5}s_{3}s_{2},
s_{4}s_{5}s_{6}s_{3},
s_{5}s_{6}s_{4}s_{2},
s_{6}s_{5}s_{4}s_{3}
\end{array}\right\}.$$
\end{description}
\end{Prop}

\paragraph{The case $\Ad E_7$, $p=2$.}
Now
$\CH^*(\Ad E_7;\mathbb{F}_2)=\mathbb{F}_2[x_2,x_6,x_{10},x_{18}]\big/\big<x_2^2,x_6^2,x_{10}^2,x_{18}^2\big>$.
The nonzero degree of $\CH^*(\Ad E_7;\mathbb{F}_2)$ is
$$\begin{array}{c|c|c|c|c|c|c|c|c|c|c|c|c}\hline\rule{0pc}{1.5pc}
\textrm{degree} & 2 & 6& 8 & 10 & 14 & 18 & 20 & 24 & 26 & 28& 32 & 36\\ [2ex]\ru
\textrm{dimension} & 1 & 1 & 1 & 1 & 1 & 2 & 1 & 1 & 1 & 1 & 1 & 1\\ \hline
\end{array}$$
The following result only computes for the generators of $\CH^*(\Ad E_7;\mathbb{F}_2)$.

\begin{Prop}[The case $E_7$, $p=2$, partially]
\label{E7mod2}The projection $\pi^*:\CH^*(G/B;\mathbb{F}_2)\to \CH^*(\Ad E_7;\mathbb{F}_2)$ satisfies
\begin{description}
\item[(1)] the set of $w\in W$ such that $\pi^*[\Sigma_w]=x_2$ is $\{s_2,s_5,s_7\}$;
\item[(2)] the set of $w\in W$ such that $\pi^*[\Sigma_w]=x_6$ is \{$s_2s_4s_3$,
$s_2s_4s_5$,
$s_3s_2s_4$,
$s_3s_4s_2$,
$s_3s_4s_5$,
$s_4s_3s_2$,
$s_4s_5s_2$,
$s_4s_5s_3$,
$s_5s_2s_4$,
$s_5s_3s_4$,
$s_5s_4s_2$,
$s_5s_4s_3$\};
\item[(3)] the set of $w\in W$ such that $\pi^*[\Sigma_w]=x_{10}$ is \{%
$ s_{2}s_{4}s_{5}s_{2}s_{4} $,
$ s_{2}s_{4}s_{5}s_{3}s_{4} $,
$ s_{2}s_{4}s_{5}s_{4}s_{2} $,
$ s_{2}s_{4}s_{5}s_{6}s_{3} $,
$ s_{2}s_{4}s_{5}s_{6}s_{4} $,
$ s_{2}s_{4}s_{5}s_{6}s_{5} $,
$ s_{2}s_{4}s_{5}s_{6}s_{7} $,
$ s_{3}s_{2}s_{4}s_{3}s_{2} $,
$ s_{3}s_{2}s_{4}s_{5}s_{2} $,
$ s_{3}s_{2}s_{4}s_{5}s_{3} $,
$ s_{3}s_{2}s_{4}s_{5}s_{4} $,
$ s_{3}s_{2}s_{4}s_{5}s_{6} $,
$ s_{3}s_{4}s_{5}s_{2}s_{4} $,
$ s_{3}s_{4}s_{5}s_{3}s_{4} $,
$ s_{3}s_{4}s_{5}s_{4}s_{3} $,
$ s_{3}s_{4}s_{5}s_{6}s_{2} $,
$ s_{3}s_{4}s_{5}s_{6}s_{4} $,
$ s_{3}s_{4}s_{5}s_{6}s_{5} $,
$ s_{3}s_{4}s_{5}s_{6}s_{7} $,
$ s_{4}s_{3}s_{2}s_{4}s_{2} $,
$ s_{4}s_{3}s_{2}s_{4}s_{3} $,
$ s_{4}s_{3}s_{2}s_{4}s_{5} $,
$ s_{4}s_{5}s_{2}s_{4}s_{2} $,
$ s_{4}s_{5}s_{2}s_{4}s_{3} $,
$ s_{4}s_{5}s_{3}s_{4}s_{2} $,
$ s_{4}s_{5}s_{3}s_{4}s_{3} $,
$ s_{4}s_{5}s_{4}s_{3}s_{2} $,
$ s_{4}s_{5}s_{6}s_{3}s_{2} $,
$ s_{4}s_{5}s_{6}s_{4}s_{2} $,
$ s_{4}s_{5}s_{6}s_{4}s_{3} $,
$ s_{4}s_{5}s_{6}s_{5}s_{2} $,
$ s_{4}s_{5}s_{6}s_{5}s_{3} $,
$ s_{4}s_{5}s_{6}s_{7}s_{2} $,
$ s_{4}s_{5}s_{6}s_{7}s_{3} $,
$ s_{5}s_{2}s_{4}s_{3}s_{2} $,
$ s_{5}s_{3}s_{4}s_{3}s_{2} $,
$ s_{5}s_{4}s_{3}s_{2}s_{4} $,
$ s_{5}s_{6}s_{2}s_{4}s_{3} $,
$ s_{5}s_{6}s_{3}s_{2}s_{4} $,
$ s_{5}s_{6}s_{3}s_{4}s_{2} $,
$ s_{5}s_{6}s_{4}s_{3}s_{2} $,
$ s_{5}s_{6}s_{5}s_{2}s_{4} $,
$ s_{5}s_{6}s_{5}s_{3}s_{4} $,
$ s_{5}s_{6}s_{5}s_{4}s_{2} $,
$ s_{5}s_{6}s_{5}s_{4}s_{3} $,
$ s_{5}s_{6}s_{7}s_{2}s_{4} $,
$ s_{5}s_{6}s_{7}s_{3}s_{4} $,
$ s_{5}s_{6}s_{7}s_{4}s_{2} $,
$ s_{5}s_{6}s_{7}s_{4}s_{3} $,
$ s_{6}s_{2}s_{4}s_{5}s_{3} $,
$ s_{6}s_{2}s_{4}s_{5}s_{4} $,
$ s_{6}s_{3}s_{2}s_{4}s_{5} $,
$ s_{6}s_{3}s_{4}s_{5}s_{2} $,
$ s_{6}s_{3}s_{4}s_{5}s_{4} $,
$ s_{6}s_{4}s_{5}s_{3}s_{2} $,
$ s_{6}s_{4}s_{5}s_{4}s_{2} $,
$ s_{6}s_{4}s_{5}s_{4}s_{3} $,
$ s_{6}s_{5}s_{2}s_{4}s_{3} $,
$ s_{6}s_{5}s_{3}s_{2}s_{4} $,
$ s_{6}s_{5}s_{3}s_{4}s_{2} $,
$ s_{6}s_{5}s_{4}s_{3}s_{2} $,
$ s_{6}s_{7}s_{2}s_{4}s_{5} $,
$ s_{6}s_{7}s_{3}s_{4}s_{5} $,
$ s_{6}s_{7}s_{4}s_{5}s_{2} $,
$ s_{6}s_{7}s_{4}s_{5}s_{3} $,
$ s_{6}s_{7}s_{5}s_{2}s_{4} $,
$ s_{6}s_{7}s_{5}s_{3}s_{4} $,
$ s_{6}s_{7}s_{5}s_{4}s_{2} $,
$ s_{6}s_{7}s_{5}s_{4}s_{3} $,
$ s_{7}s_{2}s_{4}s_{5}s_{6} $,
$ s_{7}s_{3}s_{4}s_{5}s_{6} $,
$ s_{7}s_{4}s_{5}s_{6}s_{2} $,
$ s_{7}s_{4}s_{5}s_{6}s_{3} $,
$ s_{7}s_{5}s_{6}s_{2}s_{4} $,
$ s_{7}s_{5}s_{6}s_{3}s_{4} $,
$ s_{7}s_{5}s_{6}s_{4}s_{2} $,
$ s_{7}s_{5}s_{6}s_{4}s_{3} $,
$ s_{7}s_{6}s_{2}s_{4}s_{5} $,
$ s_{7}s_{6}s_{3}s_{4}s_{5} $,
$ s_{7}s_{6}s_{4}s_{5}s_{2} $,
$ s_{7}s_{6}s_{4}s_{5}s_{3} $,
$ s_{7}s_{6}s_{5}s_{2}s_{4} $,
$ s_{7}s_{6}s_{5}s_{3}s_{4} $,
$ s_{7}s_{6}s_{5}s_{4}s_{2} $,
$ s_{7}s_{6}s_{5}s_{4}s_{3} $\};
\item[(4)] the set of $w\in W$ such that $\pi^*[\Sigma_w]=x_{18}$ is \{%
$ s_{1}s_{3}s_{2}s_{4}s_{5}s_{1}s_{3}s_{2}s_{4} $,
$ s_{1}s_{3}s_{2}s_{4}s_{5}s_{3}s_{2}s_{4}s_{1} $,
$ s_{1}s_{3}s_{2}s_{4}s_{5}s_{3}s_{2}s_{4}s_{3} $,
$ s_{1}s_{3}s_{2}s_{4}s_{5}s_{6}s_{1}s_{3}s_{2} $,
$ s_{1}s_{3}s_{2}s_{4}s_{5}s_{6}s_{1}s_{3}s_{4} $,
$ s_{1}s_{3}s_{2}s_{4}s_{5}s_{6}s_{2}s_{4}s_{3} $,
$ s_{1}s_{3}s_{2}s_{4}s_{5}s_{6}s_{2}s_{4}s_{5} $,
$ s_{1}s_{3}s_{2}s_{4}s_{5}s_{6}s_{3}s_{2}s_{1} $,
$ s_{1}s_{3}s_{2}s_{4}s_{5}s_{6}s_{3}s_{2}s_{4} $,
$ s_{1}s_{3}s_{2}s_{4}s_{5}s_{6}s_{3}s_{4}s_{1} $,
$ s_{1}s_{3}s_{2}s_{4}s_{5}s_{6}s_{3}s_{4}s_{3} $,
$ s_{1}s_{3}s_{2}s_{4}s_{5}s_{6}s_{4}s_{5}s_{2} $,
$ s_{1}s_{3}s_{2}s_{4}s_{5}s_{6}s_{5}s_{1}s_{3} $,
$ s_{1}s_{3}s_{2}s_{4}s_{5}s_{6}s_{5}s_{2}s_{4} $,
$ s_{1}s_{3}s_{2}s_{4}s_{5}s_{6}s_{5}s_{3}s_{1} $,
$ s_{1}s_{3}s_{2}s_{4}s_{5}s_{6}s_{5}s_{3}s_{4} $,
$ s_{1}s_{3}s_{2}s_{4}s_{5}s_{6}s_{5}s_{4}s_{2} $,
$ s_{1}s_{3}s_{2}s_{4}s_{5}s_{6}s_{5}s_{4}s_{3} $,
$ s_{1}s_{3}s_{2}s_{4}s_{5}s_{6}s_{7}s_{1}s_{3} $,
$ s_{1}s_{3}s_{2}s_{4}s_{5}s_{6}s_{7}s_{3}s_{1} $,
$ s_{1}s_{3}s_{2}s_{4}s_{5}s_{6}s_{7}s_{3}s_{4} $,
$ s_{1}s_{3}s_{2}s_{4}s_{5}s_{6}s_{7}s_{4}s_{3} $,
$ s_{1}s_{3}s_{4}s_{5}s_{1}s_{3}s_{2}s_{4}s_{2} $,
$ s_{1}s_{3}s_{4}s_{5}s_{3}s_{2}s_{4}s_{2}s_{1} $,
$ s_{1}s_{3}s_{4}s_{5}s_{3}s_{2}s_{4}s_{3}s_{2} $,
$ s_{1}s_{3}s_{4}s_{5}s_{4}s_{1}s_{3}s_{2}s_{4} $,
$ s_{1}s_{3}s_{4}s_{5}s_{4}s_{3}s_{2}s_{4}s_{1} $,
$ s_{1}s_{3}s_{4}s_{5}s_{4}s_{3}s_{2}s_{4}s_{3} $,
$ s_{1}s_{3}s_{4}s_{5}s_{6}s_{1}s_{3}s_{4}s_{2} $,
$ s_{1}s_{3}s_{4}s_{5}s_{6}s_{1}s_{3}s_{4}s_{3} $,
$ s_{1}s_{3}s_{4}s_{5}s_{6}s_{1}s_{3}s_{4}s_{5} $,
$ s_{1}s_{3}s_{4}s_{5}s_{6}s_{2}s_{4}s_{5}s_{3} $,
$ s_{1}s_{3}s_{4}s_{5}s_{6}s_{2}s_{4}s_{5}s_{4} $,
$ s_{1}s_{3}s_{4}s_{5}s_{6}s_{3}s_{2}s_{4}s_{2} $,
$ s_{1}s_{3}s_{4}s_{5}s_{6}s_{3}s_{2}s_{4}s_{5} $,
$ s_{1}s_{3}s_{4}s_{5}s_{6}s_{3}s_{4}s_{2}s_{1} $,
$ s_{1}s_{3}s_{4}s_{5}s_{6}s_{3}s_{4}s_{3}s_{1} $,
$ s_{1}s_{3}s_{4}s_{5}s_{6}s_{3}s_{4}s_{3}s_{2} $,
$ s_{1}s_{3}s_{4}s_{5}s_{6}s_{3}s_{4}s_{5}s_{1} $,
$ s_{1}s_{3}s_{4}s_{5}s_{6}s_{3}s_{4}s_{5}s_{3} $,
$ s_{1}s_{3}s_{4}s_{5}s_{6}s_{3}s_{4}s_{5}s_{4} $,
$ s_{1}s_{3}s_{4}s_{5}s_{6}s_{4}s_{1}s_{3}s_{2} $,
$ s_{1}s_{3}s_{4}s_{5}s_{6}s_{4}s_{3}s_{2}s_{1} $,
$ s_{1}s_{3}s_{4}s_{5}s_{6}s_{4}s_{3}s_{2}s_{4} $,
$ s_{1}s_{3}s_{4}s_{5}s_{6}s_{4}s_{5}s_{1}s_{3} $,
$ s_{1}s_{3}s_{4}s_{5}s_{6}s_{4}s_{5}s_{2}s_{4} $,
$ s_{1}s_{3}s_{4}s_{5}s_{6}s_{4}s_{5}s_{3}s_{1} $,
$ s_{1}s_{3}s_{4}s_{5}s_{6}s_{4}s_{5}s_{3}s_{4} $,
$ s_{1}s_{3}s_{4}s_{5}s_{6}s_{5}s_{2}s_{4}s_{2} $,
$ s_{1}s_{3}s_{4}s_{5}s_{6}s_{5}s_{3}s_{4}s_{2} $,
$ s_{1}s_{3}s_{4}s_{5}s_{6}s_{5}s_{4}s_{3}s_{2} $,
$ s_{1}s_{3}s_{4}s_{5}s_{6}s_{7}s_{1}s_{3}s_{4} $,
$ s_{1}s_{3}s_{4}s_{5}s_{6}s_{7}s_{2}s_{4}s_{3} $,
$ s_{1}s_{3}s_{4}s_{5}s_{6}s_{7}s_{2}s_{4}s_{5} $,
$ s_{1}s_{3}s_{4}s_{5}s_{6}s_{7}s_{3}s_{2}s_{4} $,
$ s_{1}s_{3}s_{4}s_{5}s_{6}s_{7}s_{3}s_{4}s_{1} $,
$ s_{1}s_{3}s_{4}s_{5}s_{6}s_{7}s_{3}s_{4}s_{2} $,
$ s_{1}s_{3}s_{4}s_{5}s_{6}s_{7}s_{3}s_{4}s_{5} $,
$ s_{1}s_{3}s_{4}s_{5}s_{6}s_{7}s_{4}s_{1}s_{3} $,
$ s_{1}s_{3}s_{4}s_{5}s_{6}s_{7}s_{4}s_{3}s_{1} $,
$ s_{1}s_{3}s_{4}s_{5}s_{6}s_{7}s_{4}s_{3}s_{2} $,
$ s_{1}s_{3}s_{4}s_{5}s_{6}s_{7}s_{4}s_{5}s_{2} $,
$ s_{1}s_{3}s_{4}s_{5}s_{6}s_{7}s_{4}s_{5}s_{3} $,
$ s_{1}s_{3}s_{4}s_{5}s_{6}s_{7}s_{5}s_{2}s_{4} $,
$ s_{1}s_{3}s_{4}s_{5}s_{6}s_{7}s_{5}s_{3}s_{4} $,
$ s_{1}s_{3}s_{4}s_{5}s_{6}s_{7}s_{5}s_{4}s_{2} $,
$ s_{1}s_{3}s_{4}s_{5}s_{6}s_{7}s_{5}s_{4}s_{3} $,
$ s_{2}s_{4}s_{5}s_{1}s_{3}s_{2}s_{4}s_{3}s_{1} $,
$ s_{2}s_{4}s_{5}s_{3}s_{2}s_{4}s_{1}s_{3}s_{1} $,
$ s_{2}s_{4}s_{5}s_{4}s_{1}s_{3}s_{2}s_{4}s_{1} $,
$ s_{2}s_{4}s_{5}s_{4}s_{3}s_{2}s_{4}s_{1}s_{3} $,
$ s_{2}s_{4}s_{5}s_{4}s_{3}s_{2}s_{4}s_{3}s_{1} $,
$ s_{2}s_{4}s_{5}s_{6}s_{1}s_{3}s_{2}s_{4}s_{3} $,
$ s_{2}s_{4}s_{5}s_{6}s_{1}s_{3}s_{4}s_{5}s_{1} $,
$ s_{2}s_{4}s_{5}s_{6}s_{1}s_{3}s_{4}s_{5}s_{3} $,
$ s_{2}s_{4}s_{5}s_{6}s_{1}s_{3}s_{4}s_{5}s_{4} $,
$ s_{2}s_{4}s_{5}s_{6}s_{2}s_{4}s_{1}s_{3}s_{2} $,
$ s_{2}s_{4}s_{5}s_{6}s_{2}s_{4}s_{3}s_{2}s_{1} $,
$ s_{2}s_{4}s_{5}s_{6}s_{3}s_{2}s_{4}s_{3}s_{1} $,
$ s_{2}s_{4}s_{5}s_{6}s_{3}s_{4}s_{5}s_{1}s_{3} $,
$ s_{2}s_{4}s_{5}s_{6}s_{3}s_{4}s_{5}s_{4}s_{1} $,
$ s_{2}s_{4}s_{5}s_{6}s_{4}s_{1}s_{3}s_{2}s_{4} $,
$ s_{2}s_{4}s_{5}s_{6}s_{4}s_{3}s_{2}s_{4}s_{1} $,
$ s_{2}s_{4}s_{5}s_{6}s_{4}s_{5}s_{1}s_{3}s_{4} $,
$ s_{2}s_{4}s_{5}s_{6}s_{4}s_{5}s_{3}s_{4}s_{1} $,
$ s_{2}s_{4}s_{5}s_{6}s_{4}s_{5}s_{3}s_{4}s_{3} $,
$ s_{2}s_{4}s_{5}s_{6}s_{5}s_{1}s_{3}s_{4}s_{3} $,
$ s_{2}s_{4}s_{5}s_{6}s_{5}s_{2}s_{4}s_{1}s_{3} $,
$ s_{2}s_{4}s_{5}s_{6}s_{5}s_{2}s_{4}s_{3}s_{1} $,
$ s_{2}s_{4}s_{5}s_{6}s_{5}s_{3}s_{2}s_{4}s_{3} $,
$ s_{2}s_{4}s_{5}s_{6}s_{5}s_{3}s_{4}s_{3}s_{1} $,
$ s_{2}s_{4}s_{5}s_{6}s_{5}s_{4}s_{1}s_{3}s_{2} $,
$ s_{2}s_{4}s_{5}s_{6}s_{5}s_{4}s_{3}s_{2}s_{1} $,
$ s_{2}s_{4}s_{5}s_{6}s_{5}s_{4}s_{3}s_{2}s_{4} $,
$ s_{2}s_{4}s_{5}s_{6}s_{7}s_{1}s_{3}s_{4}s_{1} $,
$ s_{2}s_{4}s_{5}s_{6}s_{7}s_{1}s_{3}s_{4}s_{3} $,
$ s_{2}s_{4}s_{5}s_{6}s_{7}s_{1}s_{3}s_{4}s_{5} $,
$ s_{2}s_{4}s_{5}s_{6}s_{7}s_{2}s_{4}s_{1}s_{3} $,
$ s_{2}s_{4}s_{5}s_{6}s_{7}s_{2}s_{4}s_{3}s_{1} $,
$ s_{2}s_{4}s_{5}s_{6}s_{7}s_{2}s_{4}s_{3}s_{2} $,
$ s_{2}s_{4}s_{5}s_{6}s_{7}s_{2}s_{4}s_{5}s_{2} $,
$ s_{2}s_{4}s_{5}s_{6}s_{7}s_{2}s_{4}s_{5}s_{4} $,
$ s_{2}s_{4}s_{5}s_{6}s_{7}s_{2}s_{4}s_{5}s_{6} $,
$ s_{2}s_{4}s_{5}s_{6}s_{7}s_{3}s_{4}s_{1}s_{3} $,
$ s_{2}s_{4}s_{5}s_{6}s_{7}s_{3}s_{4}s_{5}s_{1} $,
$ s_{2}s_{4}s_{5}s_{6}s_{7}s_{3}s_{4}s_{5}s_{3} $,
$ s_{2}s_{4}s_{5}s_{6}s_{7}s_{3}s_{4}s_{5}s_{4} $,
$ s_{2}s_{4}s_{5}s_{6}s_{7}s_{3}s_{4}s_{5}s_{6} $,
$ s_{2}s_{4}s_{5}s_{6}s_{7}s_{4}s_{1}s_{3}s_{2} $,
$ s_{2}s_{4}s_{5}s_{6}s_{7}s_{4}s_{3}s_{2}s_{1} $,
$ s_{2}s_{4}s_{5}s_{6}s_{7}s_{4}s_{5}s_{1}s_{3} $,
$ s_{2}s_{4}s_{5}s_{6}s_{7}s_{4}s_{5}s_{3}s_{1} $,
$ s_{2}s_{4}s_{5}s_{6}s_{7}s_{4}s_{5}s_{6}s_{2} $,
$ s_{2}s_{4}s_{5}s_{6}s_{7}s_{4}s_{5}s_{6}s_{3} $,
$ s_{2}s_{4}s_{5}s_{6}s_{7}s_{5}s_{1}s_{3}s_{4} $,
$ s_{2}s_{4}s_{5}s_{6}s_{7}s_{5}s_{3}s_{4}s_{1} $,
$ s_{2}s_{4}s_{5}s_{6}s_{7}s_{5}s_{4}s_{1}s_{3} $,
$ s_{2}s_{4}s_{5}s_{6}s_{7}s_{5}s_{4}s_{3}s_{1} $,
$ s_{2}s_{4}s_{5}s_{6}s_{7}s_{5}s_{6}s_{2}s_{4} $,
$ s_{2}s_{4}s_{5}s_{6}s_{7}s_{5}s_{6}s_{3}s_{4} $,
$ s_{2}s_{4}s_{5}s_{6}s_{7}s_{5}s_{6}s_{4}s_{2} $,
$ s_{2}s_{4}s_{5}s_{6}s_{7}s_{5}s_{6}s_{4}s_{3} $,
$ s_{3}s_{2}s_{4}s_{5}s_{1}s_{3}s_{2}s_{4}s_{3} $,
$ s_{3}s_{2}s_{4}s_{5}s_{3}s_{2}s_{4}s_{1}s_{3} $,
$ s_{3}s_{2}s_{4}s_{5}s_{3}s_{2}s_{4}s_{3}s_{1} $,
$ s_{3}s_{2}s_{4}s_{5}s_{4}s_{1}s_{3}s_{2}s_{4} $,
$ s_{3}s_{2}s_{4}s_{5}s_{4}s_{3}s_{2}s_{4}s_{1} $,
$ s_{3}s_{2}s_{4}s_{5}s_{4}s_{3}s_{2}s_{4}s_{3} $,
$ s_{3}s_{2}s_{4}s_{5}s_{6}s_{1}s_{3}s_{2}s_{1} $,
$ s_{3}s_{2}s_{4}s_{5}s_{6}s_{1}s_{3}s_{4}s_{5} $,
$ s_{3}s_{2}s_{4}s_{5}s_{6}s_{2}s_{4}s_{3}s_{1} $,
$ s_{3}s_{2}s_{4}s_{5}s_{6}s_{2}s_{4}s_{5}s_{1} $,
$ s_{3}s_{2}s_{4}s_{5}s_{6}s_{3}s_{2}s_{4}s_{1} $,
$ s_{3}s_{2}s_{4}s_{5}s_{6}s_{3}s_{4}s_{1}s_{3} $,
$ s_{3}s_{2}s_{4}s_{5}s_{6}s_{3}s_{4}s_{5}s_{1} $,
$ s_{3}s_{2}s_{4}s_{5}s_{6}s_{3}s_{4}s_{5}s_{3} $,
$ s_{3}s_{2}s_{4}s_{5}s_{6}s_{4}s_{1}s_{3}s_{2} $,
$ s_{3}s_{2}s_{4}s_{5}s_{6}s_{4}s_{3}s_{2}s_{1} $,
$ s_{3}s_{2}s_{4}s_{5}s_{6}s_{4}s_{3}s_{2}s_{4} $,
$ s_{3}s_{2}s_{4}s_{5}s_{6}s_{4}s_{5}s_{1}s_{3} $,
$ s_{3}s_{2}s_{4}s_{5}s_{6}s_{4}s_{5}s_{2}s_{1} $,
$ s_{3}s_{2}s_{4}s_{5}s_{6}s_{4}s_{5}s_{2}s_{4} $,
$ s_{3}s_{2}s_{4}s_{5}s_{6}s_{4}s_{5}s_{3}s_{1} $,
$ s_{3}s_{2}s_{4}s_{5}s_{6}s_{4}s_{5}s_{3}s_{2} $,
$ s_{3}s_{2}s_{4}s_{5}s_{6}s_{5}s_{1}s_{3}s_{1} $,
$ s_{3}s_{2}s_{4}s_{5}s_{6}s_{5}s_{1}s_{3}s_{2} $,
$ s_{3}s_{2}s_{4}s_{5}s_{6}s_{5}s_{2}s_{4}s_{1} $,
$ s_{3}s_{2}s_{4}s_{5}s_{6}s_{5}s_{3}s_{2}s_{1} $,
$ s_{3}s_{2}s_{4}s_{5}s_{6}s_{5}s_{3}s_{2}s_{4} $,
$ s_{3}s_{2}s_{4}s_{5}s_{6}s_{5}s_{3}s_{4}s_{1} $,
$ s_{3}s_{2}s_{4}s_{5}s_{6}s_{5}s_{3}s_{4}s_{3} $,
$ s_{3}s_{2}s_{4}s_{5}s_{6}s_{5}s_{4}s_{2}s_{1} $,
$ s_{3}s_{2}s_{4}s_{5}s_{6}s_{5}s_{4}s_{3}s_{1} $,
$ s_{3}s_{2}s_{4}s_{5}s_{6}s_{7}s_{1}s_{3}s_{1} $,
$ s_{3}s_{2}s_{4}s_{5}s_{6}s_{7}s_{1}s_{3}s_{2} $,
$ s_{3}s_{2}s_{4}s_{5}s_{6}s_{7}s_{1}s_{3}s_{4} $,
$ s_{3}s_{2}s_{4}s_{5}s_{6}s_{7}s_{3}s_{2}s_{1} $,
$ s_{3}s_{2}s_{4}s_{5}s_{6}s_{7}s_{4}s_{1}s_{3} $,
$ s_{3}s_{4}s_{1}s_{3}s_{2}s_{4}s_{5}s_{2}s_{4} $,
$ s_{3}s_{4}s_{1}s_{3}s_{2}s_{4}s_{5}s_{3}s_{2} $,
$ s_{3}s_{4}s_{1}s_{3}s_{2}s_{4}s_{5}s_{3}s_{4} $,
$ s_{3}s_{4}s_{1}s_{3}s_{2}s_{4}s_{5}s_{4}s_{2} $,
$ s_{3}s_{4}s_{1}s_{3}s_{2}s_{4}s_{5}s_{4}s_{3} $,
$ s_{3}s_{4}s_{1}s_{3}s_{2}s_{4}s_{5}s_{6}s_{2} $,
$ s_{3}s_{4}s_{1}s_{3}s_{2}s_{4}s_{5}s_{6}s_{3} $,
$ s_{3}s_{4}s_{1}s_{3}s_{2}s_{4}s_{5}s_{6}s_{5} $,
$ s_{3}s_{4}s_{1}s_{3}s_{2}s_{4}s_{5}s_{6}s_{7} $,
$ s_{3}s_{4}s_{5}s_{1}s_{3}s_{2}s_{4}s_{3}s_{2} $,
$ s_{3}s_{4}s_{5}s_{3}s_{2}s_{4}s_{1}s_{3}s_{2} $,
$ s_{3}s_{4}s_{5}s_{3}s_{2}s_{4}s_{3}s_{2}s_{1} $,
$ s_{3}s_{4}s_{5}s_{4}s_{1}s_{3}s_{2}s_{4}s_{2} $,
$ s_{3}s_{4}s_{5}s_{4}s_{1}s_{3}s_{2}s_{4}s_{3} $,
$ s_{3}s_{4}s_{5}s_{4}s_{3}s_{2}s_{4}s_{1}s_{3} $,
$ s_{3}s_{4}s_{5}s_{4}s_{3}s_{2}s_{4}s_{2}s_{1} $,
$ s_{3}s_{4}s_{5}s_{4}s_{3}s_{2}s_{4}s_{3}s_{1} $,
$ s_{3}s_{4}s_{5}s_{6}s_{1}s_{3}s_{4}s_{5}s_{2} $,
$ s_{3}s_{4}s_{5}s_{6}s_{1}s_{3}s_{4}s_{5}s_{3} $,
$ s_{3}s_{4}s_{5}s_{6}s_{2}s_{4}s_{5}s_{3}s_{1} $,
$ s_{3}s_{4}s_{5}s_{6}s_{2}s_{4}s_{5}s_{4}s_{1} $,
$ s_{3}s_{4}s_{5}s_{6}s_{3}s_{2}s_{4}s_{2}s_{1} $,
$ s_{3}s_{4}s_{5}s_{6}s_{3}s_{2}s_{4}s_{5}s_{1} $,
$ s_{3}s_{4}s_{5}s_{6}s_{3}s_{2}s_{4}s_{5}s_{4} $,
$ s_{3}s_{4}s_{5}s_{6}s_{3}s_{4}s_{1}s_{3}s_{1} $,
$ s_{3}s_{4}s_{5}s_{6}s_{3}s_{4}s_{1}s_{3}s_{2} $,
$ s_{3}s_{4}s_{5}s_{6}s_{3}s_{4}s_{5}s_{1}s_{3} $,
$ s_{3}s_{4}s_{5}s_{6}s_{3}s_{4}s_{5}s_{2}s_{1} $,
$ s_{3}s_{4}s_{5}s_{6}s_{3}s_{4}s_{5}s_{3}s_{1} $,
$ s_{3}s_{4}s_{5}s_{6}s_{3}s_{4}s_{5}s_{3}s_{2} $,
$ s_{3}s_{4}s_{5}s_{6}s_{3}s_{4}s_{5}s_{4}s_{1} $,
$ s_{3}s_{4}s_{5}s_{6}s_{3}s_{4}s_{5}s_{4}s_{2} $,
$ s_{3}s_{4}s_{5}s_{6}s_{3}s_{4}s_{5}s_{4}s_{3} $,
$ s_{3}s_{4}s_{5}s_{6}s_{4}s_{1}s_{3}s_{2}s_{1} $,
$ s_{3}s_{4}s_{5}s_{6}s_{4}s_{3}s_{2}s_{4}s_{1} $,
$ s_{3}s_{4}s_{5}s_{6}s_{4}s_{5}s_{1}s_{3}s_{1} $,
$ s_{3}s_{4}s_{5}s_{6}s_{4}s_{5}s_{1}s_{3}s_{2} $,
$ s_{3}s_{4}s_{5}s_{6}s_{4}s_{5}s_{2}s_{4}s_{1} $,
$ s_{3}s_{4}s_{5}s_{6}s_{4}s_{5}s_{3}s_{2}s_{1} $,
$ s_{3}s_{4}s_{5}s_{6}s_{4}s_{5}s_{3}s_{2}s_{4} $,
$ s_{3}s_{4}s_{5}s_{6}s_{4}s_{5}s_{3}s_{4}s_{1} $,
$ s_{3}s_{4}s_{5}s_{6}s_{4}s_{5}s_{3}s_{4}s_{3} $,
$ s_{3}s_{4}s_{5}s_{6}s_{5}s_{2}s_{4}s_{2}s_{1} $,
$ s_{3}s_{4}s_{5}s_{6}s_{5}s_{3}s_{2}s_{4}s_{2} $,
$ s_{3}s_{4}s_{5}s_{6}s_{5}s_{3}s_{4}s_{2}s_{1} $,
$ s_{3}s_{4}s_{5}s_{6}s_{5}s_{3}s_{4}s_{3}s_{2} $,
$ s_{3}s_{4}s_{5}s_{6}s_{5}s_{4}s_{3}s_{2}s_{1} $,
$ s_{3}s_{4}s_{5}s_{6}s_{7}s_{1}s_{3}s_{4}s_{2} $,
$ s_{3}s_{4}s_{5}s_{6}s_{7}s_{1}s_{3}s_{4}s_{3} $,
$ s_{3}s_{4}s_{5}s_{6}s_{7}s_{2}s_{4}s_{3}s_{1} $,
$ s_{3}s_{4}s_{5}s_{6}s_{7}s_{2}s_{4}s_{5}s_{1} $,
$ s_{3}s_{4}s_{5}s_{6}s_{7}s_{2}s_{4}s_{5}s_{4} $,
$ s_{3}s_{4}s_{5}s_{6}s_{7}s_{2}s_{4}s_{5}s_{6} $,
$ s_{3}s_{4}s_{5}s_{6}s_{7}s_{3}s_{2}s_{4}s_{1} $,
$ s_{3}s_{4}s_{5}s_{6}s_{7}s_{3}s_{2}s_{4}s_{2} $,
$ s_{3}s_{4}s_{5}s_{6}s_{7}s_{3}s_{2}s_{4}s_{5} $,
$ s_{3}s_{4}s_{5}s_{6}s_{7}s_{3}s_{4}s_{1}s_{3} $,
$ s_{3}s_{4}s_{5}s_{6}s_{7}s_{3}s_{4}s_{5}s_{1} $,
$ s_{3}s_{4}s_{5}s_{6}s_{7}s_{3}s_{4}s_{5}s_{4} $,
$ s_{3}s_{4}s_{5}s_{6}s_{7}s_{3}s_{4}s_{5}s_{6} $,
$ s_{3}s_{4}s_{5}s_{6}s_{7}s_{4}s_{1}s_{3}s_{1} $,
$ s_{3}s_{4}s_{5}s_{6}s_{7}s_{4}s_{1}s_{3}s_{2} $,
$ s_{3}s_{4}s_{5}s_{6}s_{7}s_{4}s_{5}s_{2}s_{1} $,
$ s_{3}s_{4}s_{5}s_{6}s_{7}s_{4}s_{5}s_{3}s_{1} $,
$ s_{3}s_{4}s_{5}s_{6}s_{7}s_{4}s_{5}s_{3}s_{2} $,
$ s_{3}s_{4}s_{5}s_{6}s_{7}s_{4}s_{5}s_{6}s_{2} $,
$ s_{3}s_{4}s_{5}s_{6}s_{7}s_{4}s_{5}s_{6}s_{3} $,
$ s_{3}s_{4}s_{5}s_{6}s_{7}s_{5}s_{2}s_{4}s_{1} $,
$ s_{3}s_{4}s_{5}s_{6}s_{7}s_{5}s_{2}s_{4}s_{2} $,
$ s_{3}s_{4}s_{5}s_{6}s_{7}s_{5}s_{3}s_{2}s_{4} $,
$ s_{3}s_{4}s_{5}s_{6}s_{7}s_{5}s_{3}s_{4}s_{1} $,
$ s_{3}s_{4}s_{5}s_{6}s_{7}s_{5}s_{3}s_{4}s_{3} $,
$ s_{3}s_{4}s_{5}s_{6}s_{7}s_{5}s_{4}s_{2}s_{1} $,
$ s_{3}s_{4}s_{5}s_{6}s_{7}s_{5}s_{4}s_{3}s_{1} $,
$ s_{3}s_{4}s_{5}s_{6}s_{7}s_{5}s_{4}s_{3}s_{2} $,
$ s_{3}s_{4}s_{5}s_{6}s_{7}s_{5}s_{6}s_{2}s_{4} $,
$ s_{3}s_{4}s_{5}s_{6}s_{7}s_{5}s_{6}s_{3}s_{4} $,
$ s_{3}s_{4}s_{5}s_{6}s_{7}s_{5}s_{6}s_{4}s_{2} $,
$ s_{3}s_{4}s_{5}s_{6}s_{7}s_{5}s_{6}s_{4}s_{3} $,
$ s_{4}s_{1}s_{3}s_{2}s_{4}s_{5}s_{2}s_{4}s_{1} $,
$ s_{4}s_{1}s_{3}s_{2}s_{4}s_{5}s_{2}s_{4}s_{2} $,
$ s_{4}s_{1}s_{3}s_{2}s_{4}s_{5}s_{2}s_{4}s_{3} $,
$ s_{4}s_{1}s_{3}s_{2}s_{4}s_{5}s_{3}s_{2}s_{1} $,
$ s_{4}s_{1}s_{3}s_{2}s_{4}s_{5}s_{3}s_{2}s_{4} $,
$ s_{4}s_{1}s_{3}s_{2}s_{4}s_{5}s_{3}s_{4}s_{1} $,
$ s_{4}s_{1}s_{3}s_{2}s_{4}s_{5}s_{3}s_{4}s_{3} $,
$ s_{4}s_{1}s_{3}s_{2}s_{4}s_{5}s_{4}s_{2}s_{1} $,
$ s_{4}s_{1}s_{3}s_{2}s_{4}s_{5}s_{4}s_{3}s_{1} $,
$ s_{4}s_{1}s_{3}s_{2}s_{4}s_{5}s_{4}s_{3}s_{2} $,
$ s_{4}s_{1}s_{3}s_{2}s_{4}s_{5}s_{6}s_{2}s_{1} $,
$ s_{4}s_{1}s_{3}s_{2}s_{4}s_{5}s_{6}s_{3}s_{1} $,
$ s_{4}s_{1}s_{3}s_{2}s_{4}s_{5}s_{6}s_{3}s_{4} $,
$ s_{4}s_{1}s_{3}s_{2}s_{4}s_{5}s_{6}s_{4}s_{2} $,
$ s_{4}s_{1}s_{3}s_{2}s_{4}s_{5}s_{6}s_{4}s_{5} $,
$ s_{4}s_{1}s_{3}s_{2}s_{4}s_{5}s_{6}s_{5}s_{1} $,
$ s_{4}s_{1}s_{3}s_{2}s_{4}s_{5}s_{6}s_{7}s_{1} $,
$ s_{4}s_{1}s_{3}s_{2}s_{4}s_{5}s_{6}s_{7}s_{4} $,
$ s_{4}s_{3}s_{2}s_{4}s_{5}s_{1}s_{3}s_{2}s_{1} $,
$ s_{4}s_{3}s_{2}s_{4}s_{5}s_{2}s_{4}s_{1}s_{3} $,
$ s_{4}s_{3}s_{2}s_{4}s_{5}s_{2}s_{4}s_{3}s_{1} $,
$ s_{4}s_{3}s_{2}s_{4}s_{5}s_{2}s_{4}s_{3}s_{2} $,
$ s_{4}s_{3}s_{2}s_{4}s_{5}s_{3}s_{4}s_{1}s_{3} $,
$ s_{4}s_{3}s_{2}s_{4}s_{5}s_{4}s_{1}s_{3}s_{1} $,
$ s_{4}s_{3}s_{2}s_{4}s_{5}s_{4}s_{1}s_{3}s_{2} $,
$ s_{4}s_{3}s_{2}s_{4}s_{5}s_{4}s_{3}s_{2}s_{4} $,
$ s_{4}s_{3}s_{2}s_{4}s_{5}s_{6}s_{1}s_{3}s_{1} $,
$ s_{4}s_{3}s_{2}s_{4}s_{5}s_{6}s_{1}s_{3}s_{2} $,
$ s_{4}s_{3}s_{2}s_{4}s_{5}s_{6}s_{3}s_{4}s_{1} $,
$ s_{4}s_{3}s_{2}s_{4}s_{5}s_{6}s_{3}s_{4}s_{3} $,
$ s_{4}s_{3}s_{2}s_{4}s_{5}s_{6}s_{4}s_{2}s_{1} $,
$ s_{4}s_{3}s_{2}s_{4}s_{5}s_{6}s_{4}s_{3}s_{2} $,
$ s_{4}s_{3}s_{2}s_{4}s_{5}s_{6}s_{4}s_{5}s_{1} $,
$ s_{4}s_{3}s_{2}s_{4}s_{5}s_{6}s_{4}s_{5}s_{3} $,
$ s_{4}s_{3}s_{2}s_{4}s_{5}s_{6}s_{4}s_{5}s_{4} $,
$ s_{4}s_{3}s_{2}s_{4}s_{5}s_{6}s_{5}s_{1}s_{3} $,
$ s_{4}s_{3}s_{2}s_{4}s_{5}s_{6}s_{5}s_{2}s_{1} $,
$ s_{4}s_{3}s_{2}s_{4}s_{5}s_{6}s_{5}s_{3}s_{1} $,
$ s_{4}s_{3}s_{2}s_{4}s_{5}s_{6}s_{5}s_{3}s_{4} $,
$ s_{4}s_{3}s_{2}s_{4}s_{5}s_{6}s_{5}s_{4}s_{2} $,
$ s_{4}s_{3}s_{2}s_{4}s_{5}s_{6}s_{7}s_{1}s_{3} $,
$ s_{4}s_{3}s_{2}s_{4}s_{5}s_{6}s_{7}s_{2}s_{1} $,
$ s_{4}s_{3}s_{2}s_{4}s_{5}s_{6}s_{7}s_{3}s_{1} $,
$ s_{4}s_{3}s_{2}s_{4}s_{5}s_{6}s_{7}s_{3}s_{4} $,
$ s_{4}s_{3}s_{2}s_{4}s_{5}s_{6}s_{7}s_{4}s_{1} $,
$ s_{4}s_{3}s_{2}s_{4}s_{5}s_{6}s_{7}s_{4}s_{2} $,
$ s_{4}s_{3}s_{2}s_{4}s_{5}s_{6}s_{7}s_{4}s_{3} $,
$ s_{4}s_{3}s_{2}s_{4}s_{5}s_{6}s_{7}s_{4}s_{5} $,
$ s_{4}s_{3}s_{2}s_{4}s_{5}s_{6}s_{7}s_{5}s_{4} $,
$ s_{4}s_{5}s_{1}s_{3}s_{2}s_{4}s_{3}s_{2}s_{1} $,
$ s_{4}s_{5}s_{3}s_{2}s_{4}s_{1}s_{3}s_{2}s_{1} $,
$ s_{4}s_{5}s_{4}s_{1}s_{3}s_{2}s_{4}s_{2}s_{1} $,
$ s_{4}s_{5}s_{4}s_{1}s_{3}s_{2}s_{4}s_{3}s_{1} $,
$ s_{4}s_{5}s_{4}s_{3}s_{2}s_{4}s_{1}s_{3}s_{1} $,
$ s_{4}s_{5}s_{4}s_{3}s_{2}s_{4}s_{1}s_{3}s_{2} $,
$ s_{4}s_{5}s_{4}s_{3}s_{2}s_{4}s_{3}s_{2}s_{1} $,
$ s_{4}s_{5}s_{6}s_{1}s_{3}s_{2}s_{4}s_{3}s_{2} $,
$ s_{4}s_{5}s_{6}s_{1}s_{3}s_{2}s_{4}s_{5}s_{3} $,
$ s_{4}s_{5}s_{6}s_{1}s_{3}s_{4}s_{5}s_{2}s_{1} $,
$ s_{4}s_{5}s_{6}s_{1}s_{3}s_{4}s_{5}s_{3}s_{1} $,
$ s_{4}s_{5}s_{6}s_{1}s_{3}s_{4}s_{5}s_{3}s_{2} $,
$ s_{4}s_{5}s_{6}s_{1}s_{3}s_{4}s_{5}s_{4}s_{2} $,
$ s_{4}s_{5}s_{6}s_{1}s_{3}s_{4}s_{5}s_{4}s_{3} $,
$ s_{4}s_{5}s_{6}s_{2}s_{4}s_{5}s_{1}s_{3}s_{2} $,
$ s_{4}s_{5}s_{6}s_{2}s_{4}s_{5}s_{3}s_{2}s_{1} $,
$ s_{4}s_{5}s_{6}s_{2}s_{4}s_{5}s_{4}s_{1}s_{3} $,
$ s_{4}s_{5}s_{6}s_{2}s_{4}s_{5}s_{4}s_{3}s_{1} $,
$ s_{4}s_{5}s_{6}s_{2}s_{4}s_{5}s_{4}s_{3}s_{2} $,
$ s_{4}s_{5}s_{6}s_{3}s_{2}s_{4}s_{3}s_{2}s_{1} $,
$ s_{4}s_{5}s_{6}s_{3}s_{2}s_{4}s_{5}s_{3}s_{1} $,
$ s_{4}s_{5}s_{6}s_{3}s_{2}s_{4}s_{5}s_{4}s_{3} $,
$ s_{4}s_{5}s_{6}s_{3}s_{4}s_{5}s_{1}s_{3}s_{1} $,
$ s_{4}s_{5}s_{6}s_{3}s_{4}s_{5}s_{1}s_{3}s_{2} $,
$ s_{4}s_{5}s_{6}s_{3}s_{4}s_{5}s_{4}s_{2}s_{1} $,
$ s_{4}s_{5}s_{6}s_{3}s_{4}s_{5}s_{4}s_{3}s_{1} $,
$ s_{4}s_{5}s_{6}s_{4}s_{1}s_{3}s_{2}s_{4}s_{2} $,
$ s_{4}s_{5}s_{6}s_{4}s_{1}s_{3}s_{2}s_{4}s_{3} $,
$ s_{4}s_{5}s_{6}s_{4}s_{1}s_{3}s_{2}s_{4}s_{5} $,
$ s_{4}s_{5}s_{6}s_{4}s_{3}s_{2}s_{4}s_{2}s_{1} $,
$ s_{4}s_{5}s_{6}s_{4}s_{3}s_{2}s_{4}s_{3}s_{1} $,
$ s_{4}s_{5}s_{6}s_{4}s_{3}s_{2}s_{4}s_{5}s_{1} $,
$ s_{4}s_{5}s_{6}s_{4}s_{5}s_{1}s_{3}s_{4}s_{2} $,
$ s_{4}s_{5}s_{6}s_{4}s_{5}s_{1}s_{3}s_{4}s_{3} $,
$ s_{4}s_{5}s_{6}s_{4}s_{5}s_{2}s_{4}s_{1}s_{3} $,
$ s_{4}s_{5}s_{6}s_{4}s_{5}s_{2}s_{4}s_{3}s_{1} $,
$ s_{4}s_{5}s_{6}s_{4}s_{5}s_{3}s_{2}s_{4}s_{3} $,
$ s_{4}s_{5}s_{6}s_{4}s_{5}s_{3}s_{4}s_{2}s_{1} $,
$ s_{4}s_{5}s_{6}s_{4}s_{5}s_{3}s_{4}s_{3}s_{1} $,
$ s_{4}s_{5}s_{6}s_{4}s_{5}s_{4}s_{3}s_{2}s_{4} $,
$ s_{4}s_{5}s_{6}s_{5}s_{1}s_{3}s_{4}s_{3}s_{2} $,
$ s_{4}s_{5}s_{6}s_{5}s_{2}s_{4}s_{1}s_{3}s_{2} $,
$ s_{4}s_{5}s_{6}s_{5}s_{2}s_{4}s_{3}s_{2}s_{1} $,
$ s_{4}s_{5}s_{6}s_{5}s_{3}s_{2}s_{4}s_{3}s_{2} $,
$ s_{4}s_{5}s_{6}s_{5}s_{3}s_{4}s_{3}s_{2}s_{1} $,
$ s_{4}s_{5}s_{6}s_{5}s_{4}s_{3}s_{2}s_{4}s_{2} $,
$ s_{4}s_{5}s_{6}s_{5}s_{4}s_{3}s_{2}s_{4}s_{5} $,
$ s_{4}s_{5}s_{6}s_{7}s_{1}s_{3}s_{2}s_{4}s_{3} $,
$ s_{4}s_{5}s_{6}s_{7}s_{1}s_{3}s_{4}s_{2}s_{1} $,
$ s_{4}s_{5}s_{6}s_{7}s_{1}s_{3}s_{4}s_{3}s_{1} $,
$ s_{4}s_{5}s_{6}s_{7}s_{1}s_{3}s_{4}s_{3}s_{2} $,
$ s_{4}s_{5}s_{6}s_{7}s_{1}s_{3}s_{4}s_{5}s_{2} $,
$ s_{4}s_{5}s_{6}s_{7}s_{1}s_{3}s_{4}s_{5}s_{3} $,
$ s_{4}s_{5}s_{6}s_{7}s_{2}s_{4}s_{1}s_{3}s_{2} $,
$ s_{4}s_{5}s_{6}s_{7}s_{2}s_{4}s_{3}s_{2}s_{1} $,
$ s_{4}s_{5}s_{6}s_{7}s_{2}s_{4}s_{5}s_{1}s_{3} $,
$ s_{4}s_{5}s_{6}s_{7}s_{2}s_{4}s_{5}s_{3}s_{1} $,
$ s_{4}s_{5}s_{6}s_{7}s_{2}s_{4}s_{5}s_{4}s_{2} $,
$ s_{4}s_{5}s_{6}s_{7}s_{2}s_{4}s_{5}s_{4}s_{3} $,
$ s_{4}s_{5}s_{6}s_{7}s_{2}s_{4}s_{5}s_{6}s_{2} $,
$ s_{4}s_{5}s_{6}s_{7}s_{2}s_{4}s_{5}s_{6}s_{3} $,
$ s_{4}s_{5}s_{6}s_{7}s_{3}s_{2}s_{4}s_{3}s_{1} $,
$ s_{4}s_{5}s_{6}s_{7}s_{3}s_{2}s_{4}s_{3}s_{2} $,
$ s_{4}s_{5}s_{6}s_{7}s_{3}s_{2}s_{4}s_{5}s_{3} $,
$ s_{4}s_{5}s_{6}s_{7}s_{3}s_{4}s_{1}s_{3}s_{1} $,
$ s_{4}s_{5}s_{6}s_{7}s_{3}s_{4}s_{1}s_{3}s_{2} $,
$ s_{4}s_{5}s_{6}s_{7}s_{3}s_{4}s_{5}s_{2}s_{1} $,
$ s_{4}s_{5}s_{6}s_{7}s_{3}s_{4}s_{5}s_{3}s_{1} $,
$ s_{4}s_{5}s_{6}s_{7}s_{3}s_{4}s_{5}s_{4}s_{2} $,
$ s_{4}s_{5}s_{6}s_{7}s_{3}s_{4}s_{5}s_{4}s_{3} $,
$ s_{4}s_{5}s_{6}s_{7}s_{3}s_{4}s_{5}s_{6}s_{2} $,
$ s_{4}s_{5}s_{6}s_{7}s_{3}s_{4}s_{5}s_{6}s_{3} $,
$ s_{4}s_{5}s_{6}s_{7}s_{4}s_{1}s_{3}s_{2}s_{4} $,
$ s_{4}s_{5}s_{6}s_{7}s_{4}s_{3}s_{2}s_{4}s_{1} $,
$ s_{4}s_{5}s_{6}s_{7}s_{4}s_{3}s_{2}s_{4}s_{5} $,
$ s_{4}s_{5}s_{6}s_{7}s_{5}s_{1}s_{3}s_{4}s_{2} $,
$ s_{4}s_{5}s_{6}s_{7}s_{5}s_{1}s_{3}s_{4}s_{3} $,
$ s_{4}s_{5}s_{6}s_{7}s_{5}s_{2}s_{4}s_{1}s_{3} $,
$ s_{4}s_{5}s_{6}s_{7}s_{5}s_{2}s_{4}s_{3}s_{1} $,
$ s_{4}s_{5}s_{6}s_{7}s_{5}s_{3}s_{2}s_{4}s_{3} $,
$ s_{4}s_{5}s_{6}s_{7}s_{5}s_{3}s_{4}s_{2}s_{1} $,
$ s_{4}s_{5}s_{6}s_{7}s_{5}s_{3}s_{4}s_{3}s_{1} $,
$ s_{4}s_{5}s_{6}s_{7}s_{5}s_{3}s_{4}s_{3}s_{2} $,
$ s_{4}s_{5}s_{6}s_{7}s_{5}s_{4}s_{3}s_{2}s_{4} $,
$ s_{4}s_{5}s_{6}s_{7}s_{5}s_{6}s_{2}s_{4}s_{2} $,
$ s_{4}s_{5}s_{6}s_{7}s_{5}s_{6}s_{2}s_{4}s_{3} $,
$ s_{4}s_{5}s_{6}s_{7}s_{5}s_{6}s_{3}s_{4}s_{2} $,
$ s_{4}s_{5}s_{6}s_{7}s_{5}s_{6}s_{3}s_{4}s_{3} $,
$ s_{5}s_{3}s_{4}s_{1}s_{3}s_{2}s_{4}s_{3}s_{2} $,
$ s_{5}s_{3}s_{4}s_{1}s_{3}s_{2}s_{4}s_{5}s_{2} $,
$ s_{5}s_{3}s_{4}s_{1}s_{3}s_{2}s_{4}s_{5}s_{3} $,
$ s_{5}s_{3}s_{4}s_{1}s_{3}s_{2}s_{4}s_{5}s_{6} $,
$ s_{5}s_{4}s_{1}s_{3}s_{2}s_{4}s_{3}s_{2}s_{1} $,
$ s_{5}s_{4}s_{1}s_{3}s_{2}s_{4}s_{5}s_{2}s_{1} $,
$ s_{5}s_{4}s_{1}s_{3}s_{2}s_{4}s_{5}s_{2}s_{4} $,
$ s_{5}s_{4}s_{1}s_{3}s_{2}s_{4}s_{5}s_{3}s_{1} $,
$ s_{5}s_{4}s_{1}s_{3}s_{2}s_{4}s_{5}s_{3}s_{2} $,
$ s_{5}s_{4}s_{1}s_{3}s_{2}s_{4}s_{5}s_{4}s_{3} $,
$ s_{5}s_{4}s_{1}s_{3}s_{2}s_{4}s_{5}s_{6}s_{1} $,
$ s_{5}s_{4}s_{1}s_{3}s_{2}s_{4}s_{5}s_{6}s_{2} $,
$ s_{5}s_{4}s_{1}s_{3}s_{2}s_{4}s_{5}s_{6}s_{3} $,
$ s_{5}s_{4}s_{1}s_{3}s_{2}s_{4}s_{5}s_{6}s_{4} $,
$ s_{5}s_{4}s_{1}s_{3}s_{2}s_{4}s_{5}s_{6}s_{5} $,
$ s_{5}s_{4}s_{1}s_{3}s_{2}s_{4}s_{5}s_{6}s_{7} $,
$ s_{5}s_{4}s_{3}s_{2}s_{4}s_{1}s_{3}s_{2}s_{1} $,
$ s_{5}s_{4}s_{3}s_{2}s_{4}s_{5}s_{1}s_{3}s_{1} $,
$ s_{5}s_{4}s_{3}s_{2}s_{4}s_{5}s_{1}s_{3}s_{2} $,
$ s_{5}s_{4}s_{3}s_{2}s_{4}s_{5}s_{2}s_{4}s_{1} $,
$ s_{5}s_{4}s_{3}s_{2}s_{4}s_{5}s_{3}s_{2}s_{1} $,
$ s_{5}s_{4}s_{3}s_{2}s_{4}s_{5}s_{4}s_{3}s_{1} $,
$ s_{5}s_{4}s_{3}s_{2}s_{4}s_{5}s_{6}s_{1}s_{3} $,
$ s_{5}s_{4}s_{3}s_{2}s_{4}s_{5}s_{6}s_{3}s_{2} $,
$ s_{5}s_{4}s_{3}s_{2}s_{4}s_{5}s_{6}s_{4}s_{1} $,
$ s_{5}s_{4}s_{3}s_{2}s_{4}s_{5}s_{6}s_{4}s_{3} $,
$ s_{5}s_{4}s_{3}s_{2}s_{4}s_{5}s_{6}s_{5}s_{1} $,
$ s_{5}s_{4}s_{3}s_{2}s_{4}s_{5}s_{6}s_{5}s_{3} $,
$ s_{5}s_{4}s_{3}s_{2}s_{4}s_{5}s_{6}s_{5}s_{4} $,
$ s_{5}s_{4}s_{3}s_{2}s_{4}s_{5}s_{6}s_{7}s_{1} $,
$ s_{5}s_{4}s_{3}s_{2}s_{4}s_{5}s_{6}s_{7}s_{2} $,
$ s_{5}s_{4}s_{3}s_{2}s_{4}s_{5}s_{6}s_{7}s_{4} $,
$ s_{5}s_{4}s_{3}s_{2}s_{4}s_{5}s_{6}s_{7}s_{5} $,
$ s_{5}s_{6}s_{1}s_{3}s_{2}s_{4}s_{1}s_{3}s_{2} $,
$ s_{5}s_{6}s_{1}s_{3}s_{2}s_{4}s_{3}s_{2}s_{1} $,
$ s_{5}s_{6}s_{1}s_{3}s_{2}s_{4}s_{5}s_{1}s_{3} $,
$ s_{5}s_{6}s_{1}s_{3}s_{2}s_{4}s_{5}s_{2}s_{4} $,
$ s_{5}s_{6}s_{1}s_{3}s_{2}s_{4}s_{5}s_{3}s_{1} $,
$ s_{5}s_{6}s_{1}s_{3}s_{2}s_{4}s_{5}s_{3}s_{2} $,
$ s_{5}s_{6}s_{1}s_{3}s_{2}s_{4}s_{5}s_{4}s_{3} $,
$ s_{5}s_{6}s_{1}s_{3}s_{4}s_{5}s_{1}s_{3}s_{2} $,
$ s_{5}s_{6}s_{1}s_{3}s_{4}s_{5}s_{1}s_{3}s_{4} $,
$ s_{5}s_{6}s_{1}s_{3}s_{4}s_{5}s_{2}s_{4}s_{3} $,
$ s_{5}s_{6}s_{1}s_{3}s_{4}s_{5}s_{3}s_{2}s_{1} $,
$ s_{5}s_{6}s_{1}s_{3}s_{4}s_{5}s_{3}s_{2}s_{4} $,
$ s_{5}s_{6}s_{1}s_{3}s_{4}s_{5}s_{3}s_{4}s_{1} $,
$ s_{5}s_{6}s_{1}s_{3}s_{4}s_{5}s_{3}s_{4}s_{3} $,
$ s_{5}s_{6}s_{1}s_{3}s_{4}s_{5}s_{4}s_{1}s_{3} $,
$ s_{5}s_{6}s_{1}s_{3}s_{4}s_{5}s_{4}s_{3}s_{1} $,
$ s_{5}s_{6}s_{1}s_{3}s_{4}s_{5}s_{4}s_{3}s_{2} $,
$ s_{5}s_{6}s_{2}s_{4}s_{5}s_{1}s_{3}s_{4}s_{1} $,
$ s_{5}s_{6}s_{2}s_{4}s_{5}s_{1}s_{3}s_{4}s_{3} $,
$ s_{5}s_{6}s_{2}s_{4}s_{5}s_{3}s_{4}s_{1}s_{3} $,
$ s_{5}s_{6}s_{2}s_{4}s_{5}s_{4}s_{1}s_{3}s_{2} $,
$ s_{5}s_{6}s_{2}s_{4}s_{5}s_{4}s_{3}s_{2}s_{1} $,
$ s_{5}s_{6}s_{3}s_{2}s_{4}s_{1}s_{3}s_{2}s_{1} $,
$ s_{5}s_{6}s_{3}s_{2}s_{4}s_{5}s_{1}s_{3}s_{1} $,
$ s_{5}s_{6}s_{3}s_{2}s_{4}s_{5}s_{1}s_{3}s_{4} $,
$ s_{5}s_{6}s_{3}s_{2}s_{4}s_{5}s_{2}s_{4}s_{1} $,
$ s_{5}s_{6}s_{3}s_{2}s_{4}s_{5}s_{3}s_{2}s_{1} $,
$ s_{5}s_{6}s_{3}s_{2}s_{4}s_{5}s_{3}s_{4}s_{1} $,
$ s_{5}s_{6}s_{3}s_{2}s_{4}s_{5}s_{3}s_{4}s_{3} $,
$ s_{5}s_{6}s_{3}s_{2}s_{4}s_{5}s_{4}s_{1}s_{3} $,
$ s_{5}s_{6}s_{3}s_{2}s_{4}s_{5}s_{4}s_{3}s_{2} $,
$ s_{5}s_{6}s_{3}s_{4}s_{1}s_{3}s_{2}s_{4}s_{2} $,
$ s_{5}s_{6}s_{3}s_{4}s_{1}s_{3}s_{2}s_{4}s_{3} $,
$ s_{5}s_{6}s_{3}s_{4}s_{5}s_{1}s_{3}s_{2}s_{1} $,
$ s_{5}s_{6}s_{3}s_{4}s_{5}s_{1}s_{3}s_{4}s_{2} $,
$ s_{5}s_{6}s_{3}s_{4}s_{5}s_{1}s_{3}s_{4}s_{3} $,
$ s_{5}s_{6}s_{3}s_{4}s_{5}s_{2}s_{4}s_{3}s_{1} $,
$ s_{5}s_{6}s_{3}s_{4}s_{5}s_{3}s_{2}s_{4}s_{1} $,
$ s_{5}s_{6}s_{3}s_{4}s_{5}s_{3}s_{4}s_{1}s_{3} $,
$ s_{5}s_{6}s_{3}s_{4}s_{5}s_{3}s_{4}s_{2}s_{1} $,
$ s_{5}s_{6}s_{3}s_{4}s_{5}s_{3}s_{4}s_{3}s_{1} $,
$ s_{5}s_{6}s_{3}s_{4}s_{5}s_{3}s_{4}s_{3}s_{2} $,
$ s_{5}s_{6}s_{3}s_{4}s_{5}s_{4}s_{1}s_{3}s_{1} $,
$ s_{5}s_{6}s_{3}s_{4}s_{5}s_{4}s_{1}s_{3}s_{2} $,
$ s_{5}s_{6}s_{3}s_{4}s_{5}s_{4}s_{3}s_{2}s_{4} $,
$ s_{5}s_{6}s_{4}s_{1}s_{3}s_{2}s_{4}s_{2}s_{1} $,
$ s_{5}s_{6}s_{4}s_{1}s_{3}s_{2}s_{4}s_{3}s_{1} $,
$ s_{5}s_{6}s_{4}s_{1}s_{3}s_{2}s_{4}s_{5}s_{2} $,
$ s_{5}s_{6}s_{4}s_{1}s_{3}s_{2}s_{4}s_{5}s_{3} $,
$ s_{5}s_{6}s_{4}s_{1}s_{3}s_{2}s_{4}s_{5}s_{4} $,
$ s_{5}s_{6}s_{4}s_{3}s_{2}s_{4}s_{1}s_{3}s_{1} $,
$ s_{5}s_{6}s_{4}s_{3}s_{2}s_{4}s_{1}s_{3}s_{2} $,
$ s_{5}s_{6}s_{4}s_{3}s_{2}s_{4}s_{5}s_{2}s_{1} $,
$ s_{5}s_{6}s_{4}s_{3}s_{2}s_{4}s_{5}s_{3}s_{1} $,
$ s_{5}s_{6}s_{4}s_{3}s_{2}s_{4}s_{5}s_{4}s_{1} $,
$ s_{5}s_{6}s_{4}s_{5}s_{1}s_{3}s_{2}s_{4}s_{3} $,
$ s_{5}s_{6}s_{4}s_{5}s_{1}s_{3}s_{4}s_{2}s_{1} $,
$ s_{5}s_{6}s_{4}s_{5}s_{1}s_{3}s_{4}s_{3}s_{1} $,
$ s_{5}s_{6}s_{4}s_{5}s_{1}s_{3}s_{4}s_{3}s_{2} $,
$ s_{5}s_{6}s_{4}s_{5}s_{2}s_{4}s_{1}s_{3}s_{2} $,
$ s_{5}s_{6}s_{4}s_{5}s_{2}s_{4}s_{3}s_{2}s_{1} $,
$ s_{5}s_{6}s_{4}s_{5}s_{3}s_{2}s_{4}s_{3}s_{1} $,
$ s_{5}s_{6}s_{4}s_{5}s_{3}s_{4}s_{1}s_{3}s_{1} $,
$ s_{5}s_{6}s_{4}s_{5}s_{3}s_{4}s_{1}s_{3}s_{2} $,
$ s_{5}s_{6}s_{4}s_{5}s_{4}s_{1}s_{3}s_{2}s_{4} $,
$ s_{5}s_{6}s_{4}s_{5}s_{4}s_{3}s_{2}s_{4}s_{1} $,
$ s_{5}s_{6}s_{5}s_{1}s_{3}s_{2}s_{4}s_{3}s_{2} $,
$ s_{5}s_{6}s_{5}s_{1}s_{3}s_{4}s_{1}s_{3}s_{2} $,
$ s_{5}s_{6}s_{5}s_{1}s_{3}s_{4}s_{3}s_{2}s_{1} $,
$ s_{5}s_{6}s_{5}s_{3}s_{2}s_{4}s_{1}s_{3}s_{2} $,
$ s_{5}s_{6}s_{5}s_{3}s_{4}s_{1}s_{3}s_{2}s_{1} $,
$ s_{5}s_{6}s_{5}s_{3}s_{4}s_{1}s_{3}s_{2}s_{4} $,
$ s_{5}s_{6}s_{5}s_{4}s_{1}s_{3}s_{2}s_{4}s_{1} $,
$ s_{5}s_{6}s_{5}s_{4}s_{1}s_{3}s_{2}s_{4}s_{2} $,
$ s_{5}s_{6}s_{5}s_{4}s_{1}s_{3}s_{2}s_{4}s_{3} $,
$ s_{5}s_{6}s_{5}s_{4}s_{1}s_{3}s_{2}s_{4}s_{5} $,
$ s_{5}s_{6}s_{5}s_{4}s_{3}s_{2}s_{4}s_{1}s_{3} $,
$ s_{5}s_{6}s_{5}s_{4}s_{3}s_{2}s_{4}s_{5}s_{1} $,
$ s_{5}s_{6}s_{5}s_{4}s_{3}s_{2}s_{4}s_{5}s_{2} $,
$ s_{5}s_{6}s_{7}s_{1}s_{3}s_{2}s_{4}s_{3}s_{2} $,
$ s_{5}s_{6}s_{7}s_{1}s_{3}s_{2}s_{4}s_{5}s_{2} $,
$ s_{5}s_{6}s_{7}s_{1}s_{3}s_{2}s_{4}s_{5}s_{3} $,
$ s_{5}s_{6}s_{7}s_{1}s_{3}s_{4}s_{1}s_{3}s_{2} $,
$ s_{5}s_{6}s_{7}s_{1}s_{3}s_{4}s_{3}s_{2}s_{1} $,
$ s_{5}s_{6}s_{7}s_{1}s_{3}s_{4}s_{5}s_{2}s_{4} $,
$ s_{5}s_{6}s_{7}s_{1}s_{3}s_{4}s_{5}s_{3}s_{2} $,
$ s_{5}s_{6}s_{7}s_{1}s_{3}s_{4}s_{5}s_{3}s_{4} $,
$ s_{5}s_{6}s_{7}s_{1}s_{3}s_{4}s_{5}s_{4}s_{2} $,
$ s_{5}s_{6}s_{7}s_{1}s_{3}s_{4}s_{5}s_{4}s_{3} $,
$ s_{5}s_{6}s_{7}s_{2}s_{4}s_{5}s_{1}s_{3}s_{2} $,
$ s_{5}s_{6}s_{7}s_{2}s_{4}s_{5}s_{1}s_{3}s_{4} $,
$ s_{5}s_{6}s_{7}s_{2}s_{4}s_{5}s_{3}s_{2}s_{1} $,
$ s_{5}s_{6}s_{7}s_{2}s_{4}s_{5}s_{3}s_{4}s_{1} $,
$ s_{5}s_{6}s_{7}s_{2}s_{4}s_{5}s_{3}s_{4}s_{3} $,
$ s_{5}s_{6}s_{7}s_{2}s_{4}s_{5}s_{4}s_{1}s_{3} $,
$ s_{5}s_{6}s_{7}s_{2}s_{4}s_{5}s_{4}s_{3}s_{1} $,
$ s_{5}s_{6}s_{7}s_{2}s_{4}s_{5}s_{6}s_{2}s_{4} $,
$ s_{5}s_{6}s_{7}s_{2}s_{4}s_{5}s_{6}s_{3}s_{2} $,
$ s_{5}s_{6}s_{7}s_{2}s_{4}s_{5}s_{6}s_{3}s_{4} $,
$ s_{5}s_{6}s_{7}s_{2}s_{4}s_{5}s_{6}s_{4}s_{2} $,
$ s_{5}s_{6}s_{7}s_{2}s_{4}s_{5}s_{6}s_{4}s_{3} $,
$ s_{5}s_{6}s_{7}s_{3}s_{2}s_{4}s_{1}s_{3}s_{2} $,
$ s_{5}s_{6}s_{7}s_{3}s_{2}s_{4}s_{5}s_{2}s_{1} $,
$ s_{5}s_{6}s_{7}s_{3}s_{2}s_{4}s_{5}s_{3}s_{1} $,
$ s_{5}s_{6}s_{7}s_{3}s_{2}s_{4}s_{5}s_{3}s_{4} $,
$ s_{5}s_{6}s_{7}s_{3}s_{2}s_{4}s_{5}s_{4}s_{2} $,
$ s_{5}s_{6}s_{7}s_{3}s_{2}s_{4}s_{5}s_{6}s_{2} $,
$ s_{5}s_{6}s_{7}s_{3}s_{2}s_{4}s_{5}s_{6}s_{3} $,
$ s_{5}s_{6}s_{7}s_{3}s_{4}s_{1}s_{3}s_{2}s_{1} $,
$ s_{5}s_{6}s_{7}s_{3}s_{4}s_{1}s_{3}s_{2}s_{4} $,
$ s_{5}s_{6}s_{7}s_{3}s_{4}s_{5}s_{2}s_{4}s_{1} $,
$ s_{5}s_{6}s_{7}s_{3}s_{4}s_{5}s_{3}s_{2}s_{1} $,
$ s_{5}s_{6}s_{7}s_{3}s_{4}s_{5}s_{3}s_{4}s_{1} $,
$ s_{5}s_{6}s_{7}s_{3}s_{4}s_{5}s_{3}s_{4}s_{2} $,
$ s_{5}s_{6}s_{7}s_{3}s_{4}s_{5}s_{4}s_{2}s_{1} $,
$ s_{5}s_{6}s_{7}s_{3}s_{4}s_{5}s_{4}s_{3}s_{1} $,
$ s_{5}s_{6}s_{7}s_{3}s_{4}s_{5}s_{4}s_{3}s_{2} $,
$ s_{5}s_{6}s_{7}s_{3}s_{4}s_{5}s_{6}s_{2}s_{4} $,
$ s_{5}s_{6}s_{7}s_{3}s_{4}s_{5}s_{6}s_{3}s_{2} $,
$ s_{5}s_{6}s_{7}s_{3}s_{4}s_{5}s_{6}s_{3}s_{4} $,
$ s_{5}s_{6}s_{7}s_{3}s_{4}s_{5}s_{6}s_{4}s_{2} $,
$ s_{5}s_{6}s_{7}s_{3}s_{4}s_{5}s_{6}s_{4}s_{3} $,
$ s_{5}s_{6}s_{7}s_{4}s_{1}s_{3}s_{2}s_{4}s_{1} $,
$ s_{5}s_{6}s_{7}s_{4}s_{1}s_{3}s_{2}s_{4}s_{2} $,
$ s_{5}s_{6}s_{7}s_{4}s_{1}s_{3}s_{2}s_{4}s_{3} $,
$ s_{5}s_{6}s_{7}s_{4}s_{3}s_{2}s_{4}s_{1}s_{3} $,
$ s_{5}s_{6}s_{7}s_{4}s_{3}s_{2}s_{4}s_{5}s_{2} $,
$ s_{5}s_{6}s_{7}s_{4}s_{3}s_{2}s_{4}s_{5}s_{3} $,
$ s_{5}s_{6}s_{7}s_{4}s_{3}s_{2}s_{4}s_{5}s_{4} $,
$ s_{5}s_{6}s_{7}s_{4}s_{5}s_{1}s_{3}s_{4}s_{2} $,
$ s_{5}s_{6}s_{7}s_{4}s_{5}s_{1}s_{3}s_{4}s_{3} $,
$ s_{5}s_{6}s_{7}s_{4}s_{5}s_{2}s_{4}s_{1}s_{3} $,
$ s_{5}s_{6}s_{7}s_{4}s_{5}s_{2}s_{4}s_{3}s_{1} $,
$ s_{5}s_{6}s_{7}s_{4}s_{5}s_{3}s_{4}s_{2}s_{1} $,
$ s_{5}s_{6}s_{7}s_{4}s_{5}s_{3}s_{4}s_{3}s_{1} $,
$ s_{5}s_{6}s_{7}s_{4}s_{5}s_{3}s_{4}s_{3}s_{2} $,
$ s_{5}s_{6}s_{7}s_{4}s_{5}s_{4}s_{3}s_{2}s_{4} $,
$ s_{5}s_{6}s_{7}s_{4}s_{5}s_{6}s_{2}s_{4}s_{2} $,
$ s_{5}s_{6}s_{7}s_{4}s_{5}s_{6}s_{2}s_{4}s_{3} $,
$ s_{5}s_{6}s_{7}s_{4}s_{5}s_{6}s_{3}s_{4}s_{2} $,
$ s_{5}s_{6}s_{7}s_{4}s_{5}s_{6}s_{3}s_{4}s_{3} $,
$ s_{5}s_{6}s_{7}s_{5}s_{4}s_{3}s_{2}s_{4}s_{5} $,
$ s_{6}s_{1}s_{3}s_{2}s_{4}s_{5}s_{1}s_{3}s_{2} $,
$ s_{6}s_{1}s_{3}s_{2}s_{4}s_{5}s_{2}s_{4}s_{2} $,
$ s_{6}s_{1}s_{3}s_{2}s_{4}s_{5}s_{3}s_{2}s_{1} $,
$ s_{6}s_{1}s_{3}s_{2}s_{4}s_{5}s_{3}s_{2}s_{4} $,
$ s_{6}s_{1}s_{3}s_{2}s_{4}s_{5}s_{4}s_{1}s_{3} $,
$ s_{6}s_{1}s_{3}s_{2}s_{4}s_{5}s_{4}s_{3}s_{1} $,
$ s_{6}s_{1}s_{3}s_{2}s_{4}s_{5}s_{4}s_{3}s_{2} $,
$ s_{6}s_{1}s_{3}s_{4}s_{5}s_{1}s_{3}s_{4}s_{3} $,
$ s_{6}s_{1}s_{3}s_{4}s_{5}s_{2}s_{4}s_{3}s_{2} $,
$ s_{6}s_{1}s_{3}s_{4}s_{5}s_{3}s_{2}s_{4}s_{2} $,
$ s_{6}s_{1}s_{3}s_{4}s_{5}s_{3}s_{4}s_{3}s_{1} $,
$ s_{6}s_{2}s_{4}s_{5}s_{1}s_{3}s_{2}s_{4}s_{3} $,
$ s_{6}s_{2}s_{4}s_{5}s_{1}s_{3}s_{4}s_{3}s_{1} $,
$ s_{6}s_{2}s_{4}s_{5}s_{3}s_{2}s_{4}s_{3}s_{1} $,
$ s_{6}s_{2}s_{4}s_{5}s_{3}s_{4}s_{1}s_{3}s_{1} $,
$ s_{6}s_{2}s_{4}s_{5}s_{4}s_{1}s_{3}s_{2}s_{4} $,
$ s_{6}s_{2}s_{4}s_{5}s_{4}s_{3}s_{2}s_{4}s_{1} $,
$ s_{6}s_{3}s_{2}s_{4}s_{5}s_{1}s_{3}s_{2}s_{1} $,
$ s_{6}s_{3}s_{2}s_{4}s_{5}s_{1}s_{3}s_{4}s_{3} $,
$ s_{6}s_{3}s_{2}s_{4}s_{5}s_{2}s_{4}s_{2}s_{1} $,
$ s_{6}s_{3}s_{2}s_{4}s_{5}s_{3}s_{2}s_{4}s_{1} $,
$ s_{6}s_{3}s_{2}s_{4}s_{5}s_{3}s_{4}s_{3}s_{1} $,
$ s_{6}s_{3}s_{2}s_{4}s_{5}s_{4}s_{1}s_{3}s_{1} $,
$ s_{6}s_{3}s_{2}s_{4}s_{5}s_{4}s_{1}s_{3}s_{2} $,
$ s_{6}s_{3}s_{2}s_{4}s_{5}s_{4}s_{3}s_{2}s_{4} $,
$ s_{6}s_{3}s_{4}s_{1}s_{3}s_{2}s_{4}s_{5}s_{2} $,
$ s_{6}s_{3}s_{4}s_{1}s_{3}s_{2}s_{4}s_{5}s_{3} $,
$ s_{6}s_{3}s_{4}s_{5}s_{1}s_{3}s_{4}s_{3}s_{2} $,
$ s_{6}s_{3}s_{4}s_{5}s_{2}s_{4}s_{3}s_{2}s_{1} $,
$ s_{6}s_{3}s_{4}s_{5}s_{3}s_{2}s_{4}s_{2}s_{1} $,
$ s_{6}s_{3}s_{4}s_{5}s_{3}s_{4}s_{1}s_{3}s_{1} $,
$ s_{6}s_{3}s_{4}s_{5}s_{3}s_{4}s_{3}s_{2}s_{1} $,
$ s_{6}s_{3}s_{4}s_{5}s_{4}s_{3}s_{2}s_{4}s_{2} $,
$ s_{6}s_{4}s_{1}s_{3}s_{2}s_{4}s_{5}s_{2}s_{1} $,
$ s_{6}s_{4}s_{1}s_{3}s_{2}s_{4}s_{5}s_{2}s_{4} $,
$ s_{6}s_{4}s_{1}s_{3}s_{2}s_{4}s_{5}s_{3}s_{1} $,
$ s_{6}s_{4}s_{1}s_{3}s_{2}s_{4}s_{5}s_{3}s_{4} $,
$ s_{6}s_{4}s_{1}s_{3}s_{2}s_{4}s_{5}s_{4}s_{2} $,
$ s_{6}s_{4}s_{1}s_{3}s_{2}s_{4}s_{5}s_{4}s_{3} $,
$ s_{6}s_{4}s_{3}s_{2}s_{4}s_{5}s_{1}s_{3}s_{1} $,
$ s_{6}s_{4}s_{3}s_{2}s_{4}s_{5}s_{1}s_{3}s_{2} $,
$ s_{6}s_{4}s_{3}s_{2}s_{4}s_{5}s_{2}s_{4}s_{1} $,
$ s_{6}s_{4}s_{3}s_{2}s_{4}s_{5}s_{3}s_{4}s_{1} $,
$ s_{6}s_{4}s_{3}s_{2}s_{4}s_{5}s_{4}s_{2}s_{1} $,
$ s_{6}s_{4}s_{3}s_{2}s_{4}s_{5}s_{4}s_{3}s_{1} $,
$ s_{6}s_{4}s_{3}s_{2}s_{4}s_{5}s_{4}s_{3}s_{2} $,
$ s_{6}s_{4}s_{5}s_{1}s_{3}s_{2}s_{4}s_{3}s_{2} $,
$ s_{6}s_{4}s_{5}s_{1}s_{3}s_{4}s_{3}s_{2}s_{1} $,
$ s_{6}s_{4}s_{5}s_{3}s_{2}s_{4}s_{3}s_{2}s_{1} $,
$ s_{6}s_{4}s_{5}s_{3}s_{4}s_{1}s_{3}s_{2}s_{1} $,
$ s_{6}s_{4}s_{5}s_{4}s_{1}s_{3}s_{2}s_{4}s_{2} $,
$ s_{6}s_{4}s_{5}s_{4}s_{3}s_{2}s_{4}s_{2}s_{1} $,
$ s_{6}s_{5}s_{1}s_{3}s_{2}s_{4}s_{1}s_{3}s_{2} $,
$ s_{6}s_{5}s_{1}s_{3}s_{2}s_{4}s_{3}s_{2}s_{1} $,
$ s_{6}s_{5}s_{3}s_{2}s_{4}s_{1}s_{3}s_{2}s_{1} $,
$ s_{6}s_{5}s_{3}s_{4}s_{1}s_{3}s_{2}s_{4}s_{2} $,
$ s_{6}s_{5}s_{3}s_{4}s_{1}s_{3}s_{2}s_{4}s_{3} $,
$ s_{6}s_{5}s_{3}s_{4}s_{1}s_{3}s_{2}s_{4}s_{5} $,
$ s_{6}s_{5}s_{4}s_{1}s_{3}s_{2}s_{4}s_{2}s_{1} $,
$ s_{6}s_{5}s_{4}s_{1}s_{3}s_{2}s_{4}s_{3}s_{1} $,
$ s_{6}s_{5}s_{4}s_{1}s_{3}s_{2}s_{4}s_{5}s_{1} $,
$ s_{6}s_{5}s_{4}s_{1}s_{3}s_{2}s_{4}s_{5}s_{4} $,
$ s_{6}s_{5}s_{4}s_{3}s_{2}s_{4}s_{1}s_{3}s_{1} $,
$ s_{6}s_{5}s_{4}s_{3}s_{2}s_{4}s_{1}s_{3}s_{2} $,
$ s_{6}s_{5}s_{4}s_{3}s_{2}s_{4}s_{5}s_{1}s_{3} $,
$ s_{6}s_{5}s_{4}s_{3}s_{2}s_{4}s_{5}s_{2}s_{1} $,
$ s_{6}s_{5}s_{4}s_{3}s_{2}s_{4}s_{5}s_{2}s_{4} $,
$ s_{6}s_{5}s_{4}s_{3}s_{2}s_{4}s_{5}s_{3}s_{1} $,
$ s_{6}s_{5}s_{4}s_{3}s_{2}s_{4}s_{5}s_{3}s_{2} $,
$ s_{6}s_{5}s_{4}s_{3}s_{2}s_{4}s_{5}s_{4}s_{1} $,
$ s_{6}s_{5}s_{4}s_{3}s_{2}s_{4}s_{5}s_{6}s_{2} $,
$ s_{6}s_{5}s_{4}s_{3}s_{2}s_{4}s_{5}s_{6}s_{3} $,
$ s_{6}s_{5}s_{4}s_{3}s_{2}s_{4}s_{5}s_{6}s_{7} $,
$ s_{6}s_{7}s_{1}s_{3}s_{2}s_{4}s_{5}s_{1}s_{3} $,
$ s_{6}s_{7}s_{1}s_{3}s_{2}s_{4}s_{5}s_{2}s_{4} $,
$ s_{6}s_{7}s_{1}s_{3}s_{2}s_{4}s_{5}s_{3}s_{1} $,
$ s_{6}s_{7}s_{1}s_{3}s_{2}s_{4}s_{5}s_{4}s_{2} $,
$ s_{6}s_{7}s_{1}s_{3}s_{4}s_{5}s_{1}s_{3}s_{4} $,
$ s_{6}s_{7}s_{1}s_{3}s_{4}s_{5}s_{2}s_{4}s_{2} $,
$ s_{6}s_{7}s_{1}s_{3}s_{4}s_{5}s_{2}s_{4}s_{3} $,
$ s_{6}s_{7}s_{1}s_{3}s_{4}s_{5}s_{3}s_{2}s_{4} $,
$ s_{6}s_{7}s_{1}s_{3}s_{4}s_{5}s_{3}s_{4}s_{1} $,
$ s_{6}s_{7}s_{1}s_{3}s_{4}s_{5}s_{4}s_{1}s_{3} $,
$ s_{6}s_{7}s_{1}s_{3}s_{4}s_{5}s_{4}s_{3}s_{1} $,
$ s_{6}s_{7}s_{2}s_{4}s_{5}s_{1}s_{3}s_{4}s_{1} $,
$ s_{6}s_{7}s_{2}s_{4}s_{5}s_{3}s_{4}s_{1}s_{3} $,
$ s_{6}s_{7}s_{2}s_{4}s_{5}s_{3}s_{4}s_{3}s_{1} $,
$ s_{6}s_{7}s_{2}s_{4}s_{5}s_{4}s_{3}s_{2}s_{4} $,
$ s_{6}s_{7}s_{2}s_{4}s_{5}s_{6}s_{2}s_{4}s_{3} $,
$ s_{6}s_{7}s_{2}s_{4}s_{5}s_{6}s_{2}s_{4}s_{5} $,
$ s_{6}s_{7}s_{2}s_{4}s_{5}s_{6}s_{3}s_{4}s_{3} $,
$ s_{6}s_{7}s_{2}s_{4}s_{5}s_{6}s_{3}s_{4}s_{5} $,
$ s_{6}s_{7}s_{2}s_{4}s_{5}s_{6}s_{4}s_{3}s_{2} $,
$ s_{6}s_{7}s_{2}s_{4}s_{5}s_{6}s_{4}s_{5}s_{2} $,
$ s_{6}s_{7}s_{2}s_{4}s_{5}s_{6}s_{4}s_{5}s_{3} $,
$ s_{6}s_{7}s_{2}s_{4}s_{5}s_{6}s_{5}s_{2}s_{4} $,
$ s_{6}s_{7}s_{2}s_{4}s_{5}s_{6}s_{5}s_{3}s_{4} $,
$ s_{6}s_{7}s_{2}s_{4}s_{5}s_{6}s_{5}s_{4}s_{2} $,
$ s_{6}s_{7}s_{2}s_{4}s_{5}s_{6}s_{5}s_{4}s_{3} $,
$ s_{6}s_{7}s_{3}s_{2}s_{4}s_{5}s_{1}s_{3}s_{1} $,
$ s_{6}s_{7}s_{3}s_{2}s_{4}s_{5}s_{1}s_{3}s_{2} $,
$ s_{6}s_{7}s_{3}s_{2}s_{4}s_{5}s_{1}s_{3}s_{4} $,
$ s_{6}s_{7}s_{3}s_{2}s_{4}s_{5}s_{2}s_{4}s_{1} $,
$ s_{6}s_{7}s_{3}s_{2}s_{4}s_{5}s_{3}s_{2}s_{1} $,
$ s_{6}s_{7}s_{3}s_{2}s_{4}s_{5}s_{3}s_{4}s_{1} $,
$ s_{6}s_{7}s_{3}s_{2}s_{4}s_{5}s_{4}s_{1}s_{3} $,
$ s_{6}s_{7}s_{3}s_{2}s_{4}s_{5}s_{4}s_{2}s_{1} $,
$ s_{6}s_{7}s_{3}s_{2}s_{4}s_{5}s_{4}s_{3}s_{1} $,
$ s_{6}s_{7}s_{3}s_{2}s_{4}s_{5}s_{6}s_{2}s_{4} $,
$ s_{6}s_{7}s_{3}s_{2}s_{4}s_{5}s_{6}s_{3}s_{4} $,
$ s_{6}s_{7}s_{3}s_{2}s_{4}s_{5}s_{6}s_{4}s_{2} $,
$ s_{6}s_{7}s_{3}s_{2}s_{4}s_{5}s_{6}s_{4}s_{3} $,
$ s_{6}s_{7}s_{3}s_{4}s_{1}s_{3}s_{2}s_{4}s_{5} $,
$ s_{6}s_{7}s_{3}s_{4}s_{5}s_{1}s_{3}s_{4}s_{2} $,
$ s_{6}s_{7}s_{3}s_{4}s_{5}s_{1}s_{3}s_{4}s_{3} $,
$ s_{6}s_{7}s_{3}s_{4}s_{5}s_{2}s_{4}s_{2}s_{1} $,
$ s_{6}s_{7}s_{3}s_{4}s_{5}s_{2}s_{4}s_{3}s_{1} $,
$ s_{6}s_{7}s_{3}s_{4}s_{5}s_{2}s_{4}s_{3}s_{2} $,
$ s_{6}s_{7}s_{3}s_{4}s_{5}s_{3}s_{2}s_{4}s_{1} $,
$ s_{6}s_{7}s_{3}s_{4}s_{5}s_{3}s_{2}s_{4}s_{2} $,
$ s_{6}s_{7}s_{3}s_{4}s_{5}s_{3}s_{4}s_{1}s_{3} $,
$ s_{6}s_{7}s_{3}s_{4}s_{5}s_{3}s_{4}s_{2}s_{1} $,
$ s_{6}s_{7}s_{3}s_{4}s_{5}s_{4}s_{1}s_{3}s_{1} $,
$ s_{6}s_{7}s_{3}s_{4}s_{5}s_{4}s_{1}s_{3}s_{2} $,
$ s_{6}s_{7}s_{3}s_{4}s_{5}s_{4}s_{3}s_{2}s_{1} $,
$ s_{6}s_{7}s_{3}s_{4}s_{5}s_{4}s_{3}s_{2}s_{4} $,
$ s_{6}s_{7}s_{3}s_{4}s_{5}s_{6}s_{2}s_{4}s_{2} $,
$ s_{6}s_{7}s_{3}s_{4}s_{5}s_{6}s_{2}s_{4}s_{5} $,
$ s_{6}s_{7}s_{3}s_{4}s_{5}s_{6}s_{3}s_{4}s_{2} $,
$ s_{6}s_{7}s_{3}s_{4}s_{5}s_{6}s_{3}s_{4}s_{5} $,
$ s_{6}s_{7}s_{3}s_{4}s_{5}s_{6}s_{4}s_{3}s_{2} $,
$ s_{6}s_{7}s_{3}s_{4}s_{5}s_{6}s_{4}s_{5}s_{2} $,
$ s_{6}s_{7}s_{3}s_{4}s_{5}s_{6}s_{4}s_{5}s_{3} $,
$ s_{6}s_{7}s_{3}s_{4}s_{5}s_{6}s_{5}s_{2}s_{4} $,
$ s_{6}s_{7}s_{3}s_{4}s_{5}s_{6}s_{5}s_{3}s_{4} $,
$ s_{6}s_{7}s_{3}s_{4}s_{5}s_{6}s_{5}s_{4}s_{2} $,
$ s_{6}s_{7}s_{3}s_{4}s_{5}s_{6}s_{5}s_{4}s_{3} $,
$ s_{6}s_{7}s_{4}s_{1}s_{3}s_{2}s_{4}s_{5}s_{1} $,
$ s_{6}s_{7}s_{4}s_{1}s_{3}s_{2}s_{4}s_{5}s_{4} $,
$ s_{6}s_{7}s_{4}s_{3}s_{2}s_{4}s_{5}s_{1}s_{3} $,
$ s_{6}s_{7}s_{4}s_{3}s_{2}s_{4}s_{5}s_{2}s_{1} $,
$ s_{6}s_{7}s_{4}s_{3}s_{2}s_{4}s_{5}s_{2}s_{4} $,
$ s_{6}s_{7}s_{4}s_{3}s_{2}s_{4}s_{5}s_{3}s_{1} $,
$ s_{6}s_{7}s_{4}s_{3}s_{2}s_{4}s_{5}s_{3}s_{4} $,
$ s_{6}s_{7}s_{4}s_{3}s_{2}s_{4}s_{5}s_{4}s_{1} $,
$ s_{6}s_{7}s_{4}s_{3}s_{2}s_{4}s_{5}s_{4}s_{2} $,
$ s_{6}s_{7}s_{4}s_{5}s_{1}s_{3}s_{2}s_{4}s_{3} $,
$ s_{6}s_{7}s_{4}s_{5}s_{1}s_{3}s_{4}s_{2}s_{1} $,
$ s_{6}s_{7}s_{4}s_{5}s_{1}s_{3}s_{4}s_{3}s_{1} $,
$ s_{6}s_{7}s_{4}s_{5}s_{3}s_{2}s_{4}s_{3}s_{1} $,
$ s_{6}s_{7}s_{4}s_{5}s_{3}s_{2}s_{4}s_{3}s_{2} $,
$ s_{6}s_{7}s_{4}s_{5}s_{3}s_{4}s_{1}s_{3}s_{1} $,
$ s_{6}s_{7}s_{4}s_{5}s_{3}s_{4}s_{1}s_{3}s_{2} $,
$ s_{6}s_{7}s_{4}s_{5}s_{3}s_{4}s_{3}s_{2}s_{1} $,
$ s_{6}s_{7}s_{4}s_{5}s_{4}s_{1}s_{3}s_{2}s_{4} $,
$ s_{6}s_{7}s_{4}s_{5}s_{4}s_{3}s_{2}s_{4}s_{1} $,
$ s_{6}s_{7}s_{4}s_{5}s_{4}s_{3}s_{2}s_{4}s_{2} $,
$ s_{6}s_{7}s_{4}s_{5}s_{6}s_{2}s_{4}s_{3}s_{2} $,
$ s_{6}s_{7}s_{4}s_{5}s_{6}s_{2}s_{4}s_{5}s_{2} $,
$ s_{6}s_{7}s_{4}s_{5}s_{6}s_{2}s_{4}s_{5}s_{3} $,
$ s_{6}s_{7}s_{4}s_{5}s_{6}s_{3}s_{4}s_{3}s_{2} $,
$ s_{6}s_{7}s_{4}s_{5}s_{6}s_{3}s_{4}s_{5}s_{2} $,
$ s_{6}s_{7}s_{4}s_{5}s_{6}s_{3}s_{4}s_{5}s_{3} $,
$ s_{6}s_{7}s_{4}s_{5}s_{6}s_{5}s_{2}s_{4}s_{2} $,
$ s_{6}s_{7}s_{4}s_{5}s_{6}s_{5}s_{2}s_{4}s_{3} $,
$ s_{6}s_{7}s_{4}s_{5}s_{6}s_{5}s_{3}s_{4}s_{2} $,
$ s_{6}s_{7}s_{4}s_{5}s_{6}s_{5}s_{3}s_{4}s_{3} $,
$ s_{6}s_{7}s_{5}s_{1}s_{3}s_{2}s_{4}s_{3}s_{2} $,
$ s_{6}s_{7}s_{5}s_{1}s_{3}s_{4}s_{1}s_{3}s_{2} $,
$ s_{6}s_{7}s_{5}s_{1}s_{3}s_{4}s_{3}s_{2}s_{1} $,
$ s_{6}s_{7}s_{5}s_{3}s_{2}s_{4}s_{1}s_{3}s_{2} $,
$ s_{6}s_{7}s_{5}s_{3}s_{4}s_{1}s_{3}s_{2}s_{1} $,
$ s_{6}s_{7}s_{5}s_{3}s_{4}s_{1}s_{3}s_{2}s_{4} $,
$ s_{6}s_{7}s_{5}s_{4}s_{1}s_{3}s_{2}s_{4}s_{1} $,
$ s_{6}s_{7}s_{5}s_{4}s_{1}s_{3}s_{2}s_{4}s_{2} $,
$ s_{6}s_{7}s_{5}s_{4}s_{1}s_{3}s_{2}s_{4}s_{3} $,
$ s_{6}s_{7}s_{5}s_{4}s_{1}s_{3}s_{2}s_{4}s_{5} $,
$ s_{6}s_{7}s_{5}s_{4}s_{3}s_{2}s_{4}s_{1}s_{3} $,
$ s_{6}s_{7}s_{5}s_{4}s_{3}s_{2}s_{4}s_{5}s_{1} $,
$ s_{6}s_{7}s_{5}s_{4}s_{3}s_{2}s_{4}s_{5}s_{3} $,
$ s_{6}s_{7}s_{5}s_{4}s_{3}s_{2}s_{4}s_{5}s_{4} $,
$ s_{6}s_{7}s_{5}s_{6}s_{2}s_{4}s_{5}s_{2}s_{4} $,
$ s_{6}s_{7}s_{5}s_{6}s_{2}s_{4}s_{5}s_{3}s_{2} $,
$ s_{6}s_{7}s_{5}s_{6}s_{2}s_{4}s_{5}s_{3}s_{4} $,
$ s_{6}s_{7}s_{5}s_{6}s_{2}s_{4}s_{5}s_{4}s_{2} $,
$ s_{6}s_{7}s_{5}s_{6}s_{2}s_{4}s_{5}s_{4}s_{3} $,
$ s_{6}s_{7}s_{5}s_{6}s_{3}s_{2}s_{4}s_{5}s_{2} $,
$ s_{6}s_{7}s_{5}s_{6}s_{3}s_{2}s_{4}s_{5}s_{3} $,
$ s_{6}s_{7}s_{5}s_{6}s_{3}s_{4}s_{5}s_{2}s_{4} $,
$ s_{6}s_{7}s_{5}s_{6}s_{3}s_{4}s_{5}s_{3}s_{2} $,
$ s_{6}s_{7}s_{5}s_{6}s_{3}s_{4}s_{5}s_{3}s_{4} $,
$ s_{6}s_{7}s_{5}s_{6}s_{3}s_{4}s_{5}s_{4}s_{2} $,
$ s_{6}s_{7}s_{5}s_{6}s_{3}s_{4}s_{5}s_{4}s_{3} $,
$ s_{6}s_{7}s_{5}s_{6}s_{4}s_{5}s_{2}s_{4}s_{2} $,
$ s_{6}s_{7}s_{5}s_{6}s_{4}s_{5}s_{2}s_{4}s_{3} $,
$ s_{6}s_{7}s_{5}s_{6}s_{4}s_{5}s_{3}s_{4}s_{2} $,
$ s_{6}s_{7}s_{5}s_{6}s_{4}s_{5}s_{3}s_{4}s_{3} $,
$ s_{7}s_{1}s_{3}s_{2}s_{4}s_{5}s_{6}s_{1}s_{3} $,
$ s_{7}s_{1}s_{3}s_{2}s_{4}s_{5}s_{6}s_{3}s_{1} $,
$ s_{7}s_{1}s_{3}s_{2}s_{4}s_{5}s_{6}s_{3}s_{4} $,
$ s_{7}s_{1}s_{3}s_{2}s_{4}s_{5}s_{6}s_{4}s_{3} $,
$ s_{7}s_{1}s_{3}s_{4}s_{5}s_{6}s_{1}s_{3}s_{4} $,
$ s_{7}s_{1}s_{3}s_{4}s_{5}s_{6}s_{2}s_{4}s_{3} $,
$ s_{7}s_{1}s_{3}s_{4}s_{5}s_{6}s_{2}s_{4}s_{5} $,
$ s_{7}s_{1}s_{3}s_{4}s_{5}s_{6}s_{3}s_{2}s_{4} $,
$ s_{7}s_{1}s_{3}s_{4}s_{5}s_{6}s_{3}s_{4}s_{1} $,
$ s_{7}s_{1}s_{3}s_{4}s_{5}s_{6}s_{3}s_{4}s_{2} $,
$ s_{7}s_{1}s_{3}s_{4}s_{5}s_{6}s_{3}s_{4}s_{5} $,
$ s_{7}s_{1}s_{3}s_{4}s_{5}s_{6}s_{4}s_{1}s_{3} $,
$ s_{7}s_{1}s_{3}s_{4}s_{5}s_{6}s_{4}s_{3}s_{1} $,
$ s_{7}s_{1}s_{3}s_{4}s_{5}s_{6}s_{4}s_{3}s_{2} $,
$ s_{7}s_{1}s_{3}s_{4}s_{5}s_{6}s_{4}s_{5}s_{2} $,
$ s_{7}s_{1}s_{3}s_{4}s_{5}s_{6}s_{4}s_{5}s_{3} $,
$ s_{7}s_{1}s_{3}s_{4}s_{5}s_{6}s_{5}s_{2}s_{4} $,
$ s_{7}s_{1}s_{3}s_{4}s_{5}s_{6}s_{5}s_{3}s_{4} $,
$ s_{7}s_{1}s_{3}s_{4}s_{5}s_{6}s_{5}s_{4}s_{2} $,
$ s_{7}s_{1}s_{3}s_{4}s_{5}s_{6}s_{5}s_{4}s_{3} $,
$ s_{7}s_{2}s_{4}s_{5}s_{6}s_{1}s_{3}s_{4}s_{1} $,
$ s_{7}s_{2}s_{4}s_{5}s_{6}s_{1}s_{3}s_{4}s_{3} $,
$ s_{7}s_{2}s_{4}s_{5}s_{6}s_{1}s_{3}s_{4}s_{5} $,
$ s_{7}s_{2}s_{4}s_{5}s_{6}s_{2}s_{4}s_{1}s_{3} $,
$ s_{7}s_{2}s_{4}s_{5}s_{6}s_{2}s_{4}s_{3}s_{1} $,
$ s_{7}s_{2}s_{4}s_{5}s_{6}s_{3}s_{2}s_{4}s_{3} $,
$ s_{7}s_{2}s_{4}s_{5}s_{6}s_{3}s_{4}s_{1}s_{3} $,
$ s_{7}s_{2}s_{4}s_{5}s_{6}s_{3}s_{4}s_{5}s_{1} $,
$ s_{7}s_{2}s_{4}s_{5}s_{6}s_{3}s_{4}s_{5}s_{3} $,
$ s_{7}s_{2}s_{4}s_{5}s_{6}s_{4}s_{1}s_{3}s_{2} $,
$ s_{7}s_{2}s_{4}s_{5}s_{6}s_{4}s_{3}s_{2}s_{1} $,
$ s_{7}s_{2}s_{4}s_{5}s_{6}s_{4}s_{3}s_{2}s_{4} $,
$ s_{7}s_{2}s_{4}s_{5}s_{6}s_{4}s_{5}s_{1}s_{3} $,
$ s_{7}s_{2}s_{4}s_{5}s_{6}s_{4}s_{5}s_{2}s_{4} $,
$ s_{7}s_{2}s_{4}s_{5}s_{6}s_{4}s_{5}s_{3}s_{1} $,
$ s_{7}s_{2}s_{4}s_{5}s_{6}s_{4}s_{5}s_{3}s_{2} $,
$ s_{7}s_{2}s_{4}s_{5}s_{6}s_{4}s_{5}s_{3}s_{4} $,
$ s_{7}s_{2}s_{4}s_{5}s_{6}s_{5}s_{1}s_{3}s_{4} $,
$ s_{7}s_{2}s_{4}s_{5}s_{6}s_{5}s_{2}s_{4}s_{3} $,
$ s_{7}s_{2}s_{4}s_{5}s_{6}s_{5}s_{3}s_{4}s_{1} $,
$ s_{7}s_{2}s_{4}s_{5}s_{6}s_{5}s_{4}s_{1}s_{3} $,
$ s_{7}s_{2}s_{4}s_{5}s_{6}s_{5}s_{4}s_{3}s_{1} $,
$ s_{7}s_{3}s_{2}s_{4}s_{5}s_{6}s_{1}s_{3}s_{1} $,
$ s_{7}s_{3}s_{2}s_{4}s_{5}s_{6}s_{1}s_{3}s_{2} $,
$ s_{7}s_{3}s_{2}s_{4}s_{5}s_{6}s_{1}s_{3}s_{4} $,
$ s_{7}s_{3}s_{2}s_{4}s_{5}s_{6}s_{2}s_{4}s_{3} $,
$ s_{7}s_{3}s_{2}s_{4}s_{5}s_{6}s_{3}s_{2}s_{1} $,
$ s_{7}s_{3}s_{2}s_{4}s_{5}s_{6}s_{3}s_{2}s_{4} $,
$ s_{7}s_{3}s_{2}s_{4}s_{5}s_{6}s_{3}s_{4}s_{5} $,
$ s_{7}s_{3}s_{2}s_{4}s_{5}s_{6}s_{4}s_{1}s_{3} $,
$ s_{7}s_{3}s_{2}s_{4}s_{5}s_{6}s_{4}s_{5}s_{3} $,
$ s_{7}s_{3}s_{2}s_{4}s_{5}s_{6}s_{5}s_{3}s_{4} $,
$ s_{7}s_{3}s_{2}s_{4}s_{5}s_{6}s_{5}s_{4}s_{3} $,
$ s_{7}s_{3}s_{4}s_{1}s_{3}s_{2}s_{4}s_{5}s_{6} $,
$ s_{7}s_{3}s_{4}s_{5}s_{6}s_{1}s_{3}s_{4}s_{2} $,
$ s_{7}s_{3}s_{4}s_{5}s_{6}s_{1}s_{3}s_{4}s_{3} $,
$ s_{7}s_{3}s_{4}s_{5}s_{6}s_{2}s_{4}s_{3}s_{1} $,
$ s_{7}s_{3}s_{4}s_{5}s_{6}s_{2}s_{4}s_{5}s_{1} $,
$ s_{7}s_{3}s_{4}s_{5}s_{6}s_{3}s_{2}s_{4}s_{1} $,
$ s_{7}s_{3}s_{4}s_{5}s_{6}s_{3}s_{4}s_{1}s_{3} $,
$ s_{7}s_{3}s_{4}s_{5}s_{6}s_{3}s_{4}s_{5}s_{1} $,
$ s_{7}s_{3}s_{4}s_{5}s_{6}s_{3}s_{4}s_{5}s_{2} $,
$ s_{7}s_{3}s_{4}s_{5}s_{6}s_{4}s_{1}s_{3}s_{1} $,
$ s_{7}s_{3}s_{4}s_{5}s_{6}s_{4}s_{1}s_{3}s_{2} $,
$ s_{7}s_{3}s_{4}s_{5}s_{6}s_{4}s_{5}s_{2}s_{1} $,
$ s_{7}s_{3}s_{4}s_{5}s_{6}s_{4}s_{5}s_{2}s_{4} $,
$ s_{7}s_{3}s_{4}s_{5}s_{6}s_{4}s_{5}s_{3}s_{1} $,
$ s_{7}s_{3}s_{4}s_{5}s_{6}s_{4}s_{5}s_{3}s_{4} $,
$ s_{7}s_{3}s_{4}s_{5}s_{6}s_{5}s_{2}s_{4}s_{1} $,
$ s_{7}s_{3}s_{4}s_{5}s_{6}s_{5}s_{3}s_{2}s_{4} $,
$ s_{7}s_{3}s_{4}s_{5}s_{6}s_{5}s_{3}s_{4}s_{1} $,
$ s_{7}s_{3}s_{4}s_{5}s_{6}s_{5}s_{3}s_{4}s_{3} $,
$ s_{7}s_{3}s_{4}s_{5}s_{6}s_{5}s_{4}s_{2}s_{1} $,
$ s_{7}s_{3}s_{4}s_{5}s_{6}s_{5}s_{4}s_{3}s_{1} $,
$ s_{7}s_{4}s_{1}s_{3}s_{2}s_{4}s_{5}s_{6}s_{1} $,
$ s_{7}s_{4}s_{1}s_{3}s_{2}s_{4}s_{5}s_{6}s_{4} $,
$ s_{7}s_{4}s_{3}s_{2}s_{4}s_{5}s_{6}s_{1}s_{3} $,
$ s_{7}s_{4}s_{3}s_{2}s_{4}s_{5}s_{6}s_{2}s_{1} $,
$ s_{7}s_{4}s_{3}s_{2}s_{4}s_{5}s_{6}s_{3}s_{1} $,
$ s_{7}s_{4}s_{3}s_{2}s_{4}s_{5}s_{6}s_{3}s_{4} $,
$ s_{7}s_{4}s_{3}s_{2}s_{4}s_{5}s_{6}s_{4}s_{1} $,
$ s_{7}s_{4}s_{3}s_{2}s_{4}s_{5}s_{6}s_{4}s_{2} $,
$ s_{7}s_{4}s_{3}s_{2}s_{4}s_{5}s_{6}s_{4}s_{3} $,
$ s_{7}s_{4}s_{3}s_{2}s_{4}s_{5}s_{6}s_{4}s_{5} $,
$ s_{7}s_{4}s_{3}s_{2}s_{4}s_{5}s_{6}s_{5}s_{4} $,
$ s_{7}s_{4}s_{5}s_{6}s_{1}s_{3}s_{2}s_{4}s_{3} $,
$ s_{7}s_{4}s_{5}s_{6}s_{1}s_{3}s_{4}s_{2}s_{1} $,
$ s_{7}s_{4}s_{5}s_{6}s_{1}s_{3}s_{4}s_{3}s_{1} $,
$ s_{7}s_{4}s_{5}s_{6}s_{1}s_{3}s_{4}s_{3}s_{2} $,
$ s_{7}s_{4}s_{5}s_{6}s_{1}s_{3}s_{4}s_{5}s_{2} $,
$ s_{7}s_{4}s_{5}s_{6}s_{1}s_{3}s_{4}s_{5}s_{3} $,
$ s_{7}s_{4}s_{5}s_{6}s_{2}s_{4}s_{1}s_{3}s_{2} $,
$ s_{7}s_{4}s_{5}s_{6}s_{2}s_{4}s_{3}s_{2}s_{1} $,
$ s_{7}s_{4}s_{5}s_{6}s_{2}s_{4}s_{5}s_{1}s_{3} $,
$ s_{7}s_{4}s_{5}s_{6}s_{2}s_{4}s_{5}s_{3}s_{1} $,
$ s_{7}s_{4}s_{5}s_{6}s_{3}s_{2}s_{4}s_{3}s_{1} $,
$ s_{7}s_{4}s_{5}s_{6}s_{3}s_{4}s_{1}s_{3}s_{1} $,
$ s_{7}s_{4}s_{5}s_{6}s_{3}s_{4}s_{1}s_{3}s_{2} $,
$ s_{7}s_{4}s_{5}s_{6}s_{3}s_{4}s_{5}s_{2}s_{1} $,
$ s_{7}s_{4}s_{5}s_{6}s_{3}s_{4}s_{5}s_{3}s_{1} $,
$ s_{7}s_{4}s_{5}s_{6}s_{3}s_{4}s_{5}s_{3}s_{2} $,
$ s_{7}s_{4}s_{5}s_{6}s_{4}s_{1}s_{3}s_{2}s_{4} $,
$ s_{7}s_{4}s_{5}s_{6}s_{4}s_{3}s_{2}s_{4}s_{1} $,
$ s_{7}s_{4}s_{5}s_{6}s_{4}s_{3}s_{2}s_{4}s_{2} $,
$ s_{7}s_{4}s_{5}s_{6}s_{4}s_{3}s_{2}s_{4}s_{3} $,
$ s_{7}s_{4}s_{5}s_{6}s_{4}s_{3}s_{2}s_{4}s_{5} $,
$ s_{7}s_{4}s_{5}s_{6}s_{4}s_{5}s_{2}s_{4}s_{2} $,
$ s_{7}s_{4}s_{5}s_{6}s_{4}s_{5}s_{2}s_{4}s_{3} $,
$ s_{7}s_{4}s_{5}s_{6}s_{4}s_{5}s_{3}s_{4}s_{2} $,
$ s_{7}s_{4}s_{5}s_{6}s_{4}s_{5}s_{3}s_{4}s_{3} $,
$ s_{7}s_{4}s_{5}s_{6}s_{5}s_{1}s_{3}s_{4}s_{2} $,
$ s_{7}s_{4}s_{5}s_{6}s_{5}s_{1}s_{3}s_{4}s_{3} $,
$ s_{7}s_{4}s_{5}s_{6}s_{5}s_{2}s_{4}s_{1}s_{3} $,
$ s_{7}s_{4}s_{5}s_{6}s_{5}s_{2}s_{4}s_{3}s_{1} $,
$ s_{7}s_{4}s_{5}s_{6}s_{5}s_{3}s_{2}s_{4}s_{3} $,
$ s_{7}s_{4}s_{5}s_{6}s_{5}s_{3}s_{4}s_{2}s_{1} $,
$ s_{7}s_{4}s_{5}s_{6}s_{5}s_{3}s_{4}s_{3}s_{1} $,
$ s_{7}s_{4}s_{5}s_{6}s_{5}s_{4}s_{3}s_{2}s_{4} $,
$ s_{7}s_{5}s_{4}s_{1}s_{3}s_{2}s_{4}s_{5}s_{6} $,
$ s_{7}s_{5}s_{4}s_{3}s_{2}s_{4}s_{5}s_{6}s_{1} $,
$ s_{7}s_{5}s_{4}s_{3}s_{2}s_{4}s_{5}s_{6}s_{2} $,
$ s_{7}s_{5}s_{4}s_{3}s_{2}s_{4}s_{5}s_{6}s_{4} $,
$ s_{7}s_{5}s_{4}s_{3}s_{2}s_{4}s_{5}s_{6}s_{5} $,
$ s_{7}s_{5}s_{6}s_{1}s_{3}s_{2}s_{4}s_{3}s_{2} $,
$ s_{7}s_{5}s_{6}s_{1}s_{3}s_{2}s_{4}s_{5}s_{2} $,
$ s_{7}s_{5}s_{6}s_{1}s_{3}s_{2}s_{4}s_{5}s_{3} $,
$ s_{7}s_{5}s_{6}s_{1}s_{3}s_{4}s_{1}s_{3}s_{2} $,
$ s_{7}s_{5}s_{6}s_{1}s_{3}s_{4}s_{3}s_{2}s_{1} $,
$ s_{7}s_{5}s_{6}s_{1}s_{3}s_{4}s_{5}s_{2}s_{4} $,
$ s_{7}s_{5}s_{6}s_{1}s_{3}s_{4}s_{5}s_{3}s_{2} $,
$ s_{7}s_{5}s_{6}s_{1}s_{3}s_{4}s_{5}s_{3}s_{4} $,
$ s_{7}s_{5}s_{6}s_{1}s_{3}s_{4}s_{5}s_{4}s_{2} $,
$ s_{7}s_{5}s_{6}s_{1}s_{3}s_{4}s_{5}s_{4}s_{3} $,
$ s_{7}s_{5}s_{6}s_{2}s_{4}s_{5}s_{1}s_{3}s_{2} $,
$ s_{7}s_{5}s_{6}s_{2}s_{4}s_{5}s_{1}s_{3}s_{4} $,
$ s_{7}s_{5}s_{6}s_{2}s_{4}s_{5}s_{2}s_{4}s_{2} $,
$ s_{7}s_{5}s_{6}s_{2}s_{4}s_{5}s_{3}s_{2}s_{1} $,
$ s_{7}s_{5}s_{6}s_{2}s_{4}s_{5}s_{3}s_{4}s_{1} $,
$ s_{7}s_{5}s_{6}s_{2}s_{4}s_{5}s_{4}s_{1}s_{3} $,
$ s_{7}s_{5}s_{6}s_{2}s_{4}s_{5}s_{4}s_{3}s_{1} $,
$ s_{7}s_{5}s_{6}s_{3}s_{2}s_{4}s_{1}s_{3}s_{2} $,
$ s_{7}s_{5}s_{6}s_{3}s_{2}s_{4}s_{5}s_{2}s_{1} $,
$ s_{7}s_{5}s_{6}s_{3}s_{2}s_{4}s_{5}s_{3}s_{1} $,
$ s_{7}s_{5}s_{6}s_{3}s_{2}s_{4}s_{5}s_{4}s_{3} $,
$ s_{7}s_{5}s_{6}s_{3}s_{4}s_{1}s_{3}s_{2}s_{1} $,
$ s_{7}s_{5}s_{6}s_{3}s_{4}s_{1}s_{3}s_{2}s_{4} $,
$ s_{7}s_{5}s_{6}s_{3}s_{4}s_{5}s_{2}s_{4}s_{1} $,
$ s_{7}s_{5}s_{6}s_{3}s_{4}s_{5}s_{3}s_{2}s_{1} $,
$ s_{7}s_{5}s_{6}s_{3}s_{4}s_{5}s_{3}s_{2}s_{4} $,
$ s_{7}s_{5}s_{6}s_{3}s_{4}s_{5}s_{3}s_{4}s_{1} $,
$ s_{7}s_{5}s_{6}s_{3}s_{4}s_{5}s_{3}s_{4}s_{2} $,
$ s_{7}s_{5}s_{6}s_{3}s_{4}s_{5}s_{4}s_{2}s_{1} $,
$ s_{7}s_{5}s_{6}s_{3}s_{4}s_{5}s_{4}s_{3}s_{1} $,
$ s_{7}s_{5}s_{6}s_{4}s_{1}s_{3}s_{2}s_{4}s_{1} $,
$ s_{7}s_{5}s_{6}s_{4}s_{1}s_{3}s_{2}s_{4}s_{2} $,
$ s_{7}s_{5}s_{6}s_{4}s_{1}s_{3}s_{2}s_{4}s_{3} $,
$ s_{7}s_{5}s_{6}s_{4}s_{3}s_{2}s_{4}s_{1}s_{3} $,
$ s_{7}s_{5}s_{6}s_{4}s_{3}s_{2}s_{4}s_{5}s_{2} $,
$ s_{7}s_{5}s_{6}s_{4}s_{3}s_{2}s_{4}s_{5}s_{3} $,
$ s_{7}s_{5}s_{6}s_{4}s_{3}s_{2}s_{4}s_{5}s_{4} $,
$ s_{7}s_{5}s_{6}s_{4}s_{5}s_{1}s_{3}s_{4}s_{2} $,
$ s_{7}s_{5}s_{6}s_{4}s_{5}s_{1}s_{3}s_{4}s_{3} $,
$ s_{7}s_{5}s_{6}s_{4}s_{5}s_{2}s_{4}s_{1}s_{3} $,
$ s_{7}s_{5}s_{6}s_{4}s_{5}s_{2}s_{4}s_{3}s_{1} $,
$ s_{7}s_{5}s_{6}s_{4}s_{5}s_{3}s_{2}s_{4}s_{3} $,
$ s_{7}s_{5}s_{6}s_{4}s_{5}s_{3}s_{4}s_{2}s_{1} $,
$ s_{7}s_{5}s_{6}s_{4}s_{5}s_{3}s_{4}s_{3}s_{1} $,
$ s_{7}s_{5}s_{6}s_{4}s_{5}s_{4}s_{3}s_{2}s_{4} $,
$ s_{7}s_{5}s_{6}s_{5}s_{4}s_{3}s_{2}s_{4}s_{5} $,
$ s_{7}s_{6}s_{1}s_{3}s_{2}s_{4}s_{5}s_{1}s_{3} $,
$ s_{7}s_{6}s_{1}s_{3}s_{2}s_{4}s_{5}s_{2}s_{4} $,
$ s_{7}s_{6}s_{1}s_{3}s_{2}s_{4}s_{5}s_{3}s_{1} $,
$ s_{7}s_{6}s_{1}s_{3}s_{2}s_{4}s_{5}s_{4}s_{2} $,
$ s_{7}s_{6}s_{1}s_{3}s_{4}s_{5}s_{1}s_{3}s_{4} $,
$ s_{7}s_{6}s_{1}s_{3}s_{4}s_{5}s_{2}s_{4}s_{2} $,
$ s_{7}s_{6}s_{1}s_{3}s_{4}s_{5}s_{2}s_{4}s_{3} $,
$ s_{7}s_{6}s_{1}s_{3}s_{4}s_{5}s_{3}s_{2}s_{4} $,
$ s_{7}s_{6}s_{1}s_{3}s_{4}s_{5}s_{3}s_{4}s_{1} $,
$ s_{7}s_{6}s_{1}s_{3}s_{4}s_{5}s_{4}s_{1}s_{3} $,
$ s_{7}s_{6}s_{1}s_{3}s_{4}s_{5}s_{4}s_{3}s_{1} $,
$ s_{7}s_{6}s_{2}s_{4}s_{5}s_{1}s_{3}s_{4}s_{1} $,
$ s_{7}s_{6}s_{2}s_{4}s_{5}s_{3}s_{2}s_{4}s_{3} $,
$ s_{7}s_{6}s_{2}s_{4}s_{5}s_{3}s_{4}s_{1}s_{3} $,
$ s_{7}s_{6}s_{2}s_{4}s_{5}s_{3}s_{4}s_{3}s_{1} $,
$ s_{7}s_{6}s_{2}s_{4}s_{5}s_{4}s_{3}s_{2}s_{4} $,
$ s_{7}s_{6}s_{3}s_{2}s_{4}s_{5}s_{1}s_{3}s_{1} $,
$ s_{7}s_{6}s_{3}s_{2}s_{4}s_{5}s_{1}s_{3}s_{2} $,
$ s_{7}s_{6}s_{3}s_{2}s_{4}s_{5}s_{1}s_{3}s_{4} $,
$ s_{7}s_{6}s_{3}s_{2}s_{4}s_{5}s_{2}s_{4}s_{1} $,
$ s_{7}s_{6}s_{3}s_{2}s_{4}s_{5}s_{2}s_{4}s_{2} $,
$ s_{7}s_{6}s_{3}s_{2}s_{4}s_{5}s_{3}s_{2}s_{1} $,
$ s_{7}s_{6}s_{3}s_{2}s_{4}s_{5}s_{3}s_{2}s_{4} $,
$ s_{7}s_{6}s_{3}s_{2}s_{4}s_{5}s_{3}s_{4}s_{1} $,
$ s_{7}s_{6}s_{3}s_{2}s_{4}s_{5}s_{4}s_{1}s_{3} $,
$ s_{7}s_{6}s_{3}s_{2}s_{4}s_{5}s_{4}s_{2}s_{1} $,
$ s_{7}s_{6}s_{3}s_{2}s_{4}s_{5}s_{4}s_{3}s_{1} $,
$ s_{7}s_{6}s_{3}s_{4}s_{1}s_{3}s_{2}s_{4}s_{5} $,
$ s_{7}s_{6}s_{3}s_{4}s_{5}s_{1}s_{3}s_{4}s_{2} $,
$ s_{7}s_{6}s_{3}s_{4}s_{5}s_{1}s_{3}s_{4}s_{3} $,
$ s_{7}s_{6}s_{3}s_{4}s_{5}s_{2}s_{4}s_{2}s_{1} $,
$ s_{7}s_{6}s_{3}s_{4}s_{5}s_{2}s_{4}s_{3}s_{1} $,
$ s_{7}s_{6}s_{3}s_{4}s_{5}s_{3}s_{2}s_{4}s_{1} $,
$ s_{7}s_{6}s_{3}s_{4}s_{5}s_{3}s_{4}s_{1}s_{3} $,
$ s_{7}s_{6}s_{3}s_{4}s_{5}s_{3}s_{4}s_{2}s_{1} $,
$ s_{7}s_{6}s_{3}s_{4}s_{5}s_{4}s_{1}s_{3}s_{1} $,
$ s_{7}s_{6}s_{3}s_{4}s_{5}s_{4}s_{1}s_{3}s_{2} $,
$ s_{7}s_{6}s_{3}s_{4}s_{5}s_{4}s_{3}s_{2}s_{1} $,
$ s_{7}s_{6}s_{3}s_{4}s_{5}s_{4}s_{3}s_{2}s_{4} $,
$ s_{7}s_{6}s_{4}s_{1}s_{3}s_{2}s_{4}s_{5}s_{1} $,
$ s_{7}s_{6}s_{4}s_{1}s_{3}s_{2}s_{4}s_{5}s_{4} $,
$ s_{7}s_{6}s_{4}s_{3}s_{2}s_{4}s_{5}s_{1}s_{3} $,
$ s_{7}s_{6}s_{4}s_{3}s_{2}s_{4}s_{5}s_{2}s_{1} $,
$ s_{7}s_{6}s_{4}s_{3}s_{2}s_{4}s_{5}s_{3}s_{1} $,
$ s_{7}s_{6}s_{4}s_{3}s_{2}s_{4}s_{5}s_{4}s_{1} $,
$ s_{7}s_{6}s_{4}s_{5}s_{1}s_{3}s_{2}s_{4}s_{3} $,
$ s_{7}s_{6}s_{4}s_{5}s_{1}s_{3}s_{4}s_{2}s_{1} $,
$ s_{7}s_{6}s_{4}s_{5}s_{1}s_{3}s_{4}s_{3}s_{1} $,
$ s_{7}s_{6}s_{4}s_{5}s_{3}s_{2}s_{4}s_{3}s_{1} $,
$ s_{7}s_{6}s_{4}s_{5}s_{3}s_{4}s_{1}s_{3}s_{1} $,
$ s_{7}s_{6}s_{4}s_{5}s_{3}s_{4}s_{1}s_{3}s_{2} $,
$ s_{7}s_{6}s_{4}s_{5}s_{3}s_{4}s_{3}s_{2}s_{1} $,
$ s_{7}s_{6}s_{4}s_{5}s_{4}s_{1}s_{3}s_{2}s_{4} $,
$ s_{7}s_{6}s_{4}s_{5}s_{4}s_{3}s_{2}s_{4}s_{1} $,
$ s_{7}s_{6}s_{4}s_{5}s_{4}s_{3}s_{2}s_{4}s_{2} $,
$ s_{7}s_{6}s_{5}s_{1}s_{3}s_{2}s_{4}s_{3}s_{2} $,
$ s_{7}s_{6}s_{5}s_{1}s_{3}s_{4}s_{1}s_{3}s_{2} $,
$ s_{7}s_{6}s_{5}s_{1}s_{3}s_{4}s_{3}s_{2}s_{1} $,
$ s_{7}s_{6}s_{5}s_{3}s_{2}s_{4}s_{1}s_{3}s_{2} $,
$ s_{7}s_{6}s_{5}s_{3}s_{4}s_{1}s_{3}s_{2}s_{1} $,
$ s_{7}s_{6}s_{5}s_{3}s_{4}s_{1}s_{3}s_{2}s_{4} $,
$ s_{7}s_{6}s_{5}s_{4}s_{1}s_{3}s_{2}s_{4}s_{1} $,
$ s_{7}s_{6}s_{5}s_{4}s_{1}s_{3}s_{2}s_{4}s_{2} $,
$ s_{7}s_{6}s_{5}s_{4}s_{1}s_{3}s_{2}s_{4}s_{3} $,
$ s_{7}s_{6}s_{5}s_{4}s_{1}s_{3}s_{2}s_{4}s_{5} $,
$ s_{7}s_{6}s_{5}s_{4}s_{3}s_{2}s_{4}s_{1}s_{3} $,
$ s_{7}s_{6}s_{5}s_{4}s_{3}s_{2}s_{4}s_{5}s_{1} $,
$ s_{7}s_{6}s_{5}s_{4}s_{3}s_{2}s_{4}s_{5}s_{3} $,
$ s_{7}s_{6}s_{5}s_{4}s_{3}s_{2}s_{4}s_{5}s_{4} $,
$ s_{7}s_{6}s_{5}s_{4}s_{3}s_{2}s_{4}s_{5}s_{6} $\}.
\end{description}
\end{Prop}

\noindent\textbf{Remark.}
Here we use the here we take advantage of the work of Duan and Zhao \cite{duan2014schubert}, for one choice of $w$.

\smallbreak
\noindent\textbf{Remark.}
It took roughly 2 hours to get the result on author's personal computer.

}

\end{document}